\documentclass[aop,MSNbibl,seceqn,dvips]{arximspdf}
\usepackage{mathbh}
\usepackage{graphicx}

\doi{10.1214/13-AOP851} 
\volume{42}
\issue{3}
\pubyear{2014}
\firstpage{1054}
\lastpage{1120}

\makeatletter

\newcommand{\rright}{\right}
\newcommand{\lleft}{\left}
\newcommand{\rrVert}{\Vert}
\newcommand{\rrvert}{\vert}
\newcommand{\llVert}{\Vert}
\newcommand{\llvert}{\vert}
\DeclareFontFamily{U}{mathx}{\hyphenchar\font45}
\DeclareFontShape{U}{mathx}{m}{n}{
<5> <6> <7> <8> <9> <10>
<10.95> <12> <14.4> <17.28> <20.74> <24.88>
mathx10
}{}
\DeclareSymbolFont{mathx}{U}{mathx}{m}{n}
\DeclareFontSubstitution{U}{mathx}{m}{n}
\DeclareMathAccent{\widecheck}{0}{mathx}{"71}
\DeclareMathAccent{\wideparen}{0}{mathx}{"75}

\newtheorem{thmm}{Theorem}[section]

\newtheorem{hypo}{Hypothesis}
\newtheorem{lem}{Lemma}
\newtheorem{cor}{Corollary}
\newtheorem{prop}{Proposition}[section]
\newproclaim{dfn}{Definition}[section]
\newproclaim{rem}{Remark}
\newproclaim{remfr}{Remarque}

\newcommand{\cyl}{\operatorname{cyl}}
\newcommand{\card}{\operatorname{card}}
\newcommand{\di}{\operatorname{div}}
\newcommand{\disc}{\operatorname{disc}}
\newcommand{\hyp}{\operatorname{hyp}}
\newcommand{\capa}{\operatorname{capacity}}

\newcommand{\eps}{\varepsilon}

\def\PP{\mathbb{P}}
\def\RR{\mathbb{R}}
\def\EE{\mathbb{E}}
\def\NN{\mathbb{N}}
\def\ZZ{\mathbb{Z}}
\def\SS{\mathbb{S}}
\def\E{\mathcal{E}}
\def\D{\mathcal{D}}
\def\H{\mathcal{H}}
\def\N{\mathcal{N}}
\def\I{\mathcal{I}}
\def\G{\Gamma}
\def\O{\Omega}
\def\L{\mathcal{L}}
\def\P{\mathcal{P}}
\def\B{\mathcal{B}}
\def\C{\mathcal{C}}
\def\A{\mathcal{A}}
\def\V{\mathcal{V}}
\def\M{\mathcal{M}}
\def\S{\mathcal{S}}
\def\F{\mathcal{F}}
\def\de{\delta}
\def\be{\mathbf{e}}
\def\bbf{\mathbf{f}}
\def\vv{\vec{v}}
\makeatother

\begin{document}
\begin{frontmatter}

\title{Maximal stream and minimal cutset for first passage percolation through a domain of $\RR$\lowercase{$^d$}}
\pdftitle{Maximal stream and minimal cutset for first passage percolation through a domain of Rd}
\runtitle{Maximal stream and minimal cutset}

\begin{aug}
\author[a]{\fnms{Rapha\"el} \snm{Cerf}\ead[label=e1]{rcerf@math.u-psud.fr}\ead[label=u1,url]{http://www.math.u-psud.fr/\textasciitilde cerf}}
\and
\author[b]{\fnms{Marie} \snm{Th\'eret}\corref{}\ead[label=e2]{marie.theret@univ-paris-diderot.fr}\ead[label=u2,url]{http://www.proba.jussieu.fr/\textasciitilde theret}}
\runauthor{R. Cerf and M. Th\'eret}

\affiliation{IUF and Universit\'e Paris Sud, and Universit\'e Paris
Diderot}
\address[a]{Laboratoire de Math\'ematiques, b\^atiment 425\\
IUF and Universit\'e Paris Sud \\
91405 Orsay Cedex\\
France\\
\printead{e1}\\
\printead{u1}}

\address[b]{Universit\'e Paris Diderot\\
LPMA, Site Chevaleret, case 7012\\
75205 Paris Cedex 13\\
France\\
\printead{e2}\\
\printead{u2}}
\end{aug}

\received{\smonth{4} \syear{2012}}
\revised{\smonth{2} \syear{2013}}

%
\begin{abstract}
We consider the standard first passage percolation model in the
rescaled graph $\ZZ^d/n$ for $d\geq2$ and a domain $\O$ of boundary
$\G
$ in $\RR^d$. Let $\G^1$ and $\G^2$ be two disjoint open subsets of
$\G
$, representing the parts of $\G$ through which some water can enter
and escape from $\O$. A law of large numbers for the maximal flow from
$\G^1$ to $\G^2$ in $\O$ is already known. In this paper we investigate
the asymptotic behavior of a maximal stream and a minimal cutset. A
maximal stream is a vector measure $\vec\mu_n^{\max}$ that
describes how
the maximal amount of fluid can cross $\O$. Under conditions on the
regularity of the domain and on the law of the capacities of the edges,
we prove that the sequence $(\vec\mu_n^{\max})_{n\geq1}$ converges a.s.
to the set of the solutions of a continuous deterministic problem of
maximal stream in an anisotropic network. A~minimal cutset can been
seen as the boundary of a set $E_n^{\min}$ that separates $\G^1$ from
$\G^2$ in $\O$ and whose random capacity is minimal. Under the same
conditions, we prove that the sequence $(E_n^{\min})_{n\geq1}$
converges toward the set of the solutions of a continuous deterministic
problem of minimal cutset. We deduce from this a continuous
deterministic max-flow min-cut theorem and a new proof of the law of
large numbers for the maximal flow. This proof is more natural than the
existing one, since it relies on the study of maximal streams and
minimal cutsets, which are the pertinent objects to look at.
\end{abstract}

%

\pdfkeywords{60K35, 49K20, 35Q35, 82B20, First passage percolation, continuous and discrete max-flow min-cut theorem, maximal stream, maximal flow}

\end{frontmatter}

\section{First definitions and main result}
\label{secdef}

We recall first the definitions of the random discrete model and of the
discrete objects. The continuous counterparts of the discrete objects
are briefly presented in Section~\ref{secdef2} and the main results are
presented in Section~\ref{secdef3}.

\subsection{Discrete streams, cutsets and flows}
\label{secdef1}

We use many notation introduced in \cite{Kesten:StFlour} and
\cite{Kesten:flows}. Let $d\geq2$. We consider the graph $(\mathbb{Z}^{d}_n,
\mathbb E ^{d}_n)$ having for vertices $\mathbb Z ^{d}_n = \ZZ^d/n$ and
for edges
$\mathbb E ^{d}_n$,\vspace*{1pt} the set of pairs of nearest neighbors for the standard
$L^{1}$ norm. With each edge $e$ in $\mathbb{E}^{d}_n$ we associate a random
variable $t(e)$ with values in $\mathbb{R}^{+}$. We suppose that the family
$(t(e), e \in\mathbb{E}^{d}_n)$ is independent and identically distributed,
with a common law $\Lambda$: this is the standard model of
first passage percolation on the graph $(\mathbb{Z}^d_n,
\mathbb{E}^d_n)$. We interpret $t(e)$ as the capacity of the edge $e$; it
means that $t(e)$ is the maximal amount of fluid that can go through the
edge $e$ per unit of time.

\begin{figure}

\includegraphics{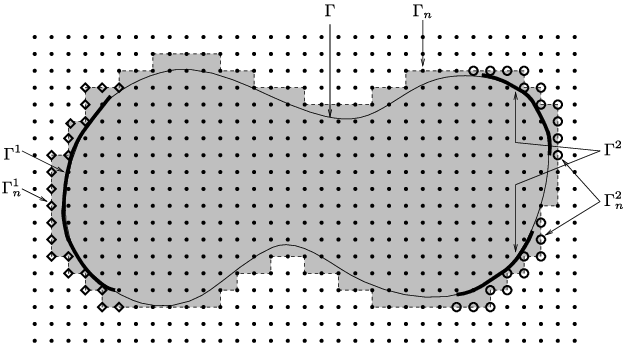}

\caption{Domain $\O$.}
\label{figdomaine}
\end{figure}

We consider an open bounded connected subset $\O$ of $\RR^d$ such that
the boundary $\G= \partial\O$ of $\O$ is piecewise
of class $\C^1$. It means that $\G$ is included in the union of a
finite number of
hypersurfaces of class $\C^1$, that is, in the union of a finite
number of
$\C^1$ submanifolds of~$\RR^d$ of codimension~$1$. Let $\G^1$, $\G^2$
be two disjoint subsets of $\G$ that are open in $\G$. We want to study
the maximal streams from $\G^1$ to $\G^2$ through $\O$ for the
capacities $(t(e), e\in\EE^d_n)$. We consider a discrete version
$(\O_n, \G_n, \G^1_n, \G^2_n)$ of $(\O, \G, \G^1,\G^2)$ defined by
\[
\cases{ %
\O_n = \bigl\{ x\in
\ZZ^d_n | d_{\infty}(x,\O) <1/n \bigr\},
\vspace*{2pt}\cr
\G_n = \bigl\{ x\in \O_n | \exists y \notin
\O_n, [ x,y ] \in\EE^d_n \bigr\},
\vspace*{2pt}\cr
\G^i_n = \bigl\{ x\in\G_n |
d_\infty \bigl(x, \G^i \bigr) <1/n, d_\infty
\bigl(x, \G^{3-i} \bigr) \geq1/n \bigr\},&\quad$\mbox{ for }i=1,2,$}
\]
where $d_{\infty}$ is the $L^{\infty}$-distance, and the segment $[
x,y]$ is the edge of endpoints $x$ and $y$; see Figure~\ref{figdomaine}. We denote by $\Pi_n$ the set of the edges with both
endpoints in~$\O_n$.

We shall study streams and flows from $\G^1_n$ to $\G^2_n$ and cutsets
between $\G^1_n$ and $\G^2_n$ in $\O_n$.
Let us define first the admissible streams from $F_1$ to $F_2$ in $C$,
for $C$ a bounded connected subset of $ \mathbb{R}^d$ and $F_1,F_2$
disjoint sets of vertices of $\ZZ^d_n$ included in $C$. We will say
that an edge
$e=[x,y]$ is included in a subset $A$ of $\mathbb{R}^{d}$, which
we denote by $e\subset A$, if the closed segment joining $x$ to $y$ is
included in $A$. Let $e=[a,b]$ be an edge of $\EE^d_n$ with endpoints
$a$ and $b$. We denote by $\langle a,b \rangle$ the oriented edge
starting at $a$ and ending at $b$. We fix next an orientation for each
edge of $\EE^d_n$. Let $(\vec\bbf_1,\ldots,\vec\bbf_n)$ be the canonical
basis of $\RR^d$. We denote by $\EE_n^{d,i}$ the set of the edges
parallel to $\vec\bbf_i$. For $e=[a,b] \in\EE^{d,i}_n$, we define
\[
\vec{e}= \vec\bbf_i \quad\mbox{and}\quad \be= \cases{ %
\langle a,b \rangle,&\quad $\mbox{if } \overrightarrow {ab} \cdot\vec
\bbf_i = +1/n,$
\vspace*{2pt}\cr
\langle b,a \rangle,&\quad $\mbox{if } \overrightarrow{ab} \cdot\vec\bbf_i
= -1/n,$ }
\]
where $\cdot$ is the scalar product on $\RR^d$ and $\overrightarrow
{ab}$ the vector of origin $a$ and endpoint~$b$. We define the set $\S
_n(F_1,F_2,C)$ of admissible ``stream functions'' as the set of
functions $f_n \dvtx \EE^d_n \rightarrow\RR$ such that:
\begin{longlist}[(iii)]
\item[(i)] \emph{the stream is inside $C$:} for each edge $e\not
\subset C$ we have $f_n(e) = 0$;
\item[(ii)] \emph{capacity constraint:} for each edge $e\in\EE^d_n$
we have
\[
\bigl|f_n(e) \bigr| \leq t(e);
\]
\item[(iii)] \emph{conservation law:} for each vertex $v \in\ZZ^d_n
\setminus(F_1 \cup F_2)$ we have
\[
\sum_{e\in\EE^d_n: \be=\langle v,\cdot \rangle} f_n(e) = \sum
_{e\in\EE^d_n: \be=\langle\cdot,v \rangle} f_n(e),
\]
\end{longlist}
where the notation $\be= \langle v,\cdot \rangle$ (resp., $\be= \langle
\cdot,v \rangle$) means that there exists $y \in\mathbb{Z}_n^d$ such that
$\be= \langle v,y \rangle$ (resp., $\be= \langle y,v \rangle$). A
function $f_n \in\S_n(F_1,F_2,C)$ is a description of a possible
stream in $C$: $|f_n(e)|$ is the amount of water that crosses $e$ per
second, and this water goes through $e$ in the direction of $f_n(e)
\be$ (thus in the direction of $\be$ is $f_n(e) >0$ and in the
direction of $-\be$ if $f_n(e)<0$). Condition (i) means that the
water does not move outside $C$; condition (ii) means that the amount
of water that can cross $e$ per second cannot exceed $t(e)$; condition
(iii) means that there is no loss of fluid in the graph. To each
stream function $f_n$ from $F_1$ to $F_2$ in $C$, we associate the
corresponding flow
\[
\operatorname{flow}^{\mathrm{disc}}_n(f_n) = \sum
_{e\subset C:
e=[a,b], a\in F_1,b\notin F_1} f_n(e) ( \mathbh{1}_{\{ \be= \langle a,b\rangle\}}
- \mathbh {1}_{\{ \be= \langle b,a \rangle\}} ).
\]
This is the amount of fluid (positive or negative) that crosses $C$
from $F_1$ to $F_2$ according to $f_n$. We define the maximal flow
$\phi
_n (F_1,F_2,C)$ from $F_1$ to $F_2$ in $C$ by
\[
\phi_n (F_1,F_2,C) = \sup \bigl\{
\operatorname{flow}^{\mathrm
{disc}}_n(f_n) |
f_n \in\S_n(F_1,F_2,C) \bigr\}
.
\]
If $D$ is a connected set of vertices of $\ZZ^d_n$ that contains two
disjoint subsets $F_1,F_2$ of $\ZZ^d_n$, we define
\[
\widehat D = D + \frac{1}{2n} [-1, 1]^d \subset
\RR^d.
\]
We define
\[
\S_n(F_1,F_2,D) = \S_n
(F_1,F_2,\widehat D) \quad\mbox{and}\quad \phi_n
(F_1,F_2,D) = \phi_n (F_1,F_2,
\widehat D).
\]

\begin{figure}

\includegraphics{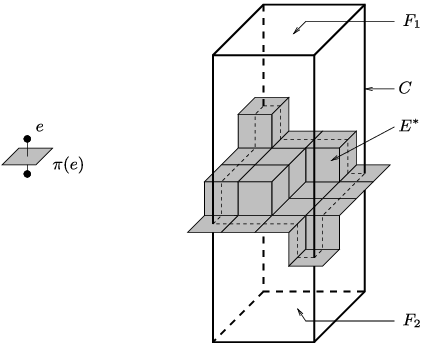}

\caption{Set of plaquettes $E^*$ corresponding to a $(F_1,F_2)$-cutset
$E$ in $C$.}
\label{figplaquette}
\end{figure}

The maximal flow $\phi_n (F_1, F_2, C)$ can be expressed differently
thanks to the (discrete) max-flow min-cut theorem; see \cite{Bollobas}.
We need some
definitions to state this result. A path on the graph $\mathbb
{Z}_n^{d}$ from the vertex $v_{0}$ to the vertex $v_{m}$ is a sequence
$(v_{0}, e_{1}, v_{1},\ldots, e_{m}, v_{m})$ of vertices $v_{0},\ldots,
v_{m}$ alternating with edges $e_{1},\ldots, e_{m}$ such that $v_{i-1}$
and $v_{i}$ are neighbors in the graph, joined by the edge $e_{i}$, for
$i$ in $\{1,\ldots, m\}$. A set $E$ of edges of $\EE^d_n$ included in $C$
is said to cut $F_1$ from $F_2$ in
$C$ if there is no path from $F_1$ to $F_2$ made of edges included in
$C$ that do not belong to~$E$. We call $E$ an $(F_1,F_2)$-cutset in $C$
if $E$ cuts $F_1$ from $F_2$ in $C$
and if no proper subset of $E$ does. With each set of edges $E \subset
\EE^d_n$ we
associate its capacity which is the random variable
\[
V(E) = \sum_{e\in E} t(e).
\]
The max-flow min-cut theorem states that
\[
\phi_n(F_1, F_2, C) = \min \bigl\{ V(E) | E
\subset\EE^d_n \mbox{ is a } (F_1,F_2)
\mbox{-cutset in }C \bigr\}.
\]
We can achieve a better understanding of what a cutset is thanks to the
following correspondence. We associate to each edge $e\in\EE^{d,i}_n$
a plaquette $\pi(e)$ defined by
\[
\pi(e) = c(e) + \frac{1}{2n} \bigl( [-1, 1]^{i-1} \times\{ 0 \}
\times[-1, 1]^{d-i} \bigr),
\]
where $c(e)$ is the middle of the edge $e$. To a set of edges $E\subset
\EE^d_n$ we associate the set of the corresponding plaquettes $E^* =
\bigcup_{e\in E} \pi(e)$. If $E$ is a $(F_1,F_2)$-cutset, then $E^*$ looks
like a ``surface'' of plaquettes that separates $F_1$ from $F_2$ in
$C$; see Figure~\ref{figplaquette}. We do not try to give a proper
definition to the term ``surface'' appearing here. In terms of
plaquettes, the discrete max-flow min-cut theorem states that the
maximal flow from $F_1$ to $F_2$ in $C$, given a local constraint on
the maximal amount of water that can cross each edge, is equal to the
minimal capacity of a ``surface'' that cuts $F_1$ from $F_2$ in $C$.

We consider now streams, cutsets and flows in $\O_n$. The set of stream
functions associated to our flow problem is $\S_n(\G^1_n, \G^2_n, \O
_n)$. We will denote by $\phi_n$ the maximal flow $\phi_n(\G^1_n, \G
^2_n, \O_n)$. To each $f_n \in\S_n(\G^1_n,\G^2_n, \O_n)$, we associate
the vector measure $\vec\mu_n$, that we call the stream itself,
defined by
\[
\vec{\mu}_n = \vec\mu_n (f_n) =
\frac{1}{n^d} \sum_{e \in\EE^d_n} f_n(e)
\vec{e} \de_{c(e)},
\]
where $c(e)$ is the center of $e$. Notice that since $f_n \in\S_n(\G
^1_n, \G^2_n, \O_n)$, the condition~(i) implies that $f_n(e)=0$ for
all $e \notin\Pi_n$; thus the sum in the previous definition is
finite. A stream $\vec\mu_n$ is a \emph{rescaled} measure version
of a
stream function $f_n$. The vector measure $\vec{\mu}_n$ is defined on
$(\RR^d, \B(\RR^d))$ where $\B(\RR^d)$ is the collection of the Borel
sets of $\RR^d$ and takes values in $\RR^d$. In fact $\vec{\mu}_n
=(\mu
_n^1,\ldots,\mu_n^d)$ where $\mu_n^i$ is a signed measure on $(\RR
^d, \B
(\RR^d))$ for all $i\in\{1,\ldots,d\}$. We define the flow corresponding
to a stream $\vec\mu_n (f_n)$ as $\operatorname{flow}^{\mathrm
{disc}}_n(f_n)$ properly rescaled,
\[
\operatorname{flow}^{\mathrm{disc}}_n(\vec\mu_n ) =
\operatorname {flow}^{\mathrm{disc}}_n \bigl(\vec\mu_n
(f_n) \bigr) = \frac{1}{n^{d-1}} \operatorname{flow}^{\mathrm{disc}}_n
(f_n).
\]
We say that $\vec\mu_n = \vec\mu_n (f_n)$ is \emph{a maximal
stream} from
$\G^1_n$ to $\G^2_n$ in $\O_n$ if and only if
%
%
\begin{equation}
\label{fluxmax1} \operatorname{flow}^{\mathrm{disc}}_n(\vec
\mu_n) = \frac{ \phi_n
}{ n^{d-1}},
\end{equation}
and for any $e=[a,b]$ such that $a\in\G^1_n$ and $b\notin\G^1_n$, we
have $f_n(e) \vec{e}\cdot\overrightarrow{ab} \geq0$, that is,
%
%
\begin{equation}
\label{fluxmax2} f_n (e) \cases{ %
\geq0, & \quad$\mbox{if } \be\mbox{ is oriented from } a \mbox{ to }b\ \bigl(\mbox{i.e. }
\be= \langle a,b\rangle \bigr),$
\vspace*{2pt}\cr
\leq0, & \quad$\mbox{if } \be\mbox{ is oriented from } b \mbox{ to } a\ \bigl(
\mbox{i.e. }\be= \langle b,a \rangle \bigr).$}
\end{equation}
The set of admissible stream functions is random since the capacity
constraint on the stream is random. Thus $ \phi_n $ is random and the
set of admissible streams (resp., maximal streams) from $\G^1_n$ to $\G
^2_n$ in $\O_n$ is random too.

\begin{figure}

\includegraphics{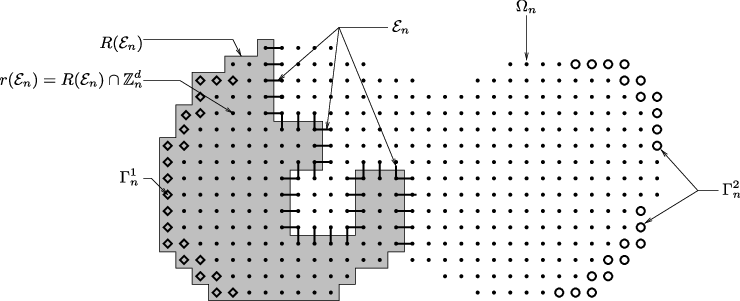}

\caption{A $(\G^1_n,\G^2_n)$-cutset $\E_n$ in $\O_n$ and the
corresponding sets $r(\E_n)$ and $R(\E_n)$.}
\label{figensemble}
\end{figure}

Let $\E_n$ be a $(\G^1_n, \G^2_n)$-cutset in $\O_n$. We say that
$\E_n$
is \emph{a minimal cutset} if and only if it realizes the minimum
%
%
\begin{equation}
\label{cutmin1} V(\E_n) = \phi_n
\end{equation}
and it has minimal cardinality, that is,
%
%
\begin{eqnarray}
\label{cutmin2} &&\card(\E_n) = \min \bigl\{ \card(\F_n
) | \F_n \mbox{ is a } \bigl(\G ^1_n,
\G^2_n \bigr)\mbox{-cutset in } \O_n \mbox{
and}
\nonumber
\\[-8pt]
\\[-8pt]
\nonumber
&&\hspace*{130pt}\qquad V(\F_n) = \phi_n \bigl(\G^1_n,
\G^2_n, \O_n \bigr) \bigr\},
\end{eqnarray}
where $\card(\E)$ denotes the cardinality of the set $\E$. We want to
see a cutset $\E_n$ as the ``boundary'' of a subset of $\O$. We define
the set $r( \E_n) \subset\ZZ^d_n$ by
\[
r(\E_n) = \bigl\{ x\in\O_n | \mbox{there exists a path
from } x \mbox{ to } \G^1_n \mbox{ in } \bigl(
\ZZ^d_n, \Pi_n \setminus
\E_n \bigr) \bigr\}.
\]
Then the edge boundary $\partial^e r(\E_n) $ of $r(\E_n)$, defined by
\[
\partial^e r(\E_n) = \bigl\{ e=[x,y] \in
\Pi_n | x \in r(\E_n) \mbox{ and } y \notin r(
\E_n) \bigr\},
\]
is exactly equal to $\E_n$. We consider a ``non discrete version''
$R(\E
_n)$ of $r(\E_n)$ defined by
\[
R (\E_n) = r (\E_n) + \frac{1}{2n}[-1,
1]^d.
\]
Notice that $\E_n = \partial^e (R(\E_n) \cap\Pi_n)$; thus the
sets $\E_n$
and $R(\E_n)$ completely define one each other; see Figure~\ref{figensemble}.
%
%
%

\begin{rem}
\label{remexistencecourant}
We want to study the asymptotic behavior of sequences of maximal
streams and minimal cutsets. For a fixed $n$ and given capacities, the
existence of at least one minimal cutset is obvious since there are
finitely many cutsets. The existence of at least one maximal stream is
not so obvious because of condition~(\ref{fluxmax2}). Under the
hypothesis that the capacities are bounded, we will prove in Section~\ref{secexistencecourant} that a maximal stream exists.
\end{rem}


\subsection{Brief presentation of the limiting objects}
\label{secdef2}

We consider a sequence $(\vec\mu_n^{\max})_{n\geq1}$ of maximal streams
and a sequence $(\E_n^{\min})_{n\geq1}$ of minimal cutsets. For each
$n$, $\vec\mu_n ^{\max}$ is a solution of a discrete random problem of
maximal flow, $\E_n^{\min}$ is a solution of a discrete random problem
of minimal cutset and by the max-flow min-cut theorem
\[
\operatorname{flow}^{\mathrm{disc}}_n \bigl(\vec
\mu_n^{\max} \bigr) = \frac
{V(\E_n^{\min})}{n^{d-1}}:= \frac
{\phi_n }{n^{d-1}}
,
\]
where $\phi_n$ stands for $\phi_n( \G^1_n, \G^2_n, \O_n)$. The
goal of
this article is to prove that:
\begin{itemize}
\item$(\vec\mu_n^{\max})_{n\geq1}$ converges in a way when $n$
goes to
infinity to a continuous stream $\vec\mu$ which is the solution of a
continuous deterministic max-flow problem to be precised;
\item$(\E_n^{\min})_{n\geq1}$ converges in a way when $n$ goes to
infinity to a continuous cutset $\E$ which is the solution of a
continuous deterministic min-cut problem to be precised;
\item these continuous deterministic max-flow and min-cut problems are
in correspondence, that is, the flow of $\vec\mu$ is equal to the
capacity of $\E$, and $\phi_n/n^{d-1}$ converges toward this constant.
\end{itemize}
We obtain these results, except that the continuous max-flow and
min-cut problems we define may have several solutions, thus we obtain
the convergence of the discrete streams $\vec\mu_n ^{\max}$
(resp., the discrete cutsets $\E_n^{\min}$) toward the set of the
solutions of a continuous deterministic max flow problem (resp., min-cut
problem). In this section, we try to present very briefly these
continuous max-flow and min-cut problems. A complete and rigorous
description will be given in Sections~\ref{secNozawa} and \ref
{secproba}. The aim of the present section is to give an intuitive idea
of the objects involved in the main theorems of Section~\ref{secdef3}.

The first quantity that has been studied is the maximal flow $\phi_n$;
however, a law of large numbers for $\phi_n$ is difficult to establish
in a general domain. It is considerably simpler in the following
situation. Let $\vv$ be a unit vector in $\RR^d$, let $Q(\vv)$ be a
unit cube centered at the origin having two faces orthogonal to $\vv$
and let
\[
F_1 = \{ x \in\partial Q | \overrightarrow{0 x} \cdot\vv<0 \},\qquad
F_2 = \{ x \in\partial Q | \overrightarrow{0 x} \cdot\vv>0 \}
\]
be, respectively, the upper half part and the lower half part of the
boundary of $Q$ in the direction $\vv$. Whenever $\EE(t(e)) < \infty$,
a subadditive argument yields the following convergence:
%
%
\begin{equation}
\label{eqdefnu} \lim_{n\rightarrow\infty} \frac{\phi_n (F_1, F_2, Q(\vv))}{n^{d-1}} = \nu(\vv)
\qquad\mbox{in } L^1,
\end{equation}
where $\nu(\vv)$ is deterministic and depends on the law of the
capacities of the edges, the dimension and $\vv$. The maximal flow
considered here is not well defined, since $F_1$ and $F_2$ are not sets
of vertices (a rigorous definition will be given in Section~\ref{secproba}), but equation (\ref{eqdefnu}) allows us to understand what
the constant $\nu(\vv)$ represents. By the max-flow min-cut theorem,
$\phi_n (F_1, F_2, Q(\vv))$ is the minimal capacity of a ``surface'' of
plaquettes that cuts $F_1$ from $F_2$ in $Q(\vv)$, thus a discrete
``surface'' whose boundary is spanned by $\partial Q(\vv)$. Thus the constant
$\nu(\vv)$ can be seen as the average asymptotic capacity of a
continuous unit surface normal to $\vv$. By symmetry we have $\nu(\vv)
= \nu(-\vv)$.

This interpretation of $\nu(\vv)$ provides in a natural way the desired
continuous deterministic min-cut problem. Indeed, if $\S$ is a ``nice''
surface (``nice'' means $\C^1$ among other things), it is natural to
define its capacity as
\[
\capa(\S) = \int_{\S\cap\O} \nu \bigl(\vv_{\S} (x)
\bigr) \,d\H^{d-1} (x),
\]
where $\H^{d-1}$ is the $(d-1)$-dimensional Hausdorff measure on $\RR
^d$, and $\vv_{\S} (x)$ is a unit vector normal to $\S$ at $x$. Exactly
as a discrete cutset $\E_n$ can be seen as the boundary of a set $R(\E
_n)$, we see $\S$ as the boundary of a set $F \subset\O$, and we
define $\capa(F) = \capa(\partial F)$. The continuous deterministic min-cut
problem we consider is the following:
\[
\phi_{\O}^{a}:= \inf \bigl\{ \capa(F) | F\subset\O,
\partial F \mbox{ is a surface separating } \G^1 \mbox{ from }
\G^2 \mbox{ in } \O \bigr\}.
\]
The above variational problem is loosely defined, since we did not give
a definition of $\capa(F)$ for all $F$, and we did not describe
precisely the admissible sets $F$: we should precise the regularity
required on $\partial F$ and what ``separating'' means. This will be
done in
Section~\ref{secproba}. We will denote by $\Sigma^a$ the set of the
continuous minimal cutsets, that is,
\[
\Sigma^a = \bigl\{ F \subset\O| F \mbox{ is ``admissible''
and } \capa(F) = \phi_{\O}^a \bigr\}.
\]
The variational problem $\phi_{\O}^{a}$ is a very good candidate to be
the continuous min-cut problem we are looking for, all the more since
it has been proved by the authors in the companion papers \cite
{CerfTheret09supc,CerfTheret09geoc} and \cite{CerfTheret09infc}
that under suitable hypotheses
\[
\lim_{n\rightarrow\infty} \frac{\phi_n }{n^{d-1}} = \phi_{\O}^{a}\qquad
\mbox{a.s.}
\]
This result is presented in Section~\ref{secproba}. By studying maximal
streams and minimal cutsets, we will give an alternative proof of this
law of large numbers for $\phi_n$.

We define now a continuous max-flow problem. A continuous stream in $\O
$ will be modeled by a vector field $\vec\sigma\dvtx  \RR^d \rightarrow
\RR
^d$ that must satisfy constraints equivalent to (i), (ii) and
(iii). For a ``nice'' stream $\vec\sigma$ (e.g., $\vec\sigma$ is
$\C^1$
on the closure $\overline{\O}$ of $\O$ and on $\RR^d \setminus
\overline{\O}$) these constraints would be:
\begin{longlist}[(iii$$')]
\item[(i$'$)]\emph{the stream is inside $\O$}: $\vec\sigma=0 $ on
$\RR^d
\setminus\overline{\O}$;
\item[(ii$'$)]\emph{capacity constraint}: $ \forall\vv\in\SS^{d-1},
\vec\sigma\cdot\vv\leq\nu(\vv)$ on $\RR^d$;
\item[(iii$'$)]\emph{conservation law}: $ \di\vec\sigma=0$ on $\O$ and
$\vec\sigma\cdot\vv_{\O} = 0$ on $\G\setminus(\G^1 \cup
\G^2)$.\vadjust{\goodbreak}
\end{longlist}
Here $\SS^{d-1}$ is the unit sphere of $\RR^d$ and $\vv_{\O}(x)$
denotes the exterior unit vector normal to $\O$ at $x$. The flow
corresponding to a ``nice'' stream $\vec\sigma$ would be
\[
\operatorname{flow}^{\mathrm{cont}}(\vec\sigma) = \int_{\G^1} -
\vec\sigma\cdot\vv_{\O} \,d\H ^{d-1}.
\]
Thus we obtain the following continuous max-flow problem:
\[
\phi_{\O}^b:= \sup\lleft\{ \operatorname{flow}^{\mathrm
{cont}}(
\vec\sigma) \left\vert %
\begin{array} {l} \vec\sigma\dvtx \RR^d
\rightarrow\RR^d \mbox{ is a stream inside } \O\mbox{ that satisfies}
\\
\mbox{the capacity constraint and the conservation law} \end{array} %
\right.
\rright\}.
\]
The above variational problem is loosely defined too, since we did not
give a definition of $\operatorname{flow}^{\mathrm{cont}}(\vec\sigma
)$ for all $\vec\sigma$,
and we did
not describe precisely the set of admissible streams $\vec\sigma$: we
should precise the regularity required on $\vec\sigma$ and adapt
conditions~(i$'$), (ii$'$) and (iii$'$) to $\vec\sigma$ in this
class of
regularity. This will be done in Section~\ref{secNozawa}. We will
denote by $\Sigma^b$ the set of the continuous maximal streams, that is,
\[
\Sigma^b = \bigl\{\vec\sigma\dvtx \RR^d \rightarrow
\RR^d | \vec\sigma \mbox{ is ``admissible''
and } \operatorname{flow}^{\mathrm
{cont}}( \vec\sigma) = \phi_{\O}^b
\bigr\}.
\]

We have also good reasons a priori to think that the variational
problem $\phi_{\O}^b$ is the max-flow problem we are looking for.
Indeed, various continuous versions of the max-flow min-cut theorem
have been proved (see, e.g., \mbox{\cite{Blum,Strang,Nozawa}}), and a main result of Nozawa's work \cite{Nozawa} is
precisely to prove that
\[
\phi_{\O}^b = \phi_{\O}^{a'},
\]
where $\phi_{\O}^{a'}$ is a variant of $\phi_{\O}^a$. Thanks to our
study of maximal flows and minimal cutsets, we will also recover this
continuous max-flow min-cut theorem in our setting.

\begin{rem}
We gave no argument a priori to justify that the sets $\Sigma^a$ and
$\Sigma^b$ are not empty. This will be a consequence of our results of
convergence. The fact that $\Sigma^b$ is not empty was already proved
by Nozawa in \cite{Nozawa}.
\end{rem}


\subsection{Main results}
\label{secdef3}

We denote by $\L^d$ the Lebesgue measure in $\RR^d$ and by $\C_b(\RR^d,
\RR)$ the set of the continuous bounded functions from $\RR^d$ to
$\RR
$. We define the distance $\mathfrak d$ on the subsets of $\RR^d$ by
\[
\forall E,F\subset\RR^d\qquad \mathfrak d(E,F) = \L^d (E
\triangle F),
\]
where $E \triangle F = (E \setminus F) \cup(F \setminus E)$
is the symmetric difference of $E$ and $F$.

We need some hypotheses on $(\O, \G^1, \G^2)$. We say that $\O$ is a
Lipschitz domain if its boundary $\G$ can be locally represented as the
graph of a Lipschitz function defined on some open ball of $\RR^{d-1}$.
We say that two $\C^1$ hypersurfaces $\S_1,\S_2$ intersect
transversally if for all $x\in\S_1 \cap\S_2$, the normal unit vector
to $\S_1$ and $\S_2$ at $x$ are not colinear. We gather here the
hypotheses we will make on $(\O, \G^1, \G^2)$
%
\renewcommand{\thehypo}{(H\arabic{hypo})}
\begin{hypo}\label{hypo1}
We suppose that $\O$ is a bounded open connected subset of $\RR^d$,
that it is a Lipschitz domain and that $\G$ is included in the union of
a finite number of oriented hypersurfaces of class $\C^1$ that
intersect each other transversally; we also suppose that $\G^1$ and
$\G
^2$ are open subsets of $\G$, that $ \inf\{ \| x-y\|, x\in\G^1
, y\in\G^2)>0$, and that their relative boundaries $\partial_{\G}
\G^1$
and $\partial_{\G} \G^2$ have null $\H^{d-1}$ measure.
\end{hypo}

We also make the following hypotheses on the law of the capacities:
%
\begin{hypo}\label{hypo2}
We suppose that the capacities of the edges are bounded by a constant
$M$, that is,
\[
\exists M < +\infty,\qquad \Lambda \bigl([0,M] \bigr) = 1.
\]
\end{hypo}
%
\begin{hypo}\label{hypo3}
We suppose that
\[
\Lambda \bigl(\{0\} \bigr) < 1 - p_c(d),
\]
where $p_c(d)$ is the critical parameter of edge Bernoulli percolation
on $(\ZZ^d, \EE^d)$.
\end{hypo}

We can now state our main results:
%
\begin{thmm}[(Law of large numbers for the maximal streams)]
\label{thmpcpal}
We suppose that the hypotheses \ref{hypo1} and \ref{hypo2} are
fulfilled. For all $n\geq1$, let $\vec\mu_n^{\max}$ be a random maximal
discrete stream from $\G^1_n$ to $\G^2_n$ in $\O_n$. Then $(\vec
\mu
_n^{\max})_{n\geq1}$ converges weakly a.s. toward the set $\Sigma^b$,
that is,
\[
\mbox{a.s.}, \forall f \in\C_b \bigl(\RR^d, \RR \bigr)\qquad \lim
_{n
\rightarrow\infty} \inf_{\vec\sigma\in\Sigma^b} \biggl\llVert \int
_{\RR^d} f \,d\vec\mu_n^{\max} - \int
_{\RR^d} f \vec\sigma \,d \L^d \biggr\rrVert = 0.
\]
\end{thmm}

\begin{thmm}[(Law of large numbers for the minimal cutsets)]
\label{thmpcpal2}
We suppose that the hypotheses \ref{hypo1}, \ref{hypo2} and \ref{hypo3}
are fulfilled. For all $n\geq1$, let $\E_n^{\min}$ be a minimal $(\G
^1_n,\G^2_n)$-cutset in $\O_n$. Then the sequence $(R(\E_n^{\min
}))_{n\geq1}$ converges a.s. for the distance $\mathfrak d$ toward the
set $\Sigma^a$, that is,
\[
a.s., \qquad\lim_{n\rightarrow\infty} \inf_{F \in\Sigma^a} \mathfrak d
\bigl(R \bigl(\E_n^{\min} \bigr), F \bigr) = 0.
\]
\end{thmm}
%
\begin{rem}
As we will see in Section~\ref{secproba}, condition \ref{hypo3} is
equivalent to $\nu\neq0$, where $\nu$ is the function defined by
equation (\ref{eqdefnu}). Thus if \ref{hypo3} is not satisfied, then
$\nu(\vv) = 0$ for all $\vv$, $\capa(F) = 0$ for every admissible
continuous cutset $F$ and the variational problem $\phi_{\O}^a$ is trivial.
\end{rem}

The two previous theorems lead to the following corollary:
%
\begin{cor}
\label{corunicite}
We suppose that hypotheses \ref{hypo1} and \ref{hypo2} are fulfilled.
If $\Sigma^b$ is reduced to a single stream $\vec\sigma$, then any
sequence of maximal streams $(\vec\mu_n^{\max})_{n\geq1}$ converges
a.s. weakly to $\vec\sigma\L^d$.\vadjust{\goodbreak} If hypothesis \ref{hypo3} is also
fulfilled and if $\Sigma^a$ is reduced to a single set $F$, then for
any sequence of minimal cutsets $(\E_n^{\min})_{n\geq1}$, the
corresponding sequence $(R(\E_n^{\min}))_{n\geq1}$ converges a.s. for
the distance $\mathfrak d$ toward $F$.
\end{cor}

\begin{rem}
We believe that the uniqueness of the maximal stream or the uniqueness
of the minimal cutset in the continuous setting may happen, or not,
depending on the domain $\O$, the sets $\Gamma^i, i=1,2$ and the
function $\nu$ (thus on the law of the capacities $\Lambda$); however
we do not handle this question here.
\end{rem}
During the proof of Theorem \ref{thmpcpal}, we prove the key
inequalities to obtain the following lemma:
%
\begin{lem}
\label{lemineg}
We suppose that hypotheses \ref{hypo1} and \ref{hypo2} are fulfilled,
and we consider the continuous variational problems $\Sigma^a$ and
$\Sigma^b$ associated to the function $\nu\dvtx  \SS^{d-1} \rightarrow
\RR
^+$. For every admissible continuous stream $\vec\sigma$, for every
admissible set $F$, we have
\[
\operatorname{flow}^{\mathrm{cont}}(\vec\sigma) \leq\capa(F).
\]
\end{lem}
The proof of Theorems \ref{thmpcpal} and \ref{thmpcpal2} relies on a
compactness argument. Combining this argument, Theorems \ref{thmpcpal},
\ref{thmpcpal2} and Corollary \ref{lemineg}, we obtain the two
following theorems:
%
\begin{thmm}[(Max-flow min-cut theorem)]
\label{cormaxflowmincut}
We suppose that the hypotheses~\ref{hypo1} and \ref{hypo2} are
fulfilled, and we consider the continuous variational problems $\Sigma
^a$ and $\Sigma^b$ associated to the function $\nu\dvtx  \SS^{d-1}
\rightarrow\RR^+$. Then there exists at least an admissible continuous
stream $\vec\sigma$ such that $\phi_{\O}^b = \operatorname
{flow}^{\mathrm{cont}}(\vec\sigma
)$, there
exists at least an admissible set $F$ such that $\phi_{\O}^a = \capa
(F)$, and we have the following max-flow min-cut theorem:
\[
\phi_{\O}^a = \phi_{\O}^b:=
\phi_{\O}.
\]
\end{thmm}

\begin{thmm}[(Law of large numbers for the maximal flows)]
\label{corLGNflow}
Suppose that hypotheses \ref{hypo1} and \ref{hypo2} are fulfilled. Then
we have
\[
\lim_{n\rightarrow\infty} \frac{\phi_n}{ n^{d-1}} = \phi_{\O}\qquad \mbox{a.s.}
\]
\end{thmm}

\begin{rem}
As will be explained in the next section, the last two theorems do not
state new results, since the continuous max-flow min-cut theorem we
obtain is a particular case of the one studied by Nozawa in \cite
{Nozawa}, and the law of large numbers for the maximal flows has been
proved by the authors in \cite
{CerfTheret09supc,CerfTheret09geoc,CerfTheret09infc} under a weaker
assumption on $\Lambda$. However, these results are recovered here by
new methods, which are more natural. Indeed, the law of large numbers
for $\phi_n$ was proved in \cite
{CerfTheret09supc,CerfTheret09geoc,CerfTheret09infc} by a study of its
lower and upper large deviations around $\phi_{\O}$. The study of the
upper large deviations \cite{CerfTheret09supc} is replaced here by the
study of a sequence of maximal streams, which is the most original part
of this article and gives a better understanding of the model. The
study of the lower large deviations \cite{CerfTheret09infc} is replaced
by the study of a sequence minimal cutsets. The techniques are the same
in both cases, but we change our point of view. To conclude, we use in
both proofs the result of polyhedral approximation presented in \cite
{CerfTheret09geoc}.
\end{rem}


\section{Background}
\label{secback}

We present now the mathematical background on which our work relies. It
is the occasion to give a proper description of the variational
problems involved in our theorems.

\begin{figure}[b]

\includegraphics{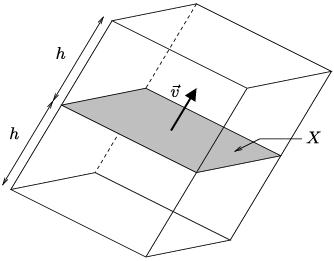}

\caption{Cylinder $\cyl(X,h)$.}
\label{figcylindre}
\end{figure}

\subsection{Some geometric tools}
\label{secgeom}

\begin{figure}

\includegraphics{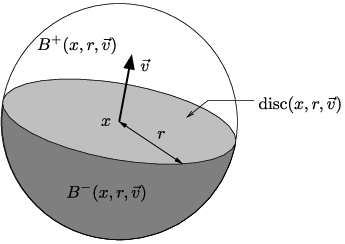}

\caption{Ball $B(x,r)$.}
\label{figboule}
\end{figure}

We start with simple geometric definitions.\vspace*{-1pt} For a subset $X$ of
$\mathbb{R}^d$, we denote by $\overline{X}$ the closure of $X$, by
$\stackrel{\circ}{X}$ the interior of $X$, by $X^c$ the set $\RR^d
\setminus X$ and by $\mathcal{H}^s (X)$ the $s$-dimensional
Hausdorff measure of $X$. The
$r$-neighborhood $\V_i(X,r)$ of $X$ for the distance $d_i$, that can be
the Euclidean distance if $i=2$ or the $L^\infty$-distance if
$i=\infty$,
is defined by
\[
\V_i (X,r) = \bigl\{ y\in\RR^d | d_i(y,X)<r
\bigr\}.
\]
If $X$ is a subset of $\RR^d$ included in an hyperplane of $\RR^d$
and of
codimension $1$ (e.g., a nondegenerate hyperrectangle), we denote by
$\hyp(X)$ the hyperplane spanned by~$X$, and we denote by $\cyl(X,
h)$ the
cylinder of basis $X$ and of height $2h$ defined by
\[
\cyl(X,h) = \bigl\{x+t \vv| x\in X, t\in [-h,h] \bigr\},
\]
where $\vv$ is one of the two unit vectors orthogonal to $\hyp(X)$ (see
Figure~\ref{figcylindre}).
For $x\in\RR^d$, $r\geq0$ and a unit vector $\vv$, we denote by
$B(x,r)$ the closed
ball centered at $x$ of radius $r$, by $\disc(x,r,\vv)$ the closed
disc centered at $x$ of radius $r$ and normal vector $v$,
and by $B^+ (x,r,\vv)$ [resp., $B^- (x,r,\vv)$] the upper (resp., lower)
half part of $B(x,r)$ where the direction is determined by $v$ (see
Figure~\ref{figboule}), that is,
\begin{eqnarray*}
B^+ (x,r,\vv) &= &\bigl\{y\in B(x,r) | \overrightarrow{xy}\cdot\vv \geq0 \bigr\},
\\
B^- (x,r,\vv) &=& \bigl\{y\in B(x,r) | \overrightarrow{xy}\cdot\vv \leq0 \bigr\}.
\end{eqnarray*}
%

%

%
We denote by $\alpha_p$ the volume of the unit ball in $\RR^p$,
$p\geq
1$. Thus $\alpha_d$ is the volume of a unit ball in $\RR^d$, and
$\alpha_{d-1}$ the $\H^{d-1}$ measure of a unit disc in $\RR^d$. We say
that a domain $\O$ of $\RR^d$ has Lipschitz boundary if its boundary
can be locally represented as the graph of a Lipschitz function defined
on some open ball of $\RR^{d-1}$. We say that a vector $\vv\neq0$
defines a rational direction if there exists a positive real number
$\lambda$ such that $\lambda\vv$ has rational coordinates. It is
equivalent to require that there exists a positive real number $\lambda
'$ such that $\lambda' \vv$ has integer coordinates. We denote by
$\SS
^{d-1}$ the unit sphere in $\RR^d$, and by $\widehat\SS^{d-1}$ the set
of the unit vectors of $\RR^d$ defining a rational direction. Notice
that $\widehat\SS^{d-1}$ is dense in $\SS^{d-1}$.

Two submanifolds $E$ and $F$ of a given finite dimensional smooth
manifold are said to intersect transversally if at every point of
intersection, their tangent spaces at that point span the tangent space
of the ambient manifold at that point; see Section~5 in \cite
{GuilleminPollack}. When a hypersurface $\S$ is piecewise of class $\C
^1$, we say that $\S$ is transversal to $\G$ if for all $x\in\S\cap
\G
$, the normal unit vectors to $\S$ and $\G$ at $x$ are not colinear; if
the normal vector to $\S$ (resp., to $\G$) at $x$ is not well defined,
this property must be satisfied by all the vectors which are limits of
normal unit vectors to $\S$ (resp., $\G$) at $y\in\S$ (resp., $y\in
\G$)
when we send $y$ to $x$---there is at most a finite number of such
limits. We say that a subset $P$ of $\RR^d$ is polyhedral if its
boundary $\partial P$ is included in the union of a finite number of
hyperplanes.

Let $E$ be a subset of $\RR^d $. We say that $E$ is $p$-rectifiable if
and only if there exists a Lipschitz function mapping some bounded
subset of $\RR^p$ onto $E$; see Definition~3.2.14 in \cite{FED}. We
define the $p$ dimensional upper (resp., lower) Minkowski content $\M
^{p,+} (E)$ [resp., $\M^{p,-} (E)$] of $E$ by
\[
\M^{p,+} (E) = \limsup_{r\rightarrow0^+} \frac{\L^d (\V
_2(E,r))}{\alpha_{d-p} r^{d-p}}\quad
\mbox{and} \quad\M^{p,-} (E) = \liminf_{r\rightarrow0^+}
\frac{\L^d (\V_2(E,r))}{\alpha_{d-p}
r^{d-p}}.
\]
If $\M^{p,+} (E) = \M^{p,-}(E)$, their common value is called the $p$
dimensional Minkowski content of $E$, which is denoted by $\M^p (E)$;
see Definition 3.2.37 in~\cite{FED}. According to Theorem 3.2.39 in
\cite{FED}, if $E$ is a closed $p$-rectifiable subset of $\RR^d$, then
its $p$ dimensional Minkowski content exists, and we have
\[
\M^p (E) = \H^p (E).
\]

We need some properties of sets of finite perimeter. We denote by $\C
_c^k (A,B)$, for $A\subset\RR^p$ and $B\subset\RR^q$, the set of
functions of class $\C^k$ defined on $\RR^p$, that takes values in $B$
and whose domain is included in a compact subset of $A$. For a subset
$F$ of $\RR^d$, we define the perimeter of $F$ in $\O$ by
\begin{eqnarray*}
\P(F, \O) &=& \sup \biggl\{\int_F \di\vec f \,d
\L^d | \vec f\in\C _c^\infty \bigl(\O,
\RR^d \bigr),\\
 &&\hspace*{28pt}\vec f (x) \in B(0,1) \mbox{ for all }x \in \O \biggr\}
,
\end{eqnarray*}
where $\di$ is the usual divergence operator. We denote by $\partial
F$ the
boundary of $F$. The reduced boundary of a set of finite perimeter $F$,
denoted by $\partial^* F$, consists of the points $x$ of $\partial
F$ such that:
\begin{itemize}
\item $| \vec\nabla\mathbh{1}_{F} | (B(x,r)) > 0$ for
any $r>0$,
\item if $\vec w_r (x) = - \vec\nabla\mathbh{1}_{F}
(B(x,r)) /
| \vec\nabla\mathbh{1}_{F} | (B(x,r))$ then, as $r$ goes to $0$,
$\vec
w_r(x)$ converges toward a unit vector $\vv_F (x)$,
\end{itemize}
where $\mathbh{1}_F$ is the indicator function of $F$, $\vec\nabla
\mathbh{1}_{F}$ is the distributional derivative of $\mathbh{1}_F$
defined by
\[
\forall\vec h \in\C^\infty_c \bigl(\RR^d,
\RR^d \bigr)\qquad \int_{\RR^d} \vec h \cdot\vec\nabla
\mathbh{1}_F \,d\L^d = - \int_{\RR^d}
\mathbh{1}_F \di\vec h \,d\L^d
\]
and $| \vec\nabla\mathbh{1}_{F} |$ is the total variation measure of
$\vec\nabla\mathbh{1}_{F}$ defined by
\begin{eqnarray*}
&&\forall A \in\B \bigl(\RR^d \bigr)\qquad | \vec\nabla
\mathbh{1}_{F} | (A) = \sup \biggl\{ \int_A
\mathbh{1}_{F} \di\vec h \,d \L^d | \vec h \in
\C^\infty_c \bigl(A, \RR^d \bigr),\\
&&\hspace*{175pt} \vec h (x) \in
B(0,1) \mbox{ for all } x\in A \biggr\}.
\end{eqnarray*}
At any point $x$ of $\partial^* F$, the vector $\vv_F (x)$ is also the
measure theoretic exterior normal to $F$ at $x$, that is,
\begin{eqnarray*}
\lim_{r\rightarrow0} r^{-d} \L^d \bigl(B^-
\bigl(x,r,\vv_F (x) \bigr) \cap F^c \bigr) &=& 0\quad \mbox{and}\\
\lim_{r\rightarrow0} r^{-d} \L^d \bigl(B^+
\bigl(x,r,\vv_F (x) \bigr) \cap F \bigr) &=& 0,
\end{eqnarray*}
where $F^c = \RR^d \setminus F$. The set of functions of bounded
variations in $\O$, denoted by $\operatorname{BV}(\O)$, is the set of all functions
$u\in L^1(\O, \RR)$ such that
\begin{eqnarray*}
| \vec\nabla u | (\O)&:=& \sup \biggl\{ \int_{\O} u \di\vec h
\,d\L^d | \vec h \in\C^\infty_c \bigl(\O,
\RR^d \bigr), \vec h (x) \in B(0,1) \mbox { for all } x\in\O \biggr\}\\
&<& \infty.
\end{eqnarray*}
By definition, a set $F$ has finite perimeter in $\O$ if and only if
$\mathbh{1}_F$ has bounded variations in $\O$,
\[
\P(F,\O) <\infty\quad\iff\quad\mathbh{1}_F \in \operatorname{BV}(\O).
\]
More details about functions of bounded variations and sets of finite
perimeters can be found in \cite{GI}.


\subsection{Continuous max-flow min-cut theorem}
\label{secNozawa}

The (discrete) max-flow min-cut theorem has been transposed into a
continuous setting by various mathematicians. We present now one of
these works on continuous max-flow min-cut theorem, the article \cite
{Nozawa} by Nozawa. Indeed, the framework chosen by Nozawa is
particularly well adapted to our model.

We give here a presentation of the part of Nozawa's paper that we will
use. We adapt some notation of Nozawa to fit within ours, and we focus
on a particular case of one of the theorems presented in \cite{Nozawa}.
We try to keep the exposition self-contained, and we refer to \cite
{Nozawa} for more details. Nozawa considers a bounded domain $\O$ of
$\RR^d$ with Lipschitz boundary $\G$, and two disjoint Borel subsets
$\G
^1$ and $\G^2$ of $\G$. A stream in $\O$ is a vector field $\vec
\sigma
\in L^\infty( \O\rightarrow\RR^d, \L^d)$. The fact that there is no
loss or creation of fluid inside $\O$ is expressed by the condition
%
%
\begin{equation}
\label{eqNozawadiv} \di\vec\sigma= 0 \qquad\mbox{on }\O,
\end{equation}
where the divergence must be understood in the distributional, that is,
$\di\vec\sigma$ is defined on $\O$ by
\[
\forall h \in\C_c^\infty(\O, \RR)\qquad \int
_{\RR^d} h \di\vec \sigma \,d\L^d = - \int
_{\RR^d} \vec\sigma\cdot\vec\nabla h \,d\L^d.
\]
Thus equation (\ref{eqNozawadiv}) means that
\[
\forall h \in\C_c^\infty(\O, \RR)\qquad \int
_{\RR^d} \vec\sigma \cdot\vec\nabla h \,d\L^d = 0.
\]

\begin{rem}
\label{remdiv1}
The divergence $\di\vec\sigma$ is defined as a distribution. Thus
it is
an abuse of notation to write $\int_{\RR^d} h \di\vec\sigma \,d\L^d$
instead of $\langle\di\vec\sigma, h \rangle$, the action of the
distribution $\di\vec\sigma$ on the function $h$. In \cite{Nozawa}
Nozawa considers in fact vector fields $\vec\sigma$ such that $\di
\vec
\sigma\in L^d(\O, \L^d)$ in the distributional sense, that is, such
that there exists a real function $G \in L^d (\O,\L^d)$ satisfying
\[
\forall h \in\C^\infty_c (\O, \RR)\qquad \int
_{\O} \vec\sigma \cdot\vec\nabla h \,d\L^d = -
\int_{\O} G h \,d\L^d.\vadjust{\goodbreak}
\]
This implies that $\di{\vec{\sigma}}$ is a distribution of order
$0$ on
$\O$, thus by the Riesz representation theorem (see Theorem 6.19 in
\cite{Rudinrc}) it corresponds to a Radon measure that we denote by
$\di\vec\sigma\L^d\vert_{\O}$ and $\di\vec\sigma\L^d\vert
_{\O}
= G \L^d\vert_{\O}$. Of course, $\di\vec\sigma=0$ on $\O$ (as
defined above) implies that such a function $G$ exists, it is the null
function on~$\O$. Thus, with a slight abuse of notation, we say that
equation (\ref{eqNozawadiv}) is equivalent to
\[
\di\vec\sigma= 0,\qquad \L^d\mbox{-a.e. on }\O,
\]
which means that the associated function $G$ in $L^d(\O, \L^d)$ is
equal to $0$ a.e. on $\O$. We will see in Section~\ref{secdiv} that for
all the vector fields $\vec\sigma$ that we will consider, $\di\vec
\sigma
$ is in fact a distribution of order $0$ on $\RR^d$ itself. Thus by the
Riesz representation theorem it is a Radon measure that we denote by
$\di\vec\sigma\L^d$. More details about distributions can be found
in \cite{Schwartz1,Schwartz2}.
\end{rem}

A stream $\vec\sigma$ from $\G^1$ to $\G^2$ in $\O$ must also satisfy
some boundary conditions: the fluid enters in $\O$ through $\G^1$, and
no fluid can cross $\G\setminus(\G^1 \cup\G^2)$. Let us
translate this in a mathematical language. According to Nozawa in~\cite
{Nozawa}, Theorem 2.1, there exists a linear mapping $\gamma$ from
$\operatorname{BV}(\O)$ to $L^1(\G\rightarrow\RR, \H^{d-1})$ such that, for any $u
\in \operatorname{BV}(\O)$,
%
%
\begin{equation}
\label{eqtrace} \lim_{\rho\rightarrow0, \rho>0} \frac{1}{\L^d(\O\cap B(x, \rho
))} \int
_{\O\cap B(x,\rho)} \bigl|u(y) - \gamma(u) (x)\bigr| \,d\L^d (y) = 0
\end{equation}
for $\H^{d-1}$-a.e. $x\in\G$. The function $\gamma(u)$ is called the
trace of $u$ on $\G$. Let $\vv_{\O}(x)$ be the exterior unit vector
normal to $\O$ at $x\in\G$. The vector $\vv_{\O}$ is defined $\H
^{d-1}$-a.e. on $\G$ and the map $x\in\G\mapsto\vv_{\O} (x)$ is
$\H
^{d-1}$-measurable. According to Nozawa in~\cite{Nozawa}, Theorem 2.3,
for every $\vec\rho= (\rho_1,\ldots,\rho_d) \dvtx \O\rightarrow\RR
^d$ such
that $\rho_i \in L^\infty(\O\rightarrow\RR, \L^d)$ for all
$i=1,\ldots,d$ and $\di\vec\rho\in L^d(\O\rightarrow\RR, \L
^d)$, there
exists $g \in L^\infty(\G\rightarrow\RR, \H^{d-1})$ defined by
\[
\forall u \in W^{1,1} (\O)\qquad \int_{\G} g \gamma(u)
\,d \H ^{d-1} = \int_{\O} \vec\rho\cdot\vec{\nabla}
u \,d \L^d + \int_{\O} u \di\vec\rho \,d
\L^d.
\]
The function $g$ is denoted by $\vec\rho\cdot\vv_{\O}$. Any
stream $\vec
{\sigma}$ satisfies the conditions required to define $\vec{\sigma}
\cdot
\vv_{\O}$, and the definition is simpler since $\di\vec{\sigma}
=0$ $\L
^d$-a.e. on $\O$,
\[
\forall u \in W^{1,1} (\O)\qquad\int_{\G} (\vec{
\sigma}\cdot\vv _{\O} ) \gamma(u) \,d \H^{d-1} = \int
_{\O} \vec{\sigma} \cdot\vec {\nabla } u \,d \L^d.
\]
We impose the following boundary conditions on any stream $\vec\sigma$
from $\G^1$ to $\G^2$ in~$\O$:
%
%
\begin{eqnarray}
\label{eqNozawabords} \vec\sigma\cdot\vv_{\O} &\leq&0 \qquad\H^{d-1}
\mbox{-a.e. on } \G ^1\quad \mbox{and}
\nonumber
\\[-8pt]
\\[-8pt]
\nonumber
\vec\sigma\cdot\vv_{\O} &=&
0\qquad \H^{d-1} \mbox{-a.e. on } \G \setminus \bigl(\G^1\cup
\G^2 \bigr).
\end{eqnarray}
Finally, Nozawa puts a local capacity constraint on any stream $\vec
\sigma$,
%
%
\begin{equation}
\label{eqNozawacons} \L^d\mbox{-a.e. on } \O, \forall\vv\in
\SS^{d-1} \qquad\vec\sigma\cdot\vv\leq\nu(\vv),
\end{equation}
where $\SS^{d-1}$ is the set of all unit vectors in $\RR^d$, and $\nu:
\RR^d \rightarrow\RR^+$ is a continuous convex function that satisfies
$\nu(\vv) = \nu(-\vv)$. In our setting this function $\nu$ is the one
we have unformally defined in equation (\ref{eqdefnu}) and that we will
properly define in Section~\ref{secproba}.

To each admissible stream, that is, to each vector field $\vec\sigma
\in
L^{\infty} (\O\rightarrow\RR^d, \L^d) $ satisfying (\ref
{eqNozawadiv}), (\ref{eqNozawabords}) and (\ref{eqNozawacons}), we
associate its flow $\operatorname{flow}^{\mathrm{cont}}(\vec\sigma
)$ defined by
\[
\operatorname{flow}^{\mathrm{cont}}(\vec\sigma) = \int_{\G^1} -
\vec\sigma\cdot\vv_{\O} \,d \H ^{d-1},
\]
which is the amount of water that enters into $\O$ along $\G^1$
according to the stream~$\vec\sigma$. Nozawa investigates the behavior
of the maximal flow over all admissible continuous streams; that is, he
considers the following continuous max-flow problem:
%
%
{\fontsize{10.3}{12.6}{\selectfont
\begin{eqnarray}
\label{eqmaxflowNozawa}\qquad && \phi_{\O}^{(M)} = \sup\lleft\{
\operatorname{flow}^{\mathrm
{cont}}(\vec\sigma) \left\vert %
\begin{array} {l}
\vec\sigma\in L^\infty \bigl(\Omega\rightarrow\RR^d, \L
^d \bigr), \di\vec\sigma=0\ \L^d\mbox{-a.e. on } \O,
\\
\vec\sigma\cdot\vv\leq\nu(\vv) \mbox{ for all } \vv\in\SS ^{d-1}\
\L^d\mbox{-a.e. on }\O,
\\
\vec\sigma\cdot\vv_{\O} \leq0\ \H^{d-1} \mbox{-a.e. on }
\G^1,
\\
\vec\sigma\cdot\vv_{\O} = 0\ \H^{d-1} \mbox{-a.e. on } \G
\setminus \bigl(\G^1\cup\G^2 \bigr) \end{array}\right.
\rright\}.
\end{eqnarray}}}
Any vector field $\vec\sigma\in L^\infty(\O\rightarrow\RR^d, \L^d)$
can be extended to $\RR^d$ by defining $\vec\sigma= 0 $ $\L
^d$-a.e. on
$\O^c$. Thus the previous variational problem can be rewritten as
%
%
\begin{eqnarray}
\label{eqmaxflowNozawabis} \qquad\phi_{\O} ^b& =&
\phi_{\O}^{(M)}
\nonumber
\\[-8pt]
\\[-8pt]
\nonumber
&= &\sup\lleft\{ \operatorname
{flow}^{\mathrm{cont}}(\vec\sigma) \left\vert %
\begin{array} {l} \vec\sigma
\in L^\infty \bigl(\RR^d \rightarrow\RR^d, \L
^d \bigr), \vec\sigma=0\ \L^d \mbox{-a.e. on }
\O^c,
\\
\di\vec \sigma=0\ \L^d\mbox{-a.e. on } \O,
\\
\vec\sigma\cdot\vv\leq\nu(\vv) \mbox{ for all } \vv\in\SS ^{d-1}\
\L^d\mbox{-a.e. on }\O,
\\
\vec\sigma\cdot\vv_{\O} \leq0\ \H^{d-1} \mbox{-a.e. on }
\G^1,
\\
\vec\sigma\cdot\vv_{\O} = 0\ \H^{d-1} \mbox{-a.e. on } \G
\setminus \bigl(\G^1\cup\G^2 \bigr) \end{array}\right.
\rright\}.
\end{eqnarray}
This variational problem is exactly the one we have informally
presented in Section~\ref{secdef2} as $\phi_{\O}^b$ and that appears in
the main results presented in Section~\ref{secdef3}. Thus we have now a
precise definition of the set of admissible streams and of the flow of
any admissible stream $\vec\sigma$. Thus the set $\Sigma^b$
appearing in
Theorem \ref{thmpcpal} is defined by
\[
\Sigma^b = \lleft\{\vec\sigma\in L^\infty \bigl(
\RR^d \rightarrow\RR^d, \L^d \bigr) \left\vert
\begin{array} {l} \vec\sigma=0\ \L^d \mbox{-a.e. on }
\O^c, \di\vec\sigma=0\ \L^d\mbox{-a.e. on } \O,
\\
\vec\sigma\cdot\vv\leq\nu(\vv) \mbox{ for all } \vv\in\SS ^{d-1}\ \L^d\mbox{-a.e. on }\O,
\\
\vec\sigma\cdot\vv_{\O} \leq0\ \H^{d-1} \mbox{-a.e. on }
\G^1,
\\
\vec\sigma\cdot\vv_{\O} = 0\ \H^{d-1} \mbox{-a.e. on } \G
\setminus \bigl(\G^1\cup\G^2 \bigr),
\\
\operatorname{flow}^{\mathrm{cont}}(\vec\sigma) = \phi_{\O}^b
\end{array} %
\right.
\rright\}.
\]
We emphasize the fact that the constant $\phi_{\O}^b$ and the set
$\Sigma^b$ depend on $\O, \G^1, \G^2$ and $\nu$.

Nozawa defines a corresponding min-cut problem. A continuous cutset is
an hypersurface included in $\O$. Such a surface is seen as the
boundary of a sufficiently regular set $S \subset\O$, that is, a set
$S$ of finite perimeter in $\O$. To express the fact that the boundary
of $S$ in $\O$, $\O\cap\partial S$, cuts $\G^1$ from $\G^2$ in $\O$,
Nozawa imposes some boundary conditions on the indicator function
$\mathbh{1}_{S}$:
\[
\gamma(\mathbh{1}_{S}) = 1 \qquad\H^{d-1} \mbox{-a.e. on }
\G^1 \quad\mbox{and}\quad \gamma(\mathbh{1}_{S}) = 0\qquad
\H^{d-1} \mbox{-a.e. on } \G ^2.
\]
It means in a weak sense that $\G^1$ is ``in'' $S$ and $\G^2$ is not
``in'' $S$. In the max-flow problem (\ref{eqmaxflowNozawa}), $\nu(\vv)$
is the local capacity of the medium in the direction~$\vv$; thus the
capacity of the surface $\O\cap\partial S$ can be defined as
\[
\int_{\O\cap\partial^* S} \nu \bigl(\vv_{S} (x) \bigr) \,d
\H^{d-1} (x).
\]
In the previous equation, the integral is taken over the reduced
boundary $\partial^* S$ of~$S$, where the exterior normal to $S$ is
defined. Nozawa
investigates the behavior of the minimal capacity of a continuous
cutset; that is, he considers the following min-cut problem:
%
%
\begin{equation}
\label{eqmincutNozawa} \phi_{\O}^{(m)} = \inf\lleft\{
\int_{\O\cap\partial^* S} \nu \bigl(\vv _{S} (x) \bigr) \,d
\H^{d-1} (x) \left\vert %
\begin{array} {l} S \subset\O,
\mathbh{1}_{S} \in \operatorname{BV}(\O),
\\
\gamma(\mathbh{1}_{S}) = 1\ \H^{d-1} \mbox{-a.e. on }
\G^1,
\\
\gamma(\mathbh{1}_{S}) = 0\ \H^{d-1} \mbox{-a.e. on }
\G^2 \end{array} %
\right.
\rright\}.\hspace*{-35pt}
\end{equation}
He obtains the following continuous max-flow min-cut theorem:
%
\begin{thmm}[(Nozawa)]
\label{thmNozawa}
We suppose that $\O$ is a bounded domain of $\RR^d$ with Lipschitz
boundary $\G$, and that $\G^1$ and $\G^2$ are two disjoint Borel
subsets of $\G$. The following equality holds:
\[
\phi_{\O}^{(M)} = \phi_{\O}^{(m)} <
\infty.
\]
Moreover, there exists a maximal continuous stream; that is, there
exists a vector field $\vec\sigma$ as required in (\ref
{eqmaxflowNozawa}) such that $\operatorname{flow}^{\mathrm
{cont}}(\vec\sigma) = \phi_{\O}^{(M)}$.
\end{thmm}

\begin{rem}
For the interested reader, we explain how to deduce Theorem~\ref
{thmNozawa} from \cite{Nozawa}. We do not define all the notation
appearing here; they come from \cite{Nozawa}. We consider the max-flow
problem $(M\Phi_2)$ and the min-cut problem $(M\G_2)$ defined in
Section~5 of \cite{Nozawa}, pages 834 and 839. As suggested in the last
remark of \cite{Nozawa}, page 841, we fix $\alpha_t = \alpha'_t = 0$
$\H
^{d-1}$-a.e. on $\G$, for all $t\in T = \NN$, and $\G_t (x) =\{0\}$ for
all $x\in\O$ and for all $t\geq1$. For $x\in\O$ we define
\[
\G_0 (x) = \G_0 = \bigl\{ \vec w \in\RR^d
| \forall\vv\in\SS ^{d-1}, \vec w \cdot\vv\leq\nu(\vv) \bigr\},
\]
that does not depend on $x$ in our setting. The set $\G_0$ is the Wulff
crystal associated to $\nu$. It is a compact convex set since $\nu$ is
convex and bounded on $\SS^{d-1}$. Since $\nu$ is convex and
continuous, it is stated in Proposition 14.1 in \cite{Cerf:StFlour} that
\[
\forall\vv\in\SS^{d-1} \qquad\nu(\vv) = \sup\{ \vv\cdot\vec w | \vec w \in
\G_0 \}.
\]
Since $\nu(-\vv) = \nu(\vv)$, we obtain
\[
\beta_{\G_0 } \bigl(-\vv_S (x), x \bigr) = \sup \bigl\{ -
\vv_S (x) \cdot\vec w | \vec w \in\G_0 \bigr\} = \nu
\bigl(-\vv_S(x) \bigr) = \nu \bigl(\vv_S (x) \bigr).
\]
In this setting $(M\G_2)$ corresponds exactly to the min-cut problem
(\ref{eqmincutNozawa}), and $(M\Phi_2)$ corresponds almost to the
max-flow problem (\ref{eqmaxflowNozawa}), except that the goal is to
maximize $+ \int_{\G^1} \vec\sigma\cdot\vv_{\O} \,d \H^{d-1}$ on
streams $\vec\sigma$ satisfying $\vec\sigma\cdot\vv_{\O} \geq
0$ $\H
^{d-1}$-a.e. on $\G^1$. Since all the others conditions on $\vec
\sigma$
are satisfied by $- \vec\sigma$, $(M\Phi_2)$ is completely
equivalent to
(\ref{eqmaxflowNozawa}). Combining Theorems 5.3 and~5.6 in
\cite
{Nozawa}, we obtain Theorem \ref{thmNozawa}.
\end{rem}

\begin{rem}
The variational problem $\phi_{\O}^{(m)}$ is not exactly the same as
$\phi_{\O}^a$, the continuous min-cut problem we have informally
presented in Section~\ref{secdef2} and that appears in the main results
of Section~\ref{secdef3}. In fact, the variational problem $\phi_{\O
}^{(m)}$ is not well posed, since the infimum may not be reached by any
admissible set $F$. Since we want to prove the convergence of a
sequence of discrete minimal cutsets to the set of minimal continuous
cutsets, we have to consider another variational problem. This is done
in the next section.
\end{rem}


\subsection{Probabilistic background}
\label{secproba}

The study of the maximal flow in first passage percolation started in
1987 with the work of Kesten \cite{Kesten:flows}. We do not give here a
complete state of the art of all the results known in this domain. We
choose to present only the results that we will rely on and that
motivate our work. For a more complete introduction to this subject we
refer to \cite{CerfTheret09geoc}, Section~3.

\begin{figure}

\includegraphics{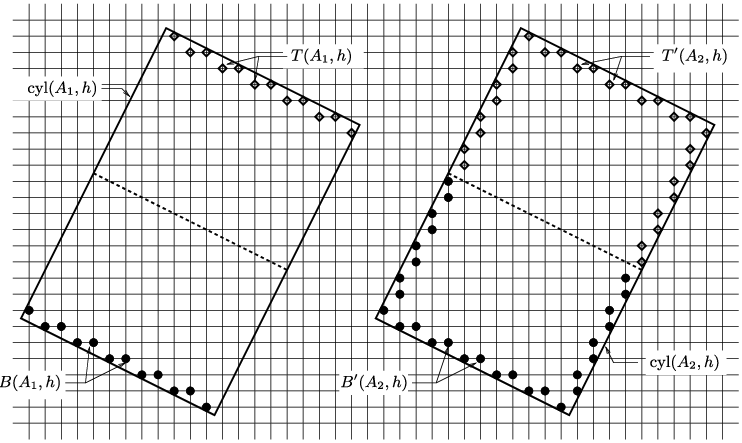}

\caption{The sets $T(A,h)$, $B(A,h)$, $T'(A,h)$ and $B'(A,h)$ in $\cyl(A,h)$.}
\label{figTB}
\end{figure}

We start with the definitions of flows in cylinders that will be useful
during the proof of Theorem \ref{thmpcpal} and the rigorous definition
of the function $\nu$ that appeared in equation (\ref{eqdefnu}). Let
$A$ be a nondegenerate hyperrectangle,
that is, a box of dimension $d-1$ in $\mathbb{R}^d$. All
hyperrectangles are
supposed to be closed in $\mathbb{R}^d$. We denote by
$\vv$ one of the two unit vectors orthogonal to $\hyp(A)$. For $h$ a
positive real number, we consider the cylinder $\cyl(A,h)$. Let
$T(A,h)$ [resp., $B(A,h)$] be the top (resp., the bottom) of $\cyl(A,h)$
with regard to the direction $\vv$ (see Figure~\ref{figTB}), that is,
\[
T(A,h) = \bigl\{ x\in\cyl(A,h) | \exists y\notin\cyl(A,h), [x,y] \in
\mathbb{E}_n^d \mbox{ and } [x,y] \cap(A+h\vv) \neq
\varnothing \bigr\}
\]
and
\[
B(A,h) = \bigl\{ x\in\cyl(A,h) | \exists y\notin\cyl(A,h), [x,y] \in
\mathbb{E}_n^d \mbox{ and }[x,y] \cap(A-h\vv)\neq
\varnothing \bigr\}.
\]
%

%
%
Let $T'(A,h)$ [resp., $B'(A,h)$] be the upper half part (resp., the lower
half part) of the boundary of $\cyl(A,h)$ (see Figure~\ref{figTB});
that is, if we denote by $z$ the center of~$A$,
\begin{eqnarray*}
&&\hspace*{-4pt}T'(A,h)\\
&&\hspace*{-4pt}\qquad = \lleft\{ x\in\cyl(A,h) \left\vert %
\begin{array} {l} \overrightarrow{zx} \cdot\vv>0 \mbox{ and}
\\
\exists y \notin\cyl(A,h), [x,y] \in\EE^d_n \mbox{ and
} [x,y]\cap\partial\cyl(A,h) \neq\varnothing \end{array} %
\right.
\rright
\}
\end{eqnarray*}
and
\begin{eqnarray*}
\hspace*{-4pt}&&B'(A,h) \\
&&\hspace*{-6pt}\qquad= \lleft\{ x\in\cyl(A,h) \left\vert %
\begin{array} {l} \overrightarrow{zx} \cdot\vv<0 \mbox{ and}
\\
\exists y \notin\cyl(A,h), [x,y] \in\EE^d_n \mbox{ and
} [x,y]\cap\partial\cyl(A,h) \neq\varnothing \end{array}
\right.
\rright
\}.
\end{eqnarray*}
For a given realization $(t(e),e\in\mathbb{E}_n^{d})$, we define the variable
$\tau_n (A,h) =\break   \tau_n(\cyl(A,h), \vv)$ by
\[
\tau_n(A,h) = \tau_n \bigl(\cyl(A,h), \vv \bigr) =
\phi_n \bigl(T'(A,h), B'(A,h), \cyl(A,h)
\bigr).
\]
The asymptotic behavior for large $n$ of the variable $\tau_n (A,h)$
properly rescaled
is well known, thanks to the almost subadditivity of
this variable. The following law of large numbers is proved in \cite
{RossignolTheret08b}:
%
\begin{thmm}[(Rossignol and Th\'eret)]
\label{thmnu}
We suppose that
\[
\int_{[0,+\infty[} x \,d\Lambda(x) < \infty.
\]
Then for each unit vector $\vv$ there exists a constant $\nu(d,
\Lambda, \vv) = \nu(\vv)$ (the dependence on $d$ and $\Lambda$ is
implicit) such that for every non degenerate hyperrectangle $A$
orthogonal to $\vv$ and for every strictly positive
constant $h$, we have\looseness=-1
\[
\lim_{n\rightarrow\infty} \frac{\tau_n(A,h)}{ n^{d-1}
\H^{d-1}(A)} = \nu(\vv) \qquad\mbox{in }
L^1.
\]\looseness=0
Moreover, if the origin of the graph belongs to $A$, or if
\[
\int_{[0,+\infty[} x^{1+{1}/{(d-1)}} \,d\Lambda(x) < \infty,
\]
then
\[
\lim_{n\rightarrow\infty} \frac{\tau_n(A,h)}{ n^{d-1}
\H^{d-1}(A)} = \nu(\vv) \qquad\mbox{a.s.}
\]
\end{thmm}
We emphasize the fact that the limit $\nu(\vv)$ depends on the
direction of $\vv$, but
neither on $h$ nor on the hyperrectangle $A$ itself. When the
capacities of the edges are bounded [hypothesis \ref{hypo2}], both
$L^1$ and a.s. convergences hold in Theorem~\ref{thmnu}. This theorem
gives the proper definition of the function $\nu$ that appeared in
equation~(\ref{eqdefnu}). The function $\nu$ is initially defined on
$\SS^{d-1}$, but we consider its homogeneous extension to $\RR^d$, that
we still denote by $\nu$, defined by
\[
\nu(\vec0) = 0 \quad\mbox{and}\quad \forall\vec w \in\RR ^d\setminus\{
\vec0 \}\qquad \nu(\vec w ) = \| w\|_2 \nu \biggl( \frac{\vec w}{ \| \vec w \|_2}
\biggr).
\]
We recall some geometric properties of the map $\nu$ that are valid
whenever $\EE(t(e))<\infty$. They have been stated in Section~4.4 of
\cite{RossignolTheret08b}. If there exists a unit vector $\vv$ such
that $\nu(\vv)=0$, then $\nu=0$ everywhere, and this happens if and
only if $\Lambda(\{0\}) \geq1-p_c(d)$, where $p_c(d)$ denotes the
critical parameter for bond percolation on $\ZZ^d$. This property has
been proved by Zhang in \cite{Zhang}. Moreover, the function $\nu:
\RR
^d \rightarrow\RR$ is convex. Since $\nu$ is finite, this implies that
$\nu$ is continuous on $\RR^d$. Moreover, $\nu$ is invariant under any
transformation of $\RR^d$ that preserves the graph $(\ZZ^d,\EE^d)$, in
particular $\nu(\vv) = \nu(-\vv)$ for all $\vv\in\RR^d$.

\begin{figure}

\includegraphics{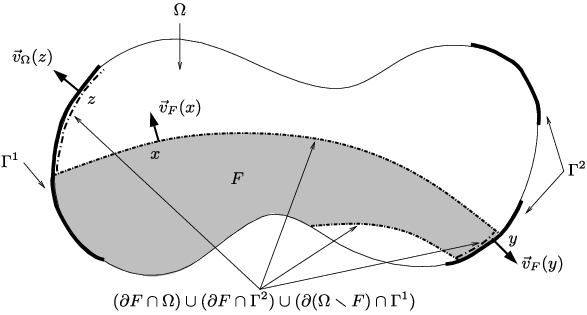}

\caption{The set $(\partial F \cap\O) \cup(\partial F \cap\G^2)
\cup(\partial(\O
\setminus F) \cap\G^1)$.}
\label{figmincut1}
\end{figure}

The asymptotic behavior of the maximal flow $\phi_n (\G^1_n,\G
^2_n,\O
_n)$ was studied in the companion papers \cite{CerfTheret09supc,CerfTheret09geoc} and \cite{CerfTheret09infc}, and the following law
of large numbers was proved:
%
\begin{thmm}[(Cerf and Th\'eret)]
\label{thmLGNphi}
We suppose that the hypotheses \ref{hypo1} and~\ref{hypo2} are
fulfilled. Then there exists a finite constant $\phi_{\O} \geq0$
defined in (\ref{eqdefcte1}) and (\ref{eqdefcte2}) such that
\[
\lim_{n\rightarrow\infty} \frac{\phi_n}{n^{d-1}} = \phi_{\O
}^{(1)}\qquad
\mbox{a.s.}
\]
Moreover, this equivalence holds:
\[
\phi_{\O}^{(1)} > 0 \quad\iff\quad\Lambda \bigl(\{ 0 \} \bigr) <
1-p_c(d).
\]
\end{thmm}
In fact the authors prove in \cite{CerfTheret09infc} that the lower
large deviations of $\phi_n /n^{d-1}$ below a constant $\phi_{\O
}^{(1)}$ are of surface order, in \cite{CerfTheret09supc} that the
upper large deviations of $\phi_n /n^{d-1}$ above a constant $\phi
_{\O
}^{(2)}$ are of volume order and finally in \cite{CerfTheret09geoc}
that $\phi_{\O}^{(1)} = \phi_{\O}^{(2)}$. The definitions of $\phi
_{\O
}^{(1)}$ and $\phi_{\O}^{(2)}$ are the following:
%
%
\begin{eqnarray}
\label{eqdefcte1} \phi_{\O}^{(1)}& =& \inf\left\{
\begin{array} {l} \displaystyle\int_{\O\cap\partial^* F} \nu \bigl(
\vv_F (x) \bigr) \,d \H^{d-1} (x) + \int
_{\G^2 \cap\partial^* F} \nu \bigl(\vv_F (x) \bigr) \,d
\H^{d-1} (x)
\\
\displaystyle\qquad{}+ \int_{\G^1 \cap\partial^* (\O\setminus F)} \nu \bigl(\vv_{\O} (x) \bigr) \,d \H
^{d-1} (x) \end{array}\right.
\left\vert\hspace*{-200pt}\phantom{\begin{array} {l}\displaystyle\int_{\O\cap\partial^* F}\\
\displaystyle\qquad{}+ \int_{\G^1 \cap\partial^* (\O\setminus F)} \end{array}}\right.
\nonumber
\\[-10pt]
\\[-10pt]
\nonumber
&&\hspace*{106pt}\left. \phantom{\begin{array} {l}\displaystyle\int_{\O\cap\partial^* F}\\
\displaystyle\qquad{}+ \int_{\G^1 \cap\partial^* (\O\setminus F)} \end{array}}F \subset\O,
\mathbh{1}_F \in \operatorname{BV}(\O)\right \},
\\[-4pt]
\label{eqdefcte2}\qquad \phi_{\O}^{(2)} &=& \inf\Biggl\{ \int
_{\O\cap\partial P} \nu \bigl(\vv _P (x) \bigr) \,d
\H^{d-1}(x) \Big\vert
\nonumber
\\[-10pt]
\\[-10pt]
\nonumber
&&\hspace*{18pt}\begin{array} {l} P\subset\,\RR^d
, \overline{\G}^1 \subset\,\stackrel {\circ} {P}, \overline{
\G}^2 \subset\,\stackrel{\circ} {\wideparen {\RR ^d
\setminus P}}
\\
P \mbox{ is polyhedral}, \partial P \mbox{ is transversal to }\G \end{array}
\Biggr\}.
\end{eqnarray}
The variational problems $\phi_{\O}^{(1)}$ and $\phi_{\O}^{(2)}$ are
continuous min-cut problems very similar to the problem $\phi_{\O
}^{(m)}$ defined by Nozawa. The variational problem $\phi_{\O}^{(1)}$
is in fact exactly the one we were looking for, that is, $\phi_{\O}^a =
\phi_{\O}^{(1)}$, where $\phi_{\O}^a$ is the continuous min-cut problem
appearing in Sections~\ref{secdef2} and \ref{secdef3}. Notice that a
condition of the type ``$\partial F$ separates $\G^1$ from $\G^2$ in
$\O$''
does not appear in $\phi_{\O}^{(1)}$, but the definition of the
capacity of $F$ is adapted: the surface that is considered as
``separating'' is in fact the surface $(\partial F \cap\O) \cup
(\partial F \cap\G
^2) \cup(\partial(\O\setminus F) \cap\G^1)$ (see Figure~\ref{figmincut1}). Thus we define for every $F \subset\O$ such that
$\mathbh{1}_F \in \operatorname{BV} (\O)$,
\begin{eqnarray*}
\capa(F) & =& \int_{\O\cap\partial^* F} \nu \bigl(\vv_F (x)
\bigr) \,d \H^{d-1} (x) + \int_{\G^2 \cap\partial^* F} \nu \bigl(
\vv_F (x) \bigr) \,d \H^{d-1} (x)
\\
& &{}+ \int_{\G^1 \cap\partial^* (\O\setminus F)} \nu \bigl(\vv_{\O
} (x) \bigr) \,d
\H^{d-1} (x),
\end{eqnarray*}
and the variational problem $\phi^{(1)}$ can be rewritten as
\[
\phi_{\O}^{a} = \phi_{\O}^{(1)} = \inf
\bigl\{\capa(F) | F \subset\O, \mathbh{1}_F \in \operatorname{BV}(\O) \bigr\}.
\]

%
%
Thus the set $\Sigma^a$ appearing in Theorem \ref{thmpcpal2} is
defined by
\[
\Sigma^a = \bigl\{ F \subset\O| \mathbh{1}_F \in \operatorname{BV}(\O)
, \capa(F) = \phi_{\O}^a \bigr\}.
\]

Let us prove that the min-cut problems $\phi_{\O}^{(1)}$, $\phi_{\O
}^{(2)}$ and $\phi_{\O}^{(m)}$ are equivalent. We claim that
%
%
\begin{equation}
\label{eqegalcte} \phi_{\O}^{(1)} \leq\phi_{\O}^{(m)}
\leq\phi_{\O}^{(2)}.
\end{equation}
Since $\phi_{\O}^{(1)}=\phi_{\O}^{(2)}$ by \cite{CerfTheret09geoc},
Theorem 11, we conclude that
\[
\phi_{\O}^{(1)} = \phi_{\O}^{(2)} =
\phi_{\O}^{(m)}.
\]
%
\begin{figure}

\includegraphics{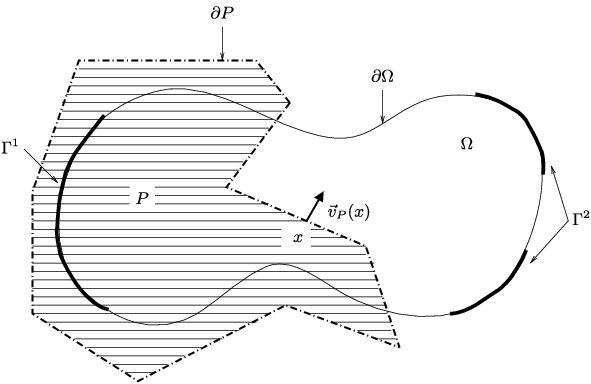}

\caption{A polyhedral set $P$ as in the definition of $\phi_{\O}^{(2)}$.}
\label{figP}
\end{figure}
Thus the three min-cut problems are equivalent. The arguments to
justify inequality (\ref{eqegalcte}) are the following. On one hand,
consider a set $P$ as in the definition~(\ref{eqdefcte2}) of $\phi
_{\O
}^{(2)}$ (see Figure~\ref{figP}), and define $S=P\cap\O$.
%
%
Since $\overline{\G}^1 \subset\stackrel{\circ}{P}$, then $\gamma
(\mathbh{1}_S) = 1$ $\H^{d-1}$-a.e. on $\G^1$ and since $\overline
{\G
}^2 \subset\,\stackrel{\circ}{\wideparen{\RR^d\setminus P}}$, then
$\gamma(\mathbh{1}_S) =0$ $\H^{d-1}$-a.e. on $\G^2$, thus $S$
satisfies all the conditions required in the definition (\ref
{eqmincutNozawa}) of $\phi_{\O}^{(m)}$ and
\[
\int_{\O\cap\partial^* S} \nu \bigl(\vv_S (x) \bigr) \,d
\H^{d-1} (x) = \int_{\O
\cap\partial P} \nu \bigl(
\vv_P (x) \bigr) \,d\H^{d-1} (x).
\]
Thus $\phi_{\O}^{(m)} \leq\phi_{\O}^{(2)}$. On the other hand consider
a set $S$ as in the definition (\ref{eqmincutNozawa}) of $\phi_{\O
}^{(m)}$. Of course $S$ satisfies the conditions required in the
definition (\ref{eqdefcte1}) of $\phi_{\O}^{(1)}$. According to the
last equality on page 809 in \cite{Nozawa}, for every set $S\subset\O$
of finite perimeter in $\O$ we have
%
%
\begin{equation}
\label{eqNozawatrace} \gamma(\mathbh{1}_S) = \mathbh{1}_{\G\cap\partial^*
S},\qquad
\H ^{d-1}\mbox{-a.e. on } \G.
\end{equation}
Thus $\gamma(\mathbh{1}_S) = 0$ $\H^{d-1}$-a.e. on $\G^2$ implies
that $\H^{d-1} (\G^2 \cap\partial^* S) =0$. By definition of the trace,
$\gamma(\mathbh{1}_{\O\setminus S}) = 1 - \gamma(\mathbh{1}_S)$
everywhere on $\G$, thus $\gamma(\mathbh{1}_S) = 1$ $\H^{d-1}$-a.e. on
$\G^1$ implies $\gamma(\mathbh{1}_{\O\setminus S})= 0$ $\H
^{d-1}$-a.e. on $\G^1$. Since $\O\setminus S$ has also finite
perimeter in $\O$, by equation (\ref{eqNozawatrace}) applied to $\O
\setminus S$ we have $\gamma(\mathbh{1}_{\O\setminus S}) =
\mathbh{1}_{\G\cap\partial^* (\O\setminus S)}$ $\H
^{d-1}$-a.e. on $\G
$, thus $\H^{d-1} (\G^1 \cap\partial^* (\O\setminus S)) =0$,
and the integrals
\[
\int_{\G^2 \cap\partial^* F} \nu \bigl(\vv_F (x) \bigr) \,d
\H^{d-1} (x)\quad \mbox {and} \quad\int_{\G^1 \cap\partial^* (\O\setminus F)} \nu \bigl(\vv
_{\O} (x) \bigr) \,d \H^{d-1} (x)
\]
vanish. We conclude that $\phi_{\O}^{(1)} \leq\phi_{\O}^{(m)}$, and
this finishes the proof of inequality~(\ref{eqegalcte}).

\begin{rem}
The simplicity of the previous argument should not hide that the real
difficulty consists in proving that $\phi_{\O}^{(1)} = \phi_{\O
}^{(2)}$. This is done in \cite{CerfTheret09geoc} by a quite
complicated process of polyhedral approximation.
\end{rem}


\section{Organization of the proof}

In Section~\ref{secmaxstream}, we study a sequence of discrete maximal
streams $(\vec\mu_n^{\max})_{n \geq1}$. We prove that from each
subsequence of $(\vec\mu_n^{\max})_{n \geq1}$ we can extract a
sub-subsequence which is weakly convergent. If we denote by $\vec\mu$
its limit, we prove that a.s. $\vec\mu= \vec\sigma\L^d$ with
$\vec
\sigma$ a continuous stream which is admissible for the max-flow
problem $\phi_{\O}^b$. Moreover, we prove that along the converging
subsequence,
%
%
\begin{equation}
\label{eqsketch1} \lim_{n\rightarrow\infty} \operatorname{flow}^{\mathrm
{disc}}_n
\bigl(\vec\mu_n^{\max} \bigr) = \operatorname{flow}^{\mathrm
{cont}}(
\vec \sigma) \qquad\mbox{a.s.}
\end{equation}
Section~\ref{secmincut} is devoted to the study of a sequence of
minimal cutsets $(\E_n^{\min})_{n\geq1}$. We prove that from each
subsequence of $(R(\E_n^{\min}))_{n \geq1}$ we can extract a
sub-subsequence which is convergent for the distance $\mathfrak d$. If
we denote by $F$ its limit, we prove that $F \subset\O$ and $\mathbh
{1}_{F} \in \operatorname{BV} (\O)$; that is, $F$ is admissible for the min-cut
problem $\phi_{\O}^a$. Moreover, we prove that along the converging
subsequence,
%
%
\begin{equation}
\label{eqsketch2} \liminf_{n\rightarrow\infty} \frac{V(\E_n^{\min})}{n^{d-1}} \geq
\capa(F)\qquad \mbox{a.s.}
\end{equation}
In Section~\ref{secccl} we establish that
%
%
\begin{equation}
\label{eqsketch3} \capa(F) \geq\operatorname{flow}^{\mathrm{cont}}(\vec\sigma).
\end{equation}
Then combining equations (\ref{eqsketch1}), (\ref{eqsketch2}) and
(\ref
{eqsketch3}) we derive the results presented in Section~\ref{secdef3}.

The most original part of our work is the study of maximal streams
presented in Section~\ref{secmaxstream}. The study of minimal cutsets
relies largely on the techniques used in~\cite{CerfTheret09infc} to
prove that the lower large deviations of $\phi_n$ are of surface order.
To complete the proofs we also use the result of polyhedral
approximation proved in \cite{CerfTheret09geoc}. In the proof of the
law of large numbers for $\phi_n$ we present here, we have replaced the
study of the upper large deviations of $\phi_n$ performed in \cite
{CerfTheret09supc} by the study of the maximal streams, which is more
natural, and we have adapted the arguments given in the study of the
lower large deviations of $\phi_n$ in \cite{CerfTheret09infc} to obtain
informations on the behavior of minimal cutsets.

Throughout the paper, we assume that hypotheses \ref{hypo1} and \ref{hypo2} are satisfied.


\section{Study of maximal streams}
\label{secmaxstream}

\subsection{Existence}
\label{secexistencecourant}

\begin{figure}

\includegraphics{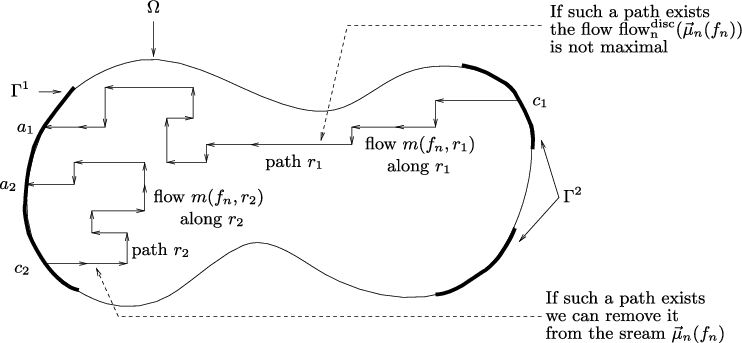}

\caption{Flow that escapes from $\O_n$ through $\G^1_n$.}
\label{figcourantmax}
\end{figure}

The existence of at least one maximal stream is not so obvious because
of condition (\ref{fluxmax2}). We will assume throughout the paper that
the capacities of the edges are bounded by a constant $M$. Under this
hypothesis, the set $\S_n (\G^1_n, \G^2_n, \O_n)$ is compact, and since
the function $f_n \in\S_n (\G^1_n, \G^2_n, \O_n) \mapsto
\operatorname{flow}^{\mathrm{disc}}_n
(\vec\mu_n
(f_n))$ is continuous, a stream $\vec\mu_n = \vec\mu_n (f_n)$ satisfying
(\ref{fluxmax1}) exists. Suppose that $\vec\mu_n$ does not satisfy
(\ref
{fluxmax2}), and let $e= [a,b]$ with $a\in\G^1_n$, $b\notin\G^1_n$,
and, for example, $\be= \langle a,b\rangle$ and $f_n (e) < 0$. Since
$f_n$ satisfies the node law and since there exists only a finite
number of self avoiding paths (i.e., paths that visit each edge at most
once) starting at $a$ in $\O_n$, then there exists a self avoiding path
$r= (a, [a,b],b,\ldots,c)$ in $\O_n$ from $a$ to a point $c$ that belongs
to $\G^1_n$ or $\G^2_n$ such that for all $e\in r$, if $r(e) =1$
(resp.,
$-1$) when $e$ is crossed by $r$ from the origin to the endpoint of
$\be
$ (resp., from the endpoint to the origin of $\be$), then $f_n(e) r(e)
<0$. Since $r$ is finite,
\[
m(f_n, r) = \inf \bigl\{- f_n(e) r(e) | e \in r \bigr
\} > 0.
\]
Consider the stream function $f_n'$ defined by
\[
f_n'(e) = \cases{ %
f_n(e),& \quad$\mbox{if }e\notin r,$
\vspace*{2pt}\cr
f_n(e) + r(e) m(f_n,r),&\quad $\mbox{if } e\in r.$}
\]
This is the stream function obtained by removing from $f_n$ a quantity
of flow $m(f_n, r)$ along $r$ from $c$ to $a$. The stream function
$f_n'$ is still admissible, since $|f_n' (e) | \leq| f_n(e)|$ for all
$e$. If $c$ belongs to $\G^2_n$ (see Figure~\ref{figcourantmax}), then
$\operatorname{flow}^{\mathrm{disc}}_n(\vec\mu_n (f_n')) =
\operatorname{flow}^{\mathrm{disc}}_n(\vec\mu_n (f_n)) + m(f_n,r)$, and
this is
not possible since $\vec\mu_n (f_n)$ satisfies (\ref{fluxmax1}). Thus
$c$ belongs to $\G^1_n$ (see Figure~\ref{figcourantmax}) and
$\operatorname{flow}^{\mathrm{disc}}_n
(\vec
\mu_n (f_n')) = \operatorname{flow}^{\mathrm{disc}}_n(\vec\mu_n
(f_n)) $. Moreover, $f_n' ([a,b]) = f_n
([a,b]) + m(f_n,r) > f_n ([a,b])$, and $m(f_n', r)=0$. We can iterate
this process finitely many times with every possible self avoiding path
$r'$ starting at $a$ until $m(f_n,r') =0$ for all $r'$. Eventually, the
stream function $f_n''$ we obtain satisfies $f_n''([a,b]) =0$. We can
do the same procedure with every edge $[a,b]$ with $a\in\G^1_n$ and
$b\notin\G^1_n$ (there is a finite number of such edges), and the
stream function $\widetilde f_n$ that we obtain at the end satisfies
(\ref{fluxmax2}) and $\operatorname{flow}^{\mathrm{disc}}_n(\vec\mu
_n (\widetilde f_n)) = \operatorname{flow}^{\mathrm{disc}}_n(\vec
\mu_n
(f_n))$. This proves the existence of a maximal stream from $\G^1_n$ to
$\G^2_n$ in $\O_n$ if we suppose that the capacities of the edges are
bounded by a constant $M$.
%

From now on $(\vec\mu_n)_{n\geq1}$ denotes a sequence of admissible
discrete streams and $(\vec\mu_n^{\max})_{n\geq1}$ a sequence of
admissible maximal discrete streams.


\subsection{Compactness}
\label{seccomp}

We prove the following property:
%
\begin{prop}
\label{propcomp} Almost surely, for $n$ large enough, the
sequence\break
$(\vec
\mu_n)_{n\geq1}$ takes its values in a deterministic weakly compact
set of measures.
\end{prop}
%
\begin{rem}
This property implies that any subsequence of $(\vec\mu_n)_{n\geq1}$
admits a sub-subsequence $(\vec\mu_{\varphi(n)})_{n\geq1}$ that is
weakly convergent, that is, such that there exists a random vector
measure $\vec\mu\dvtx  \B(\RR^d) \rightarrow\RR^d $ satisfying
\[
\forall f \in\C_b \bigl(\RR^d, \RR \bigr)\qquad \lim
_{n\rightarrow\infty} \int_{\RR^d} f \,d\vec
\mu_{\varphi(n)} = \int_{\RR^d} f \,d\vec\mu.
\]
The choice of the sub-subsequence $(\vec\mu_{\varphi(n)})_{n\geq
1}$ is
random, that is, the function $\varphi$ may depend on the realization
of the capacities.
\end{rem}

\begin{pf*}{Proof of Proposition \ref{propcomp}} For the rest of
this section, we consider a fixed realization of the capacities. Let
$\vec\mu_n = (\mu_n^1,\ldots, \mu_n^d)$ be an admissible discrete
stream on
$(\ZZ^d_n, \EE^d_n)$. For all $n\geq1$, the support of $\vec\mu
_n$ is
included in the compact set $\overline{\V_\infty(\O, 1)}$. Hence the
admissible discrete streams are tight. Moreover, for $i =1,\ldots,d$,
we have
\[
\bigl| \mu^i_n \bigr| \bigl(\overline{\V_\infty(\O, 1)}
\bigr) \leq \frac
{1}{n^d} \sum_{e \in\Pi_n}
\bigl|f_n(e)\bigr| \leq\frac{M \card(\Pi
_n)}{n^d}.
\]
Since
%
%
\begin{eqnarray}
\label{eqPin} \card(\Pi_n) & \leq&2d \card(\O_n) \leq2d
n^d \L^d \biggl(\O+ \frac{1}{2n} [-1, 1] \biggr)\nonumber\\
&\leq&2d n^d \L^d \bigl(\V _\infty \bigl(\O,
n^{-1} \bigr) \bigr)
\\
& \leq&2d n^d \L^d \bigl(\V_\infty(\O, 1 ) \bigr)
,\nonumber
\end{eqnarray}
we conclude that for all $i=1,\ldots,d$, for all $n\geq1$,
\[
\bigl| \mu^i_n \bigr| \bigl(\overline{\V_\infty(\O, 1)}
\bigr) \leq2d M \L^d \bigl(\V_\infty(\O, 1) \bigr).
\]
Thus the admissible discrete streams are uniformly bounded for the
total variation distance. The conclusion follows from Prohorov's
theorem; see, for example, Theorem 8.6.2 in volume II of \cite{Bogachev}.
\end{pf*}

\begin{rem}
Since all the measures $\vec\mu_n$ have a support included in the same
compact, the weak convergence $(\vec\mu_n)_{n\geq1}$ is characterized
by the convergence of $\vec\mu_n (f)$ for all $f$ in any of the
following classes of functions: the continuous bounded functions, the
continuous functions with compact support or the continuous functions
that goes to zero at infinity.
\end{rem}

From now on, we consider a measure $\vec{\mu}$ which is the weak
limit of
a subsequence of $(\vec{\mu}_n)_{n\geq1}$, and we study its properties.
Notice that $\vec{\mu}$ is a priori random, so some of its properties
will be proved for all events, and others only a.s.


\subsection
{Absolute continuity with respect to
Lebesgue measure}
\label{secAC}

In this section, we prove that $\vec{\mu}$ is absolutely continuous with
respect to $\L^d$, the Lebesgue measure on~$\RR^d$.
%
\begin{prop}
\label{propAC}
If $\vec\mu$ is the weak limit of a subsequence of $(\vec\mu
_n)_{n\geq
1}$, where $\vec\mu_n$ is an admissible stream for all $n\geq1$, then
there exists a random vector field $\vec{\sigma} \dvtx \RR^d
\rightarrow\RR
^d$ such that $\vec{\mu} = \vec{\sigma}\L^d$, $\vec{\sigma}
\in L^{\infty}
(\RR^d\rightarrow\RR^d, \L^d)$ and $\vec\sigma= 0$ $\L^d$-a.e.
on~$\O^c$.
\end{prop}

\begin{pf}
For the rest of this
section, we consider a fixed realization of the capacities. Let $\mu
_n^i = \mu_n^{i,+} - \mu_n^{i,-}$ be the Hahn--Jordan decomposition of
the signed measure $\mu_n^i$. Then $\mu_n^{i,+}$ and $\mu_n^{i,-}$ are
positive measures on $(\RR^d, \B(\RR^d) )$, respectively, the positive
and negative part of $\mu_n^i$. By the same arguments as in Proposition
\ref{propcomp}, we see that the sequences $(\mu_n^{i,+})_{n\geq1}$ and
$(\mu_n^{i,-})_{n\geq1}$ take their value in a weakly compact set.
Thus up to extraction we can suppose that $\vec{\mu}_n
\rightharpoonup\vec
{\mu}$, $\mu_n^{i,+} \rightharpoonup\mu^{i,+}$ and $\mu_n^{i,-}
\rightharpoonup\mu^{i,-}$ for all $i=1,\ldots,d$, where $\mu^{i,+}$ and
$\mu^{i,-}$ are positive measures. If we write $\vec{\mu} = (\mu
^1,\ldots,\mu
^d)$, we have $\mu^i = \mu^{i,+} - \mu^{i,-}$, but this may not be the
Hahn--Jordan decomposition of $\mu^i$ since it may not be minimal. Let
$B(x,r)$ be the ball centered at $x$ of radius $r>0$. We have
\begin{eqnarray*}
\mu_n^{i,+} \bigl(B(x,r) \bigr) &\leq&\bigl|
\mu_n^i\bigr | \bigl(B(x,r) \bigr) \leq\frac
{1}{n^d} \sum
_{e \in\EE^{d,i}_n, c(e) \in B(x,r)}\bigl |f_n(e)\bigr|\\
& \leq& \frac{M \card(\{ e \in\EE^{d,i}_n | c(e) \in B(x,r) \})}{n^d},
\end{eqnarray*}
and we remark as in Section~\ref{seccomp} that
\begin{eqnarray*}
\card \bigl( \bigl\{ e \in\EE^{d,i}_n | c(e) \in B(x,r)
\bigr\} \bigr)& \leq &\card \bigl( \ZZ^d_n \cap B
\bigl(x,r + n^{-1} \bigr) \bigr)
\\
& \leq& n^d \L^d \biggl( \ZZ^d_n
\cap B \bigl(x, r+n^{-1} \bigr) + \frac
{1}{2n}
[-1,1]^d \biggr)
\\
& \leq & n^d \L^d \bigl(B \bigl(x,r+ 2 n^{-1}
\bigr) \bigr),
\end{eqnarray*}
whence
\[
\mu_n^{i,+} \bigl(B(x,r) \bigr) \leq M \L^d
\bigl(B \bigl(x, r+2 n^{-1} \bigr) \bigr).
\]
With the help of Portmanteau's theorem (see, e.g., Theorem 8.2.3 in
\cite{Bogachev}) we obtain that
\[
\mu^{i,+} \bigl(B(x,r) \bigr) \leq M \L^d \bigl(B(x,r)
\bigr).
\]
Let next $A$ be a Borel subset of $\RR^d$. Since the Lebesgue measure
$\L^d$ is outer regular, for $\eps>0$ there exists an open set $O$
such taht $A \subset O$ and $\L^d (O \setminus A) < \eps$. By the
Vitali covering theorem for Radon measures (see Theorem 2.8 in \cite
{Mattila}), there exists a countable family $(B_j, j \in J)$ of
disjoint closed balls such that:
\begin{itemize}
\item$\forall j \in J\ B_j \subset O$;
\item$\mu^{i,+} (O \setminus\bigcup_{j\in J} B_j) = 0$.
\end{itemize}
Thus
\[
\mu^{i,+} (O) = \sum_{j \in J}
\mu^{i,+} (B_j) \leq M \sum_{j\in J}
\L^d(B_j) = M \L^d \biggl(\bigcup
_{j\in J} B_j\biggr) \leq M \L^d (O),
\]
whence
\[
\mu^{i,+} (A) \leq\mu^{i,+} (O) \leq M \L^d (O)
\leq M \bigl( \L^d (A) + \eps \bigr).
\]
Sending $\eps$ to $0$, we obtain that
\[
\mu^{i,+} (A) \leq M \L^d (A).
\]
We conclude that $\mu^{i,+}$ is absolutely continuous with respect to
$\L^d$. The same holds for $\mu^{i,-}$, for all $i=1,\ldots,d$, thus
$\vec
{\mu}$ is absolutely continuous with respect to $\L^d$; that is, there
exists $\vec{\sigma} \in L^1 (\RR^d\rightarrow\RR^d, \L^d)$ such
that $\vec
{\mu} = \vec{\sigma} \L^d$. We use the notation $\vec{\sigma} =
(\sigma
^1,\ldots, \sigma^d)$. Moreover we have proved that for all
$i=1,\ldots,d$,
\[
\forall A\in\B \bigl(\RR^d \bigr) \qquad\int_A \bigl|
\sigma^i\bigr| \,d \L^d \leq\mu ^{i,+} (A) +
\mu^{i,-} (A) \leq2M \L^d (A),
\]
which implies that $|\sigma^i| \leq2M $ $\L^d$-a.e. and thus that
$\vec
{\sigma}$ belongs to $L^\infty(\RR^d \rightarrow\RR^d, \L^d)$.
Finally, we notice that for $n \geq1$, the support of $\vec\mu_n$ is
included in $\V_\infty(\O, 1/n)$, thus the support of $\vec\mu$ is
included in $\bigcap_{n\geq1} \V_\infty(\O, 1/n) = \overline{\O}$. This
implies that $\vec\sigma=0$ $\L^d$-a.e. on $\RR^d \setminus
\overline{\O}$, thus on $\O^c$ since $\L^d (\partial\O)=0$.
\end{pf}


\subsection{Divergence and boundary conditions}
\label{secdiv}

We study the divergence of $\vec{\sigma}$. We recall that divergence must
be understood in the distributional sense. By definition, for every
function $ h \in\C_c^\infty(\RR^d, \RR)$, we have
%
%
\begin{equation}
\label{eqdefdiv} \int_{\RR^d} h \di\vec{\sigma} \,d
\L^d = - \int_{\RR^d} \vec {\sigma} \cdot\vec{
\nabla} h \,d\L^d.
\end{equation}
We first prove the following result:
%
\begin{prop}
\label{propdivgen}
If $\vec\mu= \vec\sigma\L^d$ is the weak limit of a subsequence
$(\vec
\mu_{\varphi(n)})_{n\geq1}$ of $(\vec\mu_n)_{n\geq1}$, then for every
function $h\in\C_c^\infty(\RR^d, \RR)$ we have
%
%
\begin{eqnarray}
\label{eqresdiv} \int_{\RR^d} h \di\vec{\sigma} \,d
\L^d& =& - \int_{\RR^d} \vec {\sigma} \cdot\vec{
\nabla} h \,d\L^d
\nonumber
\\[-8pt]
\\[-8pt]
\nonumber
& =& \lim_{n\rightarrow\infty} \frac{1}{\varphi(n)^{d-1}}
\sum_{x \in\Gamma^1_{\varphi(n)} \cup
\Gamma^2_{\varphi(n)}} h(x) \widehat f_{\varphi(n)} (x),
\end{eqnarray}
where for all $x \in(\Gamma^1_n \cup\Gamma^2_n)$, $\widehat f _n (x)$
is the amount of water that appears at $x$ according to the stream $f_n$:
\[
\widehat f_n (x) = \sum_{e= \langle\cdot, x \rangle}
f_n(e) - \sum_{e= \langle x, \cdot\rangle} f_n(e)
.
\]
\end{prop}

\begin{pf}
The idea of the
proof is the following: we interpret $\di\vec\sigma$ as the limit
of a
discrete divergence, which we can control thanks to the node law
satisfied by the stream function $f_n$. We consider again a fixed
realization of the capacities. We consider a subsequence of $(\vec\mu
_n)_{n\geq1}$ converging toward $\vec\mu$, but we still denote this
subsequence by $(\vec\mu_n)_{n\geq1}$ to simplify the notation. Since
$\vec{\nabla} h \in\C_b(\RR^d, \RR^d)$, we see that
%
%
\begin{equation}
\label{eqlimdiv}\quad \int_{\RR^d} \vec{\sigma} \cdot\vec{\nabla} h
\,d\L^d = \lim_{n\rightarrow\infty} \int_{\RR^d}
\vec{\nabla} h \cdot d\vec {\mu}_n = \lim_{n\rightarrow\infty}
\frac{1}{n^d} \sum_{e \in\EE_n^d} f_n(e)
\vec{e}\cdot\vec{\nabla} h \bigl(c(e) \bigr).
\end{equation}
We study the sum appearing in the previous equality. Let $i \in\{
1,\ldots,d \}$. By Taylor's theorem, we know that for all
$x=(x_1,\ldots,x_d),y=(y_1,\ldots,y_d)\in\RR^d$ such that $x_j=y_j$
for all
$j\neq i$ we have
\[
h(x) - h(y) = \partial_i h (y) (x_i - y_i) +
g (x,y),
\]
and since $h$ is in $\C^2_c(\RR^d, \RR)$ we know that $
|g(x,y)| \leq K (x_i - y_i)^2$, where $K = K(h) = \| h \|_{W^{2,\infty
}} /2$. For $e\in\EE^{d,i}_n$, let us denote by $l_i(e)$ [resp.,
$r_i(e)$] the endpoint at the origin (resp., the end) of $\be$ according
to the orientation
\begin{figure}[b]

\includegraphics{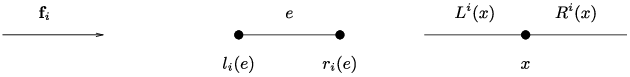}

\caption{Correspondence between edges and vertices.}
\label{figarete}
\end{figure}
chosen on $\EE^{d,i}_n$; see Figure~\ref{figarete}.
Conversely, for $x\in\ZZ^d_n$, we denote by $L^i(x)$ [resp., $R^i(x)$]
the edge of $\EE^{d,i}_n$ which ends at $x$ (resp., starts at
$x$).\break
%
%
We obtain
\begin{eqnarray*}
&&\sum_{e \in\EE^d_n} f_n(e) \vec{e}\cdot\vec{
\nabla} h \bigl(c(e) \bigr)\\
&&\qquad = \sum_{i=1}^d
\sum_{e \in\EE^{d,i}_n} f_n(e) \partial_i
h \bigl(c(e) \bigr)
\\
&&\qquad = \sum_{i=1}^d \sum
_{e \in\EE^{d,i}_n} f_n(e) n \bigl( \bigl[ h
\bigl(r_i(e) \bigr) - h \bigl(l_i(e) \bigr) \bigr] -
\bigl[ g \bigl(r_i(e), c(e) \bigr) - g \bigl(c(e), l_i(e)
\bigr) \bigr] \bigr)
\\
&&\qquad =  n \sum_{x \in\ZZ^d_n} h(x) \sum
_{i=1}^d \bigl[ f_n
\bigl(L_i(x) \bigr) - f_n \bigl(R_i(x) \bigr)
\bigr] + \alpha_n (h,f,\O),
\end{eqnarray*}
\noindent\hspace*{-4pt} where by inequality (\ref{eqPin}) we have
%
%
\begin{equation}
\label{eqalpha} \bigl|\alpha_n (h,f,\O)\bigr| \leq n \frac{2K}{(2 n)^2} M
\card(\Pi_n) \leq d K M \L^d \bigl(\V_\infty(\O,
1) \bigr) n^{d-1}.
\end{equation}
Since the stream satisfies the node law, we have for all $x \notin
(\Gamma^1_n \cup\Gamma^2_n)$ that
\[
\sum_{i=1}^d \bigl[ f_n
\bigl(L_i(x) \bigr) - f_n \bigl(R_i(x) \bigr)
\bigr] = 0.
\]
For all $x \in(\Gamma^1_n \cup\Gamma^2_n)$, let us denote by
$\widehat f _n (x)$ the amount of water that appears at $x$ according
to the stream $f_n$, that is,
%
%
\begin{equation}
\widehat f _n (x) = \sum_{i=1}^d
\bigl[ f_n \bigl(R_i(x) \bigr) - f_n
\bigl(L_i(x) \bigr) \bigr] = \sum_{e = \langle\cdot, x \rangle}
f_n(e) - \sum_{e =
\langle x, \cdot\rangle} f_n(e)
.
\end{equation}
Then we have proved that
\[
\int_{\RR^d} \vec{\nabla} h \cdot d\vec{\mu}_n =
- \frac{1}{n^{d-1}} \sum_{x \in\Gamma^1_n \cup\Gamma^2_n} h(x) \widehat
f_n(x) + \frac
{\alpha_n (h,f, \O)}{n^{d}}.
\]
According to equations (\ref{eqlimdiv}) and (\ref{eqalpha}), this
implies equation (\ref{eqresdiv}), and thus Proposition \ref
{propdivgen} is proved.
\end{pf}

We now deduce from Proposition \ref{propdivgen} that $\di\vec\sigma$
and $\vec{\sigma} \cdot\vv_{\O}$ satisfy the conditions required
in \cite
{Nozawa}. Remember that divergence is understood in the distributional
sense. The meaning of $\vec{\sigma} \cdot\vv_{\O}$ is the one
given by
Nozawa in \cite{Nozawa} that we have recalled in Section~\ref{secNozawa}.
%
\begin{cor}
\label{cordiv}
If $\vec\mu= \vec\sigma\L^d$ is the limit of a subsequence of
$(\vec\mu
_n)_{n\geq1}$, then it satisfies $\di\vec\sigma= 0 $ $\L^d$-a.e. on
$\O$ and $\vec\sigma\cdot\vv_{\O} = 0 $ $\H^{d-1}$-a.e. on $\G
\setminus(\G^1 \cup\G^2)$. Moreover, if for all $n\geq1$,
$\widehat f_n(x):= \sum_{e = \langle\cdot, x \rangle} f_n(e) -
\sum_{e = \langle x, \cdot\rangle} f_n(e) \geq0$ for all $x\in\G^1_n$
then $\vec\sigma\cdot\vv_{\O} \leq0 $ $\H^{d-1}$-a.e. on $\G^1$.
\end{cor}

\begin{rem}
By definition, the last condition is satisfied by a sequence of maximal
flows $(\vec\mu_n^{\max})_{n\geq1}$.
\end{rem}

\begin{pf*}{Proof of Corollary \ref{cordiv}} We consider a fixed
realization of the capacities. We prove first that $\di\vec{\sigma}
= 0$
on $\O$ in terms of distributions. Indeed, for every function $h\in\C
_c^\infty(\O, \RR)$, $h$ is null on $\G^1_n \cup\G^2_n$, for all $n$.
Thus by Proposition \ref{propdivgen},
\[
\int_{\RR^d} h \di\vec{\sigma} \,d \L^d = - \int
_{\RR^d} \vec {\sigma} \cdot\vec{\nabla} h \,d\L^d =
0.
\]
As explained in Remark \ref{remdiv1}, we rewrite this equality as $\di
\vec\sigma= 0$ $\L^d$-a.e. on $\O$.
We now study the boundary conditions satisfied by $\vec\sigma$. As
explained in Section~\ref{secNozawa}, $\vec{\sigma} \cdot\vv_{\O
}$ is an
element of $L^{\infty} (\G, \H^{d-1})$ characterized by
%
%
\begin{equation}
\label{eqdefbord} \forall u \in W^{1,1} (\O)\qquad \int
_{\G} (\vec{\sigma}\cdot\vv _{\O
} ) \gamma(u) \,d
\H^{d-1} = \int_{\O} \vec{\sigma} \cdot\vec {\nabla}
u \,d \L^d.
\end{equation}
In fact $\vec{\sigma} \cdot\vv_{\O}$ is characterized by
%
%
\begin{equation}
\label{eqdefbord2} \forall u \in\C_c^\infty \bigl(
\RR^d, \RR \bigr)\qquad \int_{\G} (\vec{\sigma }\cdot
\vv_{\O} ) u\, d \H^{d-1} = \int_{\RR^d} \vec{
\sigma} \cdot\vec{\nabla} u \,d \L^d.
\end{equation}
Let us prove that the conditions (\ref{eqdefbord}) and (\ref
{eqdefbord2}) are equivalent. We recall that $W^{1,1} (\O)$ is the set
of functions $u \dvtx \O\rightarrow\RR$ satisfying $u \in L^1(\O)$, and
for all $i\in\{ 1,\ldots,d \}$, there exists $g_i \in L^1 (\O)$ such that
\[
\forall h \in\C^\infty_c (\O, \RR) \qquad\int
_{\O} u \partial_i h \,d \L^d = - \int
_{\O} g_i h \,d \L^d.
\]
By definition $\| \partial_i u \|_{L^1(\O)} = \| g_i \|_{L^1(\O)}$.
The norm
on the Sobolev space $W^{1,1} (\O)$ is given by
\[
\forall u \in W^{1,1} (\O)\qquad \|u \|_{W^{1,1}(\O)} = \| u \|
_{L^1(\O)} + \sum_{i=1}^d \|
\partial_i u \|_{L^1(\O)}.
\]
The set of functions $\{ f|_{\O}, f \in\C_c^\infty(\RR^d, \RR)
\}$ is dense into $W^{1,1} (\O)$ with respect to the norm $\| \cdot\|
_{W^{1,1}(\O)}$ (see, e.g., \cite{Nozawa}, page 809). Let $u
\in W^{1,1} (\O) $ and $(u_n)_{n\geq1} $ be a sequence of functions in
$\C_c^\infty(\RR^d, \RR) $ such that $\tilde u_n = u_n|_{\O}$
converges toward $u$ in $W^{1,1 }(\O)$. Then $u$ and $\tilde u_n$
belong to $\operatorname{BV}(\O)$ for all $n\geq1$, and $\tilde u_n$ converges toward
$u$ in $\operatorname{BV}(\O)$ in the sense given by Nozawa in \cite{Nozawa}, page
808;\vspace*{1pt} that is, $\tilde u_n$ converges toward $u$ in $L^1(\O)$ and $|
\vec
\nabla\tilde u_n |(\O)$ converges toward $| \vec\nabla u | (\O)$. Then
by~\cite{Nozawa}, Theorem 2.1 (that comes from \cite{GI}), we know that
$\gamma( \tilde u_n) = u _n |_{\G} $ converges toward $\gamma(u)$ in
$L^1 (\G)$. Since $\vec\sigma\cdot\vv_{\O}$ is in $L^\infty(\G)$,
this implies that
%
%
\begin{equation}
\label{eqcarbord1} \lim_{n\rightarrow\infty} \int_{\G} (
\vec\sigma\cdot\vv_{\O
}) u_n\, d \H^{d-1} = \int
_{\G} (\vec\sigma\cdot\vv_{\O}) \gamma(u) \,d
\H^{d-1}.
\end{equation}
Moreover, since $\vec\sigma\in L^{\infty} (\O, \RR^d)$, the convergence
of $\tilde u _n$ toward $u$ in $W^{1,1} (\O)$ implies
%
%
\begin{equation}
\label{eqcarbord2} \lim_{n\rightarrow\infty} \int_{\O}
\vec\sigma\cdot\vec \nabla\tilde u_n \,d\L^d = \int
_{\O} \vec\sigma\cdot\vec\nabla u \,d\L^d.
\end{equation}
Finally, $\vec{\sigma} = 0$ $\L^d$-a.e. on $\O^c$ implies that
%
%
\begin{equation}
\label{eqcarbord3} \int_{\O} \vec\sigma\cdot\vec\nabla\tilde
u_n \,d\L^d = \int_{\O} \vec\sigma
\cdot\vec\nabla u_n \,d\L^d = \int_{\RR^d}
\vec \sigma \cdot\vec\nabla u_n \,d\L^d.
\end{equation}
Combining equations (\ref{eqcarbord1}), (\ref{eqcarbord2}) and (\ref
{eqcarbord3}), we conclude that properties (\ref{eqdefbord}) and (\ref
{eqdefbord2}) are equivalent. According to (\ref{eqresdiv}), we obtain
that for all $u \in\C_c^\infty(\RR^d, \RR)$,
\[
\int_\Gamma(\vec{\sigma} \cdot\vv_{\O}) u\, d
\H^{d-1} = \lim_{n\rightarrow\infty} - \frac{1}{n^{d-1}} \sum
_{x \in\G_n^1 \cup
\G
^2_n} u(x) \widehat f_n(x).
\]
On one hand, let $u$ be a function of $ \in\C_c^\infty((\G^1 \cup
\G
^2)^c, \RR)$, that is, $u$ is defined on $\RR^d$, takes values in
$\RR
$, is of class $\C^\infty$ and its domain is contained in a compact
subset of $(\G^1 \cup\G^2)^c$. Then for $n$ large enough, $u$ is null
on $\G^1_n \cup\G^2_n$, thus
\[
\int_\Gamma(\vec{\sigma} \cdot\vv_{\O}) u \,d
\H^{d-1} = 0,
\]
and we conclude that $\vec{\sigma} \cdot\vv_\O=0$, $\H
^{d-1}$-a.e. on
$\G\setminus(\G^1 \cup\G^2)$. On the other hand, if $u \in\C
_c^\infty( (\G\setminus\G^1)^c,\RR)$, then for $n$ large enough,
\[
\int_\Gamma(\vec{\sigma} \cdot\vv_{\O}) u\, d
\H^{d-1} = \lim_{n\rightarrow\infty} - \frac{1}{n^{d-1}} \sum
_{x \in\G_n^1 } u(x) \widehat f_n(x),
\]
and if for all $n\geq1$ we have $\widehat f_n (x)\geq0$ for all $x\in
\G^1_n$, then we conclude that $\vec\sigma\cdot\vv_{\O} \leq0$
$\H
^{d-1}$-a.e. on $\G^1$.
\end{pf*}

\begin{rem}
\label{remdiv2}
Notice that combining equations (\ref{eqdefdiv}) and (\ref
{eqdefbord2}), we obtain that
%
%
\begin{equation}
\label{eqdefdiv2} \forall h \in C_c^\infty \bigl(
\RR^d, \RR \bigr) \qquad\int_{\RR^d} h \di\vec {\sigma }
\,d \L^d = - \int_\Gamma h (\vec{\sigma} \cdot
\vv_{\O}) \,d \H ^{d-1}.
\end{equation}
This implies that $\di\vec\sigma$ is a distribution of order $0$ on
$\RR^d$,
\[
\forall K \mbox{ compact of } \RR^d, \exists C_K,
\forall h\in\C^\infty_c (K, \RR)\qquad \biggl\llvert \int
_{\RR^d} h \di\vec\sigma \,d\L^d \biggr\rrvert \leq
C_K \| h \|_\infty,
\]
since we can choose $C_K = C = \| \vec\sigma\cdot\vv_{\O} \|
_\infty\H
^{d-1} (\G)$ for any compact $K$. Thus by the Riesz representation
theorem (see Theorem 6.19 in \cite{Rudinrc}) we know that it is a Radon
measure that we denote by $\di\vec\sigma\L^d$, and this measure is
completely characterized by equation (\ref{eqdefdiv2}), that is,
%
%
\begin{equation}
\label{eqdefdiv3} \di\vec\sigma\L^d = - (\vec\sigma\cdot
\vv_{\O}) \H ^{d-1}\vert _{\G}.
\end{equation}
\end{rem}


\subsection{Capacity constraint}
\label{secconstr}

In this section, we prove the following proposition:
%
\begin{prop}
\label{propcap1}
If $\vec\mu= \vec\sigma\L^d$ is the limit of a subsequence of
$(\vec\mu
_n )_{n\geq1}$, then almost surely we have
%
%
\begin{equation}
\label{eqpropcap1} \L^d\mbox{-a.e. on } \RR^d, \forall
\vv\in\SS^{d-1} \qquad\vec {\sigma} \cdot\vv\leq\nu(\vv).
\end{equation}
\end{prop}
\begin{pf}
We explain first the idea of the proof. The convergence $\vec\mu_n
\rightharpoonup\vec\sigma\L^d$ implies that
\[
\int_D \vec\sigma\cdot\vv \,d\L^d = \lim
_{n\rightarrow\infty} \int_D \,d\vec\mu_n
\cdot\vv
\]
for every Borel set $D$ such that $\L^d(\partial D)=0$. On one hand, using
Lebesgue differentiation theorem, we know that for $\L^d$-a.e. $x$,
\[
\frac{1}{\L^d (D(x,\eps))}\int_{D(x,\eps)} \vec\sigma\cdot\vv \,d
\L^d
\]
converges toward $\vec\sigma(x) \cdot\vv$ when $\eps$ goes to zero,
where $D(x,\eps)$ is a ``nice'' sequence of neighborhoods of $x$ of
diameter $\eps$. To conclude that $\vec\sigma\cdot\vv$ is bounded by
$\nu(\vv)$, it remains to compare $\int_D \,d\vec\mu_n \cdot\vv$ with
$\nu(\vv)$. Proposition \ref{propflux} states that when $D$ is a
cylinder of height $h$ in the direction $\vv$, $\int_D \,d \vec\mu_n
\cdot\vv$ is close to $h \Psi(\vec\mu_n, D, \vv)/n^{d-1}$, where $
\Psi(\vec\mu_n, D, \vv)$ is the amount of fluid that crosses $D$ from
the lower half part to the upper half part of its boundary in the
direction $\vv$ according to the stream $\vec\mu_n$. Since $\tau_n (D,
\vv)$ is the maximal value of such a flow, $ \Psi(\vec\mu_n, D,
\vv)
\leq\tau_n (D, \vv)$, and we can conclude the proof by using the
convergence of the rescaled flow $\tau_n (D, \vv)$ toward $\nu(\vv)$.
The key argument---and the less intuitive---is Proposition \ref
{propflux}. In fact, if $ \vec l $ is a $\C^1$ vector field on $D$ with
null divergence and such that $\vec l \cdot\vv_{D} = 0$ $\H^{d-1}$-a.e.
on the vertical faces of $D$ (the ones who are not normal to $\vv$), if
we denote by $B$ the basis of $D$, then by Fubini theorem we have
\[
\int_D \vec l \cdot\vv \,d\L^d = \int
_0^h \biggl(\int_{B + u \vv}
\vec l \cdot \vv\, d \H^{d-1} \biggr) \,du
\]
and we have for all u
\[
\int_{B + u \vv} \vec l \cdot\vv\, d \H^{d-1} = \int
_B \vec l \cdot\vv\, d \H^{d-1}
\]
since by the Gauss--Green theorem we get
\[
\int_{\partial D} \vec l \cdot\vv_D \,d
\H^{d-1} = \int_D \di\vec l \,d\L ^d =
0.
\]
We obtain that
\[
\int_D \vec l \cdot\vv \,d\L^d = h \int
_{B } \vec l \cdot\vv\, d \H ^{d-1}
\]
and $\int_{B } \vec l \cdot\vv\, d \H^{d-1} $ is indeed the flow that
goes from the bottom to the top of $D$ according to $\vec l$. In the
proof of Proposition \ref{propflux}, we adapt this argument to a
discrete stream $\vec\mu_n$, and we consider a cylinder flat enough
(i.e., $h$ small enough) to control the amount of fluid that enters in
$D$ or escapes from $D$ through its vertical faces.

\emph{Step \textup{1:} From $\vec\sigma(x) \cdot\vv$ to $\int_{cyl(x_p+ p^{-1} A, p^{-1} h)} \vec\sigma\cdot\vv
\,d\L^d$}.
Since the functions $\vv\in\SS^{d-1} \mapsto\vec{\sigma}(x)
\cdot\vv$
(for a fixed realization and a fixed $x$) and $\vv\in\SS^{d-1}
\mapsto
\nu(\vv)$ are continuous, property (\ref{eqpropcap1}) is equivalent to
%
%
\begin{equation}
\label{eqcap2} \L^d\mbox{-a.e. on } \RR^d, \forall\vv
\in\widehat\SS^{d-1}\qquad\vec{\sigma} \cdot\vv\leq\nu(\vv),
\end{equation}
where $\widehat\SS^{d-1}$ denotes the set of all the unit vectors of
$\RR^d$ that define a rational direction.
Theorem \ref{thmnu} states that for every cylinder $\cyl(A,h)$ with
$A$ a non degenerate hyperrectangle normal to $\vv$ and $h>0$, we have
%
%
\begin{equation}
\label{eqcap1} \lim_{n\rightarrow\infty} \frac{\tau_n (\cyl(A,h), \vv)}{
n^{d-1} \H
^{d-1} (A)} = \nu(\vv)\qquad
\mbox{a.s.}
\end{equation}
Thus these convergences hold a.s. simultaneously for all cylinders
$\cyl
(A,h)$ whose vertices have rational coordinates. For the rest of this
section, we consider a fixed realization of the capacities on which
these convergences happen.

According to Definition 7.9 in \cite{Rudinrc}, we say that a sequence
$(V_p)_{p\geq1}$ of Borel sets in $\RR^d$ \emph{shrinks to a point
$x\in\RR^d$ nicely} if there exists a number $\alpha>0$ and a
sequence of positive real numbers $(r_p)_{p\geq1}$ satisfying $\lim_{p\rightarrow\infty} r_p = 0$, and for all $p\geq1$,
\[
V_p \subset B(x,r_p)\quad \mbox{and} \quad \L^d
(V_p ) \geq \alpha\L^d \bigl(B(x,r_p)
\bigr).
\]
We need the Lebesgue differentiation theorem on $\RR^d$ (see Theorem
7.10 in~\cite{Rudinrc}):
%
\begin{thmm}
\label{thmBesicyl}
Let $g$ be a Borel function in $ L^1 (\RR^d, \L^d)$. To each $x\in
\RR
^d$, associate a sequence $(V_p(x))_{p\geq1}$ of Borel sets in $\RR^d$
that shrinks to $x$ nicely. Then for $\L^d$-a.e. $x\in\RR^d$,
\[
\lim_{p \rightarrow\infty} \frac{1}{\L^d (V_p(x))} \int_{V_p(x)}
g \,d \L^d = g(x).
\]
\end{thmm}
To each point $x\in\RR^d$, we associate a deterministic sequence $(x_p
(x))_{p\geq1}$ of points of $\RR^d$ that have rational coordinates and
satisfying $\| x_p(x) - x \| \leq1/p$. Then for every nondegenerate
cylinder $D$ of center $0$, the sequence of Borel sets $(x_p(x) +
p^{-1} D)_{p\geq1}$ shrinks to $x$ nicely. We apply Theorem \ref
{thmBesicyl} to the function $\vec{\sigma}$ (coordinate by
coordinate) to
obtain that for $\L^d$-a.e. $x$, for every cylinder $D= \cyl(A,h)$ of
center $0$ and whose vertices have rational coordinates (we say that
$D$ is a rational cylinder), we have
\[
\lim_{p \rightarrow\infty} \frac{1}{p^{-d} \L^d(D)} \int_{x_p(x)+p^{-1} D}
\vec{\sigma} \,d\L^d = \vec{\sigma}(x).
\]
From now on, we consider a fixed $x$ [thus a fixed sequence
$(x_p(x))_{p\geq1} $ that we denote by $(x_p)_{p\geq1}$] such that the
previous convergence holds for every rational cylinder $D$.

Let $\vv$ be a vector in $\widehat\SS^{d-1}$ and $\eta$ a positive
real number. There exists a positive real number $\lambda$ such that
$\lambda\vv$ has integer coordinates. If $(\vec\bbf_1,\ldots,
\vec\bbf
_d)$ is the canonical basis on $\RR^d$, suppose for instance $\lambda'
\vv\cdot\vec\bbf_1 \neq0$, then $(\lambda' \vv, \vec\bbf_2,
\ldots, \vec
\bbf_d)$ is a basis of $\RR^d$. Adapting slightly the Gram--Schmidt
process, we can obtain an orthogonal basis $(\vv, \vec u_2,\ldots,
\vec u_d)$
of $\RR^d$ such that all the vectors $\vec u_i$, $i=2,\ldots,d$ have integer
coordinates. Thus there exists a non degenerate hyperrectangle $A$ of
center $0$, normal to $\vv$ and whose vertices have rational
coordinates. Then for every positive rational number $h$ [thus $\cyl
(A,h)$ is a rational cylinder], there exists $p_0 (x,A,h, \eta)
<\infty
$ such that for all $p\geq p_0$ we get
%
%
\begin{eqnarray}
\label{eqcap3} \qquad &&\biggl\llvert \frac{1}{ \L^d(\cyl(x_p+ p^{-1} A, p^{-1} h))} \int_{\cyl
(x_p+p^{-1} A, p^{-1} h)}
\vec{\sigma} \cdot\vv \,d\L^d - \vec {\sigma }(x) \cdot\vv \biggr
\rrvert
\nonumber
\\[-8pt]
\\[-8pt]
\nonumber
&&\qquad \leq\eta.
\end{eqnarray}
The value of $h$ will be fixed later; see step $3$ below.

\emph{Step \textup{2:} From $\int_{\cyl(x_p+ p^{-1} A, p^{-1} h)}
\vec
{\sigma} \cdot\vv \,d\L^d$ to $\int_{\cyl(x_p+ p^{-1} A, p^{-1} h)}
\,d \mu_n \cdot\vv$}.
Let $h$ be a fixed positive number and $p$ be a fixed integer, $p\geq
1$. Up to extraction of a subsequence $\vec\mu_n \rightharpoonup
\vec
\sigma\L^d$, and since $\L^d (\partial\cyl(x+\eps A, \eps h)) = 0$,
by Portmanteau's theorem we have
\[
\int_{\cyl(x_p + p^{-1} A, p^{-1} h)} \vec{\sigma} \cdot\vv \,d\L^d = \lim
_{n\rightarrow\infty} \int_{ \cyl(x_p+p^{-1} A, p^{-1} h)} \,d\vec{
\mu}_n \cdot\vv.
\]
Thus we obtain that for all $n$ large enough (how large depending on
$x,A,\break h, \eta, p$)
%
%
\begin{eqnarray}
\label{eqcap4}&& \biggl\llvert \int_{\cyl(x_p+p^{-1} A,p^{-1} h)} \vec{\sigma} \cdot
\vv \,d\L^d - \int_{ \cyl(x_p+p^{-1} A, p^{-1} h)} \,d\vec{
\mu}_n \cdot \vv \biggr\rrvert
\nonumber
\\[-8pt]
\\[-8pt]
\nonumber
&&\qquad\leq\eta\L^d \bigl(\cyl
\bigl(x_p+p^{-1} A, p^{-1} h \bigr) \bigr).
\end{eqnarray}

\emph{Step \textup{3:} From $\int_{\cyl(x_p+p^{-1} A, p^{-1} h)}
\,d \mu_n \cdot\vv$ to $\Psi(\vec\mu_n, \cyl(x_p+p^{-1} A, p^{-1}
h), \vv)$}.
For any $h>0$, any nondegenerate hyperrectangle $A$, we denote by
$\Psi
(\vec{\mu}_n,\break  \cyl(A,h), \vv)$ the flow that crosses $\cyl
(A,h)$ from
the lower half part of its boundary to the upper half part of its
boundary in the direction $\vv$ according to the stream $\vec{\mu
}_n$ on
$(\ZZ^d_n, \EE^d_n)$,
\begin{eqnarray*}
&&\Psi \bigl(\vec{\mu}_n, \cyl(A,h), \vv \bigr) \\
&&\qquad= \sum
_{e\in\cyl(A,h),
e=[a,b], a\in B'(A,h), b\notin B'(A,h)} f_n(e) ( \mathbh {1}_{\{ \be= \langle a,b\rangle\}} -
\mathbh{1}_{\{ \be= \langle b,a
\rangle\}} ),
\end{eqnarray*}
where, if $z$ is the center of $\cyl(A,h)$ and $\vv$ is normal to $A$,
the set $B'(A,h)$ is defined by
\begin{eqnarray*}
&&B'(A,h) = \Biggl\{ x \in\cyl(A,h) \cap\ZZ^d_n
\Big\vert\\
&&\hspace*{60pt}\begin{array} {l} \overrightarrow{zx} \cdot\vv<0 \mbox{ and}
\\
\exists y \notin\cyl(A,h), [x,y] \in\EE^d_n \mbox{ and
} [x,y] \cap\partial\cyl(A,h) \neq\varnothing \end{array} %
\Biggr\}.
\end{eqnarray*}
We state the following property:
%
\begin{prop}
\label{propflux}
Let $A$ be a nondegenerate hyperrectangle normal to a unit vector $\vv
$, and let $\eta>0$. There exists $h_0 (A,\eta)$ such that for $0<
h\leq h_0$, for $n$ large enough (how large depending on everything
else), we have
\[
\biggl\llvert \int_{ \cyl(A,h)} \,d\vec{\mu}_n \cdot
\vv- \frac{2 h
\Psi(\vec
{\mu}_n, \cyl( A, h), \vv) }{n^{d-1}} \biggr\rrvert \leq\eta\L ^d \bigl(\cyl(A,h)
\bigr),
\]
and for all $\eps>0$ and $y\in\RR^d$, we have
\[
h_0 (y + \eps A, \eta) = \eps h_0 (A,\eta).
\]
\end{prop}
Before proving Proposition \ref{propflux}, we end the proof of
Proposition \ref{propcap1}. We apply Proposition \ref{propflux} to the
hyperrectangle $A$ to obtain an $h_0 (A, \eta)$ and then we use the
rescaling property in Proposition \ref{propflux} to obtain that for all
$h\leq h_0 (A, \eta)$, for all $p\geq1 $, and for all $n$ large enough
(how large depending on everything else), we have
%
%
\begin{eqnarray}\qquad
\label{eqcap5}&& \biggl\llvert \int_{\cyl(x_p+p^{-1} A,p^{-1} h)} \,d\vec{
\mu}_n \cdot \vv- \frac{2 p^{-1} h \Psi(\vec{\mu}_n, \cyl(x_p +p^{-1} A, p^{-1} h),
\vv)
}{n^{d-1}} \biggr\rrvert
\nonumber
\\[-8pt]
\\[-8pt]
\nonumber
&&\qquad\leq\eta
\L^d \bigl(\cyl \bigl(x_p+p^{-1}
A,p^{-1} h \bigr) \bigr).
\end{eqnarray}

\emph{Step \textup{4:} Conclusion}.
Let us combine the previous steps. Theorem \ref{thmnu} states that for
every cylinder $\cyl(A,h)$ with $A$ a nondegenerate hyperrectangle
normal to $\vv$ and $h>0$, we have
\[
\lim_{n\rightarrow\infty} \frac{\tau_n (\cyl(A,h), \vv)}{
n^{d-1} \H
^{d-1} (A)} = \nu(\vv)\qquad \mbox{a.s.}
\]
Thus these convergences hold a.s. simultaneously for all the rational
cylinders, that is, cylinders with rational vertices [like the
cylinders $(\cyl(a_i,h),\break  i \in\I)$]. We consider a fixed realizations
of the capacities on which these convergences occur. We consider a
point $x\in\O$ as explained in step $1$, a vector $\vv\in\hat\SS
^{d-1}$ and a nondegenerate hyperrectangle $A$ normal to $\vv$, of
center $0$ and whose vertices have rational coordinates. We fix $\eta
>0$. We choose a positive rational $\widehat h_0 \leq h_0 (A, \eta)$ as
given in Proposition \ref{propflux} in step $3$, then $p_0 (x, A,
\widehat h_0, \eta) $ as defined in step $1$, and combining
inequalities (\ref{eqcap3}), (\ref{eqcap4}) and (\ref{eqcap5}) applied
with $p = p_0$, we get that for $n$ large enough (as large as required
in steps $2$ and $3$), we have
\[
\biggl\llvert \vec{\sigma} (x) \cdot\vv- \frac{ 2 p_0^{-1} \widehat h_0
\Psi(\vec{\mu}_n, \cyl(x_{p_0} +p_0^{-1} A, p_0^{-1} \widehat
h_0), \vv)
}{n^{d-1} \L^d (\cyl(x_{p_0}+p_0^{-1} A,p_0^{-1} \widehat h_0))} \biggr\rrvert \leq3
\eta.
\]
Since by maximality of $\tau$ we know that $\Psi(\vec{\mu}_n, \cyl
(x_{p_0} +p_0^{-1} A, p_0^{-1} \widehat h_0), \vv) \leq\tau_n (x_{p_0}
+p_0^{-1} A, p_0^{-1} \widehat h_0)$ we obtain, for all $n$ large enough,
\begin{eqnarray*}
\vec{\sigma} (x) \cdot\vv&\leq&\frac{ 2 p_0^{-1} \widehat h_0\tau_n
(x_{p_0} +p_0^{-1} A, p_0^{-1} \widehat h_0) }{n^{d-1} \L^d (\cyl(
x_{p_0}+ p_0^{-1} A,p_0^{-1} \widehat h_0))} + 3 \eta\\
&= &\frac{\tau
_n (x_{p_0}+ p_0^{-1} A, p_0^{-1} \widehat h_0)}{n^{d-1} \H^{d-1}
(p_0^{-1} A)} + 3
\eta.
\end{eqnarray*}
Since the cylinder $\cyl(x_{p_0}+ p_0^{-1} A, p_0^{-1} \widehat h_0)$
is rational, we get, when $n$ goes to infinity,
\[
\vec{\sigma}(x) \cdot\vv\leq\nu(\vv) + 3 \eta.
\]
Thus $\vec{\sigma} (x) \cdot\vv\leq\nu(\vv)$, and (\ref
{eqcap2}) and
Proposition \ref{propcap1} are proved.
\end{pf}

\begin{pf*}{Proof of Proposition \ref{propflux}} We give first the
idea of the proof. We recall what it would be if we would consider a
continous regular stream $ \vec l $, that is, a $\C^1$ vector field
on $D
= \cyl(A,h)$, with null divergence and such that $\vec l \cdot\vv
_{D} =
0$ $\H^{d-1}$-a.e. on the vertical faces of $D$ (the ones who are not
normal to $\vv$). If we denote by $B$ the basis of $D$, then by
Fubini's theorem we would get
\begin{eqnarray*}
\int_D \vec l \cdot\vv \,d\L^d & = &\int
_{0}^{2h} \biggl(\int_{B + u \vv}
\vec l \cdot\vv\, d \H^{d-1} \biggr) \,du
\\
& =& 2h \times\mbox{``flow from the top to the bottom of }D\mbox{ according to }
\vec l \mbox {.''}
\end{eqnarray*}
We have to adapt this argument to $\vec\mu_n$. The set $B + u\vv$
is a
continuous cutset that separates the top from the bottom of $D$. The
equivalent discrete cutset is, roughly speaking, the set of edges
\[
\E_u = \{ e\subset D | e \cap B + u \vv\neq\varnothing\}.
\]
The flow that crosses $B + u \vv$ according to $\vec\mu_n$ is $\sum_{e\in\E_u} f_n(e)$, and it is almost equal to $\Psi(\vec\mu_n,
D, \vv
)$ up to an error which is due to the flow that can cross the vertical
faces of $D$; thus we can control it if the height of $D$ is small
enough. Then we get almost
\[
\int_{0}^{2h} \sum_{e\in\E_u}
f_n(e) \,du = 2h \Psi(\vec\mu_n, D, \vv).
\]
As in the continuous case, the left-hand side of the previous equality
is almost the integral of the stream over $D$,
\begin{eqnarray*}
\int_0^{2h} \sum_{e\in\E_u}
f_n(e) \,du &=& \sum_{e\subset D}
f_n(e) \int_0^{2h}
\mathbh{1}_{e \in\E_u} \,du = \sum_{e\subset
D}
f_n(e) \frac{1}{n} \vec{e}\cdot\vv\\
&=& n^{d-1} \int
_D \,d\vec\mu_n \cdot\vv,
\end{eqnarray*}
up to a small error that appears for edges located near the boundary of $D$.

We begin now the proof. We will use another property, Proposition \ref
{propflux1}, that will be proved after the end of the proof of
Proposition \ref{propflux}. We give first the expression of $\Psi
(\vec
{\mu}_n, \cyl(A,h), \vv)$ in terms of $f_n$. Let $E$ be a
$(B'(A,h),\break
T'(A,h))$-cutset in $\cyl( A, h)$. We define $s(E) \subset\ZZ^d_n $ by
\[
s(E) = \lleft\{ y\in\ZZ^d_n \cap\cyl(A, h) \left\vert
\begin{array} {l} \mbox{there exists a path from }y \mbox{ to }
B'( A, h)
\\
\mbox{made of edges in } \bigl(\EE^d_n \cap\cyl( A, h)
\bigr) \setminus E \end{array} %
\right.
\rright\}.
\]
The set $s(E)$ is the connected component of $B'( A, h)$ in $(\ZZ^d_n,
\EE^d_n \setminus E) \cap\cyl( A, h) $. We consider a ``non
discrete version'' $S(E)$ of $s(E)$, defined by
\[
S(E) = \biggl(s(E) + \frac{1}{2n}[-1, 1] ^d \biggr)\cap\cyl( A,
h)
\]
%
\begin{figure}

\includegraphics{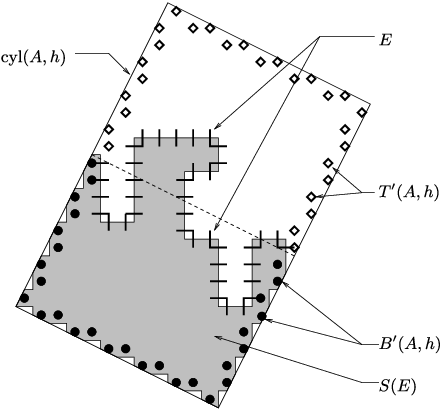}

\caption{A $(B'(A,h), T'(A,h))$-cutset $E$ in $\cyl(A,h)$ and the
corresponding set $S(E)$.}
\label{figensemblecyl}
\end{figure}
\hspace*{-2pt}[this is a subset of $\RR^d$ included in $\cyl( A, h)$; see Figure~\ref{figensemblecyl}].
%
%
For each edge $e\in E$, $c(e)$ belongs to $\partial S(E)$ and the
exterior unit vector normal to $S(E)$ at $c(e)$, which we denote by
$\vv
_{S(E)} (c(e))$, is equal to $\vec{e}$ or $-\vec{e}$, thus $\vec
{e}\cdot\vv_{S(E)}
(c(e))$ equals $+1$ or $-1$. If $\vec{e}\cdot\vv_{S(E)} (c(e)) = +1$
(resp., $-1$), then $\be= \langle a,b \rangle$ with $a \in s(E)$
[resp.,
$b\in s(E)$] and $b$ in the connected component of $T'( A, h)$ in $(\ZZ
^d_n, \EE^d_n \setminus E) \cap\cyl( A, h)$ (resp., $a$ in this
component). Indeed, $E$ is minimal thus if we remove $e$ from $E$ we
create a path from $B'(A,h)$ to $T'(A,h)$ that contains $e$. By the
node law, we know that $\Psi(\vec{\mu}_n, \cyl( A, h), \vv)$ is
equal to
the flow that crosses $E$ according to $\vec{\mu}_n$, that is,
%
%
\begin{equation}
\label{eqcalculPsi} \Psi \bigl(\vec{\mu}_n, \cyl( A, h), \vv \bigr)
= \sum_{e\in E} f_n(e) \vec{e} \cdot
\vv_{S(E)} \bigl(c(e) \bigr).
\end{equation}
We construct now several such cutsets inside $\cyl( A, h)$. By
symmetry, we can suppose that all the coordinates of $\vv$ are
nonnegative. Let $x$ be the center of $A$. Let $u \in\RR$. We define
the hypersurface $\P(u)$ by
\[
\P(u) = \bigl\{ y \in\RR^d | \overrightarrow{xy} \cdot\vv= u- h
\bigr\}.
\]
For each edge $e$ such that $\be= \langle a,b \rangle$, we define
$\tilde e = [a,b[$, the segment that includes~$b$, the endpoint of $\be
$, but excludes $a$, its origin. We define the set of edges $E_n(u)$ by
\[
E_n(u) = \bigl\{ e\in\EE^d_n | e \subset
\cyl( A, h) \mbox{ and } \tilde e \cap\P(u) \neq\varnothing \bigr\};
\]
see Figure~\ref{figcutsets}. We define also the set of edges $F_n$ by
\[
F_n = \bigl\{ e \in\EE^d_n | e \subset
\cyl( A, h) \cap\V_2 \bigl(\cyl( \partial A, h), 2d/n \bigr) \bigr\},
\]
which is the set of the edges in $\cyl( A, h)$ that are near the faces
of the cylinder that are normal to $A$, and the set $\widetilde E_n(u)$ by
\[
\widetilde E_n(u) = \bigl\{ e\in E_n(u) | e \not
\subset\V_2 \bigl(\cyl( \partial A, h), 4d/n \bigr) \bigr\},
\]
which is the set of the edges of $E_n(u)$ that are not too close from
the faces of the cylinder that are normal to $A$.
%
\begin{figure}

\includegraphics{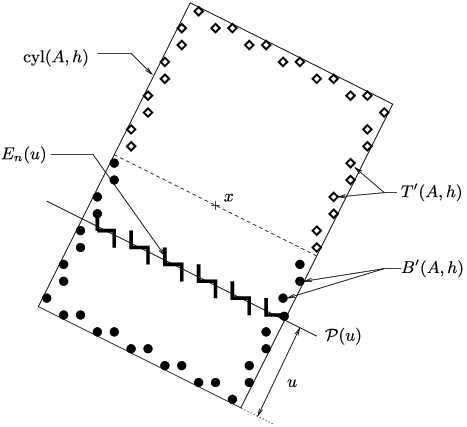}

\caption{The set $E_n(u)$.}\label{figcutsets}
\end{figure}
%
%
We need the following property:
%
\begin{prop}
\label{propflux1}
For all $u \in[1/n, 2h - 1/n]$, $E_n(u) \cup F_n$ contains a
$(B'(A,h), T'(A,h))$-cutset in $\cyl( A, h)$. We denote such a cutset
by $\widehat E_n(u)$. Necessarily $\widetilde E_n(u)$ is included in
$\widehat E_n(u)$ (whichever way we construct it), and
\[
\forall e\in\widetilde E_n(u)\qquad \vec{e}\cdot\vv_{S(\widehat E_n(u))}
\bigl(c(e) \bigr) = +1.
\]
\end{prop}
We consider such a cutset $\widehat E_n(u)$ for a given $u$ in $[1/n, 2
h - 1/n]$. Using equation (\ref{eqcalculPsi}) we obtain as in Section~\ref{seccomp} that
%
%
\begin{eqnarray}
\label{eqPsisum} &&\biggl\vert\Psi \bigl(\vec{\mu}_n, \cyl( A, h), \vv
\bigr) - \sum_{e\in E_n(u)} f_n(e)\biggr \vert
\nonumber
\\
& &\qquad= \biggl\vert\sum_{e \in\widehat E_n(u)} f_n (e) \vec{e}
\cdot\vv _{S(\widehat E_n(u))} \bigl(c(e) \bigr) - \sum_{e\in E_n(u)}
f_n(e) \biggr\vert
\nonumber
\\
&&\qquad = \biggl\vert\sum_{e \in\widehat E_n(u) \setminus
\widetilde E_n(u)} f_n(e) \vec{e}
\cdot\vv_{S(\widehat E_n (u)) } \bigl(c(e) \bigr) - \sum_{e\in E_n(u) \setminus\widetilde E_n(u)}
f_n(e)\biggr \vert
\nonumber
\\[-8pt]
\\[-8pt]
\nonumber
&&\qquad \leq M \bigl[ \card \bigl(\widehat E_n(u) \setminus
\widetilde E_n(u) \bigr) + \card \bigl(E_n(u)
\setminus \widetilde E_n(u) \bigr) \bigr]
\\
&&\qquad \leq4d M \card \bigl(\ZZ^d_n \cap\V_2
\bigl(\cyl( \partial A, h), 4d n^{-1} \bigr) \bigr)
\nonumber
\\
&&\qquad \leq4d M n^d \L^d \bigl( \V_2 \bigl(\cyl(
\partial A, h), 5d n^{-1} \bigr) \bigr)
\leq C M\H^{d-2} ( \partial A) \bigl(h+ 5
dn^{-1} \bigr) n^{d-1}\nonumber
\end{eqnarray}
for a constant $C$. Let us consider the quantity
\[
\gamma= \int_{0}^{2 h} \biggl( \sum
_{e\in E_n(u)} f_n(e) \biggr) \,du.
\]
On one hand, inequality (\ref{eqPsisum}) states that
\begin{eqnarray*}
&& \bigl\vert2 h \Psi \bigl(\vec{\mu}_n, \cyl( A, h), \vv \bigr) - \gamma
\bigr\vert
\\
&&\qquad \leq2 h CM\H^{d-2} ( \partial A) \bigl(h+ 5d n^{-1} \bigr)
n^{d-1}
\\
&&\qquad\quad{} + \int_{[0, 1/n]\cup[ 2 h -1/n, 2 h]} \biggl\llvert \sum
_{e\in E_n(u)} f_n(e) -\Psi \bigl(\vec{
\mu}_n, \cyl( A, h), \vv \bigr) \biggr\rrvert \,du.
\end{eqnarray*}
Moreover, there exists a constant $C'(d,A)$ such that
\begin{eqnarray*}
\forall u \in\RR\qquad\card \bigl(E_n(u) \bigr) &\le& C'
(d,A) n^{d-1} \quad\mbox{and}\\
 \bigl\vert\psi \bigl(\vec{\mu}_n,
\cyl(A, h), \vv \bigr)\bigr\vert&\leq& C'(d,A) M n^{d-1}.
\end{eqnarray*}
We obtain the second inequality by noticing that the set of edges $E_n
( h)$ separates $B'(A,h)$ from $T'(A,h)$ in $\cyl( A, h)$, and the
first inequality bounds its cardinal. We conclude that
%
%
\begin{eqnarray}
\label{eqgamma1}&& \bigl\vert2 h \Psi \bigl(\vec{\mu}_n, \cyl( A, h), \vv
\bigr) - \gamma\bigr \vert\nonumber\\
&&\qquad \leq 2 C M h \bigl(h+5dn^{-1} \bigr)
\H^{d-2} (\partial A) n^{d-1}+ 4 C' (d,A) M
n^{d-2}
\\
&&\qquad \leq 2C M h^2 \H^{d-2} (\partial A) n^{d-1} +
K_1 (d,A,h,M) n^{d-2}.\nonumber
\end{eqnarray}
On the other hand, we have
\[
\gamma= \sum_{e \in\cyl( A, h)} f_n(e) \int
_0^{2 h} \mathbh {1}_{\{ e\in E_n(u) \}} \,du.
\]
For all $e\in\cyl( A, h)$, $e\notin E_n(u)$ if $u \notin[0, 2 h]$,
and we have
\[
\int_{0}^{2 h} \mathbh{1}_{\{ e\in E_n(u) \}} \,du =
\int_{\RR} \mathbh{1}_{\{ e\in E_n(u) \}} \,du =
\frac{1}{n} \vec{e}\cdot\vv.
\]
Indeed for all $e \in\cyl(A,h)$, if $\be= \langle a,b \rangle$, then
$\overrightarrow{xa} \cdot\vv\leq\overrightarrow{xb} \cdot\vv$
(remember that we supposed that the coordinates of $\vv$ are non
negative) and
\[
e \in E_n(u)\quad \iff\quad\overrightarrow{xa} \cdot\vv< u-h \leq
\overrightarrow{xb} \cdot\vv\quad\iff\quad u \in\,] h+ \overrightarrow{xa} \cdot\vv, h+
\overrightarrow{xb} \cdot\vv],
\]
and the measure of the above interval is
\[
\overrightarrow{xb} \cdot\vv- \overrightarrow{xa} \cdot\vv= \frac{1}{n}
\vec{e}\cdot\vv.
\]
Thus
%
%
\begin{eqnarray}
\label{eqgamma2}&& \biggl\vert\gamma- n^{d-1} \int_{\cyl( A, h)}
\,d \vec{\mu}_n \cdot \vv \biggr\vert
\nonumber
\\
& &\qquad= \biggl\llvert \frac{1}{n} \sum_{e \in\cyl( A, h)}
f_n(e) \vec{e} \cdot\vv- \frac{1}{n} \mathop{\sum
_{e\in\EE^d_n,}}_{ c(e) \in\cyl( A, h)
} f_n(e) \vec{e}\cdot\vv \biggr\rrvert
\nonumber
\\[-8pt]
\\[-8pt]
\nonumber
&&\qquad \leq\frac{M}{n} \card \bigl( \bigl\{ e\in\EE^d_n
| e \cap \partial\cyl(A, h) \neq\varnothing \bigr\} \bigr)
\nonumber
\\
&&\qquad \leq K_2 (d,A,h,M) n^{d-2}.\nonumber
\end{eqnarray}
Combining inequalities (\ref{eqgamma1}) and (\ref{eqgamma2}) we obtain
%
%
\begin{eqnarray}
\label{eqgamma3} \quad&& \biggl\llvert \int_{\cyl(A, h)} \,d\vec{
\mu}_n \cdot\vv- \frac{2
h \Psi
(\vec{\mu}_n, \cyl( A, h), \vv)}{n^{d-1}} \biggr\rrvert
\nonumber
\\[-8pt]
\\[-8pt]
\nonumber
&&\qquad \leq2CM h^2 \H^{d-2} (\partial A) + \bigl(K_1(d,A,h,M)
+ K_2(d,A,h,M) \bigr) n^{-1}.
\end{eqnarray}
We define
\[
h_0 (A,\eta) = \frac{\eta\H^{d-1} (A)}{4 C M \H^{d-2} (\partial
A)}.
\]
We deduce from inequality (\ref{eqgamma3}) that all $h\leq h_0$, for
all $n$ we have
\begin{eqnarray*}
&&\biggl\llvert \int_{\cyl(A, h)} \,d\vec{\mu}_n \cdot
\vv- \frac{2 h
\Psi
(\vec{\mu}_n, \cyl( A, h), \vv)}{n^{d-1}} \biggr\rrvert\\
&&\qquad \leq \frac
{\eta\L^d(\cyl(A,h))}{2} +
\bigl(K_1(A,h,M) + K_2(A,h,M) \bigr) n^{-1},
\end{eqnarray*}
and thus for $n$ large enough (how large depending on $A,h,M$) we
obtain the desired inequality. Moreover for all $\eps>0$, $y\in\RR^d$,
we immediatly obtain that
\[
h_0 (y + \eps A, \eta) = \frac{\eta\H^{d-1} (y+\eps A)}{ 4CM \H
^{d-2} (\partial(y+\eps A))} = \frac{\eta\eps^{d-1} \H^{d-1} (
A)}{ 4CM \eps^{d-2} \H^{d-2} (\partial A)} =
\eps h_0 (A,\eta).
\]
This ends the proof of Proposition \ref{propflux}.
\end{pf*}

\begin{pf*}{Proof of Proposition \ref{propflux1}} First of
all, we
prove that for all $u\in[1/n, 2 h - 1/n]$, $E_n(u)$ separates the
bottom from the top of $\cyl( A, h)$. Let us consider a self-avoiding
path $r$ from $B(A,h)$ to $T(A,h)$ in $\cyl( A, h)$. The path $r$
admits a continuous parametrization $r = (r_t)_{t\in[0,1]}$. Let $x$
be the center of $A$. The two sets
\begin{eqnarray*}
V_1 (u)& = &\bigl\{ y\in\RR^d | \overrightarrow{xy} \cdot
\vv< u - h \bigr\} \cap\cyl( A, h),
\\
V_2 (u)& =& \bigl\{ y\in\RR^d | \overrightarrow{xy} \cdot
\vv\geq u - h \bigr\} \cap\cyl( A, h)
\end{eqnarray*}
form a partition of $\cyl( A, h)$. The path $r$ starts in $V_1(u)$ and
ends in $V_2(u)$. Indeed, $B( A, h) \subset\V_2 ( A -h\vv, n^{-1} )
\subset V_1(u)$ because $u\geq n^{-1}$ and $T( A, h) \subset\V_2 (
A+h\vv, n^{-1} ) \subset V_2(u)$ because $u\leq2 h - n^{-1}$. Since
$r$ is continuous, there exists $t_0 \in[0,1]$ such that
\[
t_0 = \sup \bigl\{ t \in[0,1] | r_t \in
V_1(u) \bigr\}.
\]
We define the point $z=r_{t_0}$. It is obvious that $z\in\P(u)$; see
Figure~\ref{figcoupure1}. If $z \notin\ZZ^d_n $, then $z$ belongs only
to one edge $e \subset r \subset\cyl( A, h)$, and $z$ is not an
extreme point of $e$ so $z \in\P(u)$ implies that $e \in E_n(u)$. If
$z\in\ZZ^d_n$, then $z$ belongs to exactly two edges $e_1$ and $e_2$
that are included in $r$. By the definition of $t_0$, we know that one
of these edges, say $e_1$ for example, is included in the adherence of
$V_1(u)$, and the other one, $e_2$, is included in $V_2(u)$. Since all
the coordinates of $\vv$ are nonnegative, we conclude that $e_2 \in
E_n(u)$. This proves that $E_n(u)$ separates $T( A, h)$ from $B( A, h)$
in $\cyl(A, h)$.
%
\begin{figure}

\includegraphics{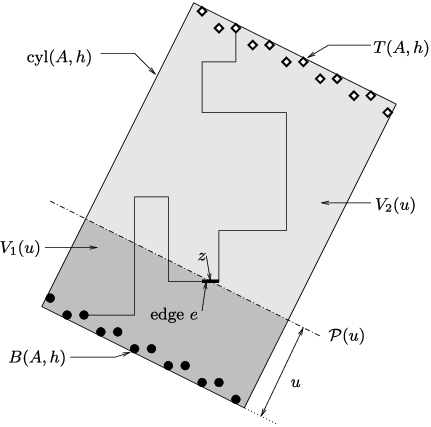}

\caption{The set $E_n(u)$ separates $B(A,h)$ from $T(A,h)$ in $\cyl(A,h)$.}
\label{figcoupure1}
\end{figure}

%

We deduce easily that $E_n(u)\cup F_n$ separates $T'( A, h)$ from $B'(
A, h)$ in $\cyl(A, h)$. Indeed, consider a path $\hat r$ from $T'(A,
h)$ to $B'( A, h)$ in $\cyl( A, h)$. If the starting point (resp., the
endpoint) of $\hat r$ belongs to $T'( A, h) \setminus T( A, h)$
[resp., $B'( A,h) \setminus B( A, h)$], then the first (resp., last)
edge of $r$ belongs to $F_n$. Otherwise, $\hat r$ is a path from $T(A,
h)$ to $B( A, h)$ in $\cyl( A, h)$, and we have proved that it must
contain at least one edge of $E_n(u)$.

\begin{figure}

\includegraphics{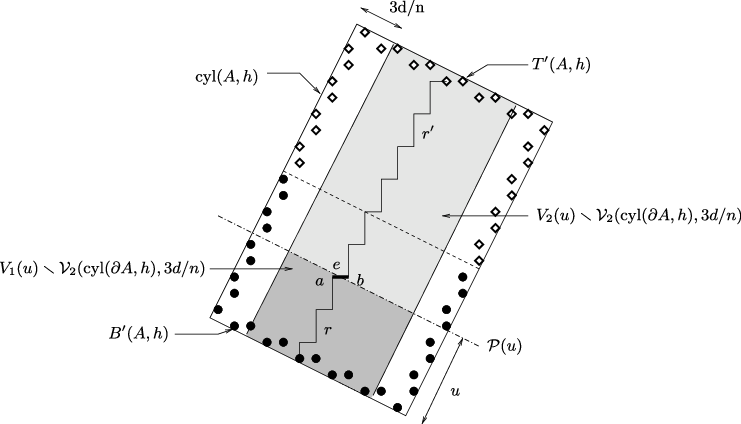}

\caption{Construction of a path from $B(A,h)$ to $T(A,h)$.}
\label{figcoupure2}
\end{figure}

We consider an edge $e$ of $\widetilde E_n(u)$, $\be= \langle a,b
\rangle$. Then $a\in V_1(u)$ and $b \in V_2(u)$. Moreover $e \not
\subset\V_2 (\cyl( \partial A, h), 4d/n)$ implies that $d_2 (a, \cyl(
\partial A, h))> 3d/n$. The set
\[
\D= V_1 (u) \setminus\V_2 \bigl(\cyl(\partial A, h)
, 3d/n \bigr)
\]
is a parallelepiped; thus the graph $\D\cap(\ZZ^d_n, \EE^d_n)$ is
connected. Let $r$ be a path from $a$ to $B(A,h)$ included in $\D$; see
Figure~\ref{figcoupure2}.
%
%
In the same way, there exists a path $r'$ from $b$ to $T( A, h)$ that
is included in $V_2(u)\setminus\V_2 (\cyl(\partial A, h), 3d/n)$.
Thus $r \cup e \cup r'$ is a path from $B( A, h)$ to $T( A, h)$ that
does not contain edges of $F_n \cup( E_n(u) \setminus\{ e \} )$,
and we conclude that $F_n \cup( E_n(u) \setminus\{ e \} )$ does
not separate $B( A, h)$ from $T( A, h)$ in $\cyl( A, h)$. This implies
that $e$ must belong to any cutset $\widehat E_n(u)$ with the
properties given in Proposition \ref{propflux1}. Moreover, we have
proved that $a\in S(\widehat E_n(u))$ and $b \notin S(\widehat
E_n(u))$, and this implies that $\vv_{S(\widehat E_n(u))} (c(e)) =
\vec{e}
$, so Proposition \ref{propflux1} is proved.
\end{pf*}


\subsection{Maximality}
\label{secmax}

We recall that
\[
\operatorname{flow}^{\mathrm{cont}}(\vec{\sigma}) = \int_{\G^1}
- \vec{\sigma} \cdot\vv_{\O
} \,d \H^{d-1}
\]
and
\[
\operatorname{flow}^{\mathrm{disc}}_n(\vec\mu_n) =
\frac
{1}{n^{d-1}} \sum_{e\in\Pi_n: e=[ab]
, a\in\G^1_n, b\notin\G^1_n} f_n (e) (
\mathbh{1}_{\{ \be=
\langle a,b\rangle\}} -\mathbh{1}_{\{ \be= \langle b,a\rangle\}} ).
\]
To complete the proof of Theorem \ref{thmpcpal}, we must prove that the
limit $\vec\mu= \vec\sigma\L^d$ of a subsequence of the
sequence of
maximal discrete flows $(\vec\mu_n^{\max})_{n\geq1}$ satisfies
%
%
\begin{equation}
\label{eqfluxmax} \operatorname{flow}^{\mathrm{cont}}(\vec\sigma) =
\phi_{\O}^b \qquad\mbox{a.s.}
\end{equation}
In this section, we prove the following result:
%
\begin{prop}
\label{propmax2}
Let $(\vec\mu_{n})_{n\geq1}$ be a sequence of admissible discrete
streams. If a subsequence $(\vec\mu_{\varphi(n)})_{n\geq1}$ converges
weakly toward a measure $\vec\mu= \vec\sigma\L^d$ with $\vec
\sigma
\in L^{\infty} (\RR^d \rightarrow\RR^d, \L^d) $, then
\[
\lim_{n \rightarrow\infty} \operatorname{flow}^{\mathrm
{disc}}_{\varphi(n)}(
\vec\mu_{\varphi(n)}) = \operatorname {flow}^{\mathrm{cont}} (\vec \sigma).
\]
\end{prop}
%
\begin{rem}
We will deduce equation (\ref{eqfluxmax}) from Proposition \ref
{propmax2} in Section~\ref{secccl}, using our study of minimal cutsets
and Lemma \ref{lemineg}.
\end{rem}

\begin{pf*}{Proof of Proposition \ref{propmax2}} The idea of the
proof is very similar to the one of Proposition \ref{propcap1}. Suppose
$\vec\sigma$ is very regular---$\C^1$ for example. By the Gauss--Green
theorem, we know that for all sets $E$ with finite perimeter,
\[
\int_{\partial E} \vec\sigma\cdot\vv_E \,d
\H^{d-1} = \int_E \di \vec \sigma \,d
\L^d = 0.
\]
If $\partial E = \G^1 \cup\S\cup\widehat\S$, where $\widehat\S
\subset\G
\setminus(\G^1 \cup\G^2)$ and $\S= \partial E \cap\O$,
then we obtain
\[
\int_{\S} \vec\sigma\cdot\vv_E \,d
\H^{d-1} = - \int_{\G^1} \vec\sigma \cdot
\vv_{\O} \,d \H^{d-1} + 0 = \operatorname{flow}^{\mathrm
{cont}}(
\vec\sigma).
\]
We can choose $E$ such that $\S$ is polyhedral: it allows us to cover
(up to a small volume) a neighborhood of $\S$ by a union of cylinders
$D_i$ of height $h$ and oriented in the direction $\vv_i$, where $\vv_i
= \vv_E$ on the face of $\S$ that $D_i$ crosses. As explained in the
sketch of the proof of Proposition \ref{propcap1}, $\int_{D_i} \vec
\sigma\cdot\vv_i \,d\H^{d-1}$ is very close to $h \int_{\S\cap D_i}
\vec\sigma\cdot\vv_i \,d\H^{d-1}$. Since $\vec\mu_n
\rightharpoonup\vec
\sigma\L^d$, $\int_{D_i} \vec\sigma\cdot\vv_i \,d\L^d$ is the limit
of $\int_{D_i} \,d\vec\mu_n \cdot\vv_i$. By Proposition \ref
{propflux}, we know that $\int_{D_i} \,d\vec\mu_n \cdot\vv_i$ is very
close to $\Psi(\vec\mu_n, D_i, \vv_i)$. Finally, we notice that
the flow
that crosses $\O_n$ from $\G^1_n$ to $\G^2_n$ is the flow that crosses
the $D_i$ up to a small error, and thus $\operatorname{flow}^{\mathrm
{disc}}_n(\vec\mu_n)$ is close to
$\sum_i \Psi(\vec\mu_n, D_i, \vv_i)$ properly rescaled. This is exactly
the idea we follow to compare $\operatorname{flow}^{\mathrm
{disc}}_n(\vec\mu_n)$ to $\operatorname{flow}^{\mathrm{cont}}(\vec
\sigma)$.
However, $\vec\sigma$ is not regular enough to allow a direct
application of the Gauss--Green theorem. The easiest way to get round
this problem is to come back to the definition of the divergence $\di
\vec\sigma$ to compare $\operatorname{flow}^{\mathrm{cont}}(\vec
\sigma)$ to a sum of the type
$\sum_i
\int_{D_i} \vec\sigma\cdot\vv_i \,d\L^d$.

From now on we consider a fixed realization of the capacities. We
consider a subsequence of $(\vec\mu_n)_{n\geq1}$ converging toward
$\vec
\mu$, but we still denote this subsequence by $(\vec\mu_n)_{n\geq
1}$ to
simplify the notation.

\emph{Step \textup{1:} From $\operatorname{flow}^{\mathrm
{cont}}(\vec\sigma)$ to $\sum_{i=1}^{\N} \int_{\cyl(A_i,h)} \vec\sigma\cdot\vv_i \,d\L^d$}.
For $\A$ a subset of $\RR^d$, we denote by $\overline{\A}$ its closure
and by $\stackrel{\circ}{\A}$ its interior. Let $P$ be a closed
polyhedral set of $\RR^d$ such that
\[
\overline{\G}^1 \subset\stackrel{\circ} {P}, \qquad\overline{
\G}^2 \subset\,\stackrel{\circ} { \wideparen{\RR^d
\setminus P}} \quad\mbox{and}\quad \partial P \mbox{ is transversal to } \G.
\]
The construction of such a set $P$ is made in Section~5 of \cite
{CerfTheret09supc}. The idea of the construction is the following. For
each $x\in\overline{\G}^1$, let $C_x$ be a closed cube of center $x$
and of positive size but small enough so that $d_2 (C_x, \G^2) \geq
d_2(\G^1, \G^2) /2$. The cubes $(C_x)_{x\in\overline{\G}^1}$ can be
chosen carefully so that their boundaries are transversal to~$\G$. Of course,
\[
\overline{\G}^1 \subset\bigcup_{x \in\overline{\G}^1}
\stackrel {\circ} {C}_x,
\]
and by compactness of $\overline{\G}^1$ we know that there exists a
finite subcovering of $\overline{ \G}^1$, say
\[
\overline{\G}^1 \subset\bigcup_{i=1}^p
\stackrel{\circ} {C}_{x_i}.
\]
We can take $P= \bigcup_{i=1}^p C_{x_i}$. By construction $\partial P$
is a
polyhedral hypersurface that is transversal to $\G$ and that does not
intersect $\overline{\G}^1 $ nor $\overline{\G}^2 $, thus $d_2
(\partial P, \G
^1 \cup\G^2) >0$. In the same way, for any $\zeta>0$, we can construct
a set $\O'$ satisfying $ \O\subset\O' \subset\V_2 (\O, \zeta)$ and
such that $ \partial\O' $ is polyhedral and transversal to
$\partial P $. We fix a
positive real number $\widehat\eta>0$. Since $\partial P$ is
transversal to
$\G$, there exists $\eps(\widehat\eta) >0 $ such that $\H^{d-1}
(\partial P
\cap(\V_2 (\O, \eps) \setminus\O)) \leq\widehat\eta$. We
consider a set $\O'$ corresponding to $\eps(\widehat\eta)$ as
described previously. Thus $\O'$ depends on $\O$, $P$ and $\widehat
\eta
$, and we have
%
%
\begin{equation}
\label{eqhateta} \H^{d-1} \bigl(\partial P \cap \bigl(\O'
\setminus\O \bigr) \bigr) \leq\widehat \eta.
\end{equation}
We need the following property (see Figure~\ref{figcylAi}):

%
\begin{figure}

\includegraphics{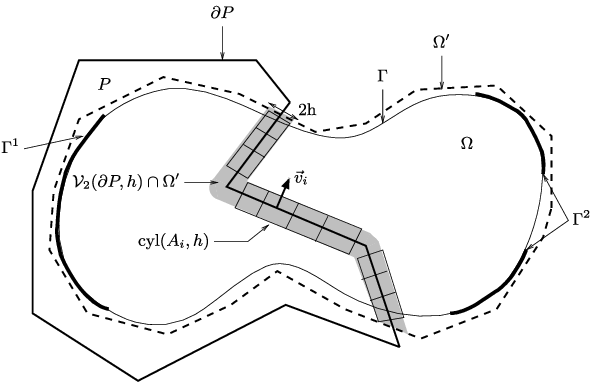}

\caption{The cylinders $(\cyl(A_i,h), i=1,\ldots,\N)$.}
\label{figcylAi}
\end{figure}

%
\begin{prop}
\label{propAi}
Let $\eta>0$. There exists a finite family of hyperrectangles
$A_1,\ldots,A_{\N}$ (depending on $\O, P, \widehat\eta, \eta$) of
disjoint interiors included in $\partial P \cap\O'$, a positive real number
$h_1 (\O, P, \eta, \widehat\eta)$ and a constant $C (\O, P,
\widehat
\eta)$ such that for all $h\leq h_1$, we have:
\begin{itemize}
\item$\H^{d-1}  ((\partial P\cap\O') \setminus
(\bigcup_{i=1}^{\N} A_i  )  ) \leq\eta$;
\item$\L^d  ( (\V_2(\partial P, h) \cap\O') \setminus
\bigcup_{i=1}^{\N} \cyl(A_i,h) ) \leq2 C \eta h$;
\item$\forall i \in\{ 1,\ldots,\N\}, \cyl(A_i,h) \subset\O'$;
\item$\forall i \in\{ 1,\ldots,\N\}, \forall x \in\cyl(A_i,h),
d_2(x, \partial P \setminus A_i) > d_2 (x, A_i) $.
\end{itemize}
\end{prop}

%
%
We admit this proposition for the time being. We fix a positive $\eta$.
For each $i\in\{ 1,\ldots,\N\}$, let $\vv_i$ be the exterior normal unit
vector to $P$ along the face of $\partial P$ on which $A_i$ is, thus
$\vv_i$
is normal to $A_i$.
We have explained in Remark \ref{remdiv2} that $\di\vec\sigma\L
^d =
- (\vec\sigma\cdot\vv_{\O}) \H^{d-1}\vert_{\G}$. Thus for any
function $\varphi\in W^{1,1}(\RR^d)$, we have
%
%
\begin{eqnarray}
\label{eqmaxbis1} \int_{\O} \vec\sigma\cdot\vec\nabla\varphi
\,d\L^d &=& \int_{\RR^d} \vec\sigma\cdot\vec\nabla
\varphi \,d \L^d = - \int_{\RR^d} \varphi \di\vec
\sigma \,d\L^d
\nonumber
\\[-8pt]
\\[-8pt]
\nonumber
& =& + \int_{\G} (\vec\sigma\cdot
\vv_{\O}) \gamma(\varphi\vert_{\O}) \,d \H^{d-1}
\end{eqnarray}
(this corresponds exactly to the definition of $\vec\sigma\cdot\vv
_{\O
}$ given in \cite{Nozawa}; cf. equation~(\ref{eqdefbord})). For a
positive $h\leq h_1$, we define the function $\varphi_h$ by
\[
\varphi_h (x) = \zeta \biggl( \frac{d_2 (x, P^c)}{h} \biggr) + \zeta
\biggl( \frac{h- d_2 (x,P)}{h} \biggr),
\]
where $\zeta(r) = r \mathbh{1}_{[0,1[} + \mathbh{1}_{[1,+\infty[}$.
Then $\varphi_h = 2$ on $P \cap\V_2(\partial P, h)^c$ and $\varphi
_h= 0$ on
$P^c \cap\V_2(\partial P, h)^c$, $\varphi_h$ is Lipschitz and has compact
support included in $P \cup\V_2(\partial P, h)$, in particular
$\varphi_h
\in W^{1,1} (\RR^d)$.
On one hand, we know that $\vec\sigma\cdot\vv_{\O} = 0 $ $\H
^{d-1}$-a.e. on $\G\setminus(\G^1\cup\G^2)$, and there exists
$h_2(\O, P) = d_2 (\partial P, \G^1\cup\G^2)/2>0$ such that for $h
\leq h_2
$ we have $\gamma(\varphi_h\vert_{\O} ) = \varphi_h = 2$ $\H
^{d-1}$-a.e. on $\G^1$ and $\gamma(\varphi_h\vert_{\O} ) = \varphi
_h =
0$ $\H^{d-1}$-a.e. on $\G^2$, thus
\[
\int_{\G} (\vec\sigma\cdot\vv_{\O}) \gamma(
\varphi_h\vert _{\O}) \,d \H ^{d-1} = 2 \int
_{\G^1} (\vec\sigma\cdot\vv_{\O}) \,d\H^{d-1}
= -2 \operatorname{flow}^{\mathrm{cont}}(\vec\sigma).
\]
On the other hand, we know that
\[
\vec\nabla\varphi_h (\cdot) = \mathbh{1}_{[0,1[} \biggl(
\frac
{d_2(\cdot, P^c)}{h} \biggr) h^{-1} \vec\nabla d_2 \bigl(
\cdot, P^c \bigr) - \mathbh{1}_{[0,1[} \biggl(
\frac{h-d_2(\cdot, P)}{h} \biggr) h^{-1} \vec\nabla d_2(\cdot, P),
\]
thus $\vec\nabla\varphi_h =0$ $\L^d$-a.e. on $\V_2(\partial P,
h)^c$, $ \| \vec
\nabla\varphi_h \|_\infty\leq h^{-1} $, and for all $i\in\{
1,\ldots,\N\}$ we have on $\cyl(A_i,h)$
\[
\vec\nabla\varphi_h = - h^{-1} \vv_i.
\]
For all $h\leq\min(h_1, h_2)$, equation (\ref{eqmaxbis1}) applied to
$\varphi_h$ gives
\begin{eqnarray*}
&&\operatorname{flow}^{\mathrm{cont}}(\vec\sigma)\\
&&\qquad = - \frac{1}{2} \int
_{\G} (\vec\sigma\cdot \vv_{\O
}) \gamma(
\varphi_h \vert_{\O}) \,d \H^{d-1}
\\
&&\qquad = -\frac{1}{2} \int_{\O} \vec\sigma\cdot\vec\nabla
\varphi _h \,d \L^d
\\
&&\qquad = -\frac{1}{2} \int_{\V_2 (\partial P, h)\cap\O} \vec\sigma \cdot\vec
\nabla\varphi_h \,d \L^d
\\
&& \qquad =-\frac{1}{2} \sum_{i=1}^{\N}
\int_{\cyl(A_i, h)} \vec\sigma \cdot \biggl( - \frac{1}{h}
\vv_i \biggr) \,d\L^d\\
&&\qquad\quad{} - \frac{1}{2} \int
_{(\V
_2 (\partial P, h) \cap\O)\setminus\bigcup_{i=1}^{\N} \cyl
(A_i,h)} \vec \sigma\cdot\vec\nabla\varphi_h \,d
\L^d.
\end{eqnarray*}
Thus
%
%
\begin{eqnarray}
\label{eqmax12} &&\Biggl\llvert \operatorname{flow}^{\mathrm{cont}}(\vec\sigma) -
\frac
{1}{2 h}\sum_{i=1}^{\N} \int
_{\cyl
(A_i,h)} \vec\sigma\cdot\vv_i \,d\L^d
\Biggr\rrvert\nonumber \\
&&\qquad \leq\frac
{1}{2} \| \vec\sigma\|_\infty\| \vec\nabla
\varphi_h\|_\infty\L ^d \Biggl( \bigl(
\V_2(\partial P, h) \cap\O \bigr)\setminus\bigcup
_{i=1}^{\N} \cyl(A_i,h) \Biggr)
\nonumber
\\[-8pt]
\\[-8pt]
\nonumber
&&\qquad \leq\frac{1}{2} \| \vec\sigma\|_\infty\| \vec\nabla\varphi
_h\| _\infty\L^d \Biggl( \bigl(\V_2(
\partial P, h) \cap\O' \bigr)\setminus \bigcup
_{i=1}^{\N} \cyl(A_i,h) \Biggr)
\\
&&\qquad \leq C \eta\|\vec\sigma\|_\infty.\nonumber
\end{eqnarray}

\emph{Step \textup{2:} From $\int_{\cyl(A_i,h)} \vec\sigma
\cdot\vv_i
\,d\L^d$ to $\int_{\cyl(A_i,h)} \,d \vec\mu_n \cdot\vv_i$}.
As in Section~\ref{secconstr} if $\vec\mu_n \rightharpoonup\vec
\sigma\L
^d$ (up to extraction), since $\L^d (\partial\cyl(A_i,h)) = 0$, we
know by Portmanteau theorem that for all $i\in\{ 1,\ldots,\N\}$, for all
$n$ large enough (how large depending on $P,A_i,h,\eta$) we have
\[
\biggl\llvert \int_{\cyl(A_i,h)} \vec\sigma\cdot\vv_i \,d
\L^d - \int_{\cyl(A_i,h)} \,d\vec\mu_n \cdot
\vv_i \biggr\rrvert \leq\frac
{\eta
h}{ \N},
\]
and we conclude that for all $n$ large enough (how large depending on
$P,h,\eta$), we have
%
%
\begin{equation}
\label{eqmax3} \Biggl\llvert \sum_{i=1}
^{\N}\int_{\cyl(A_i,h)} \vec\sigma\cdot \vv_i \,d
\L^d - \sum_{i=1}^{\N} \int
_{\cyl(A_i,h)} \,d\vec\mu_n \cdot\vv_i \Biggr
\rrvert \leq\eta h.
\end{equation}

\emph{Step \textup{3:} From $\int_{\cyl(A_i,h)} \,d\vec\mu_n
\cdot\vv
_i$ to $\Psi(\vec\mu_n, \cyl(A_i,h), \vv_i)$}.
As in Section~\ref{secconstr} we use Proposition \ref{propflux} to
obtain that for all $i\in\{ 1,\ldots,\N\}$, there exists $\widetilde h_i
(P, \eta) >0$ such that for all $h\leq\widetilde h_i$, for all $n$
large enough (how large depending on $P,h,\eta$), we have
\[
\biggl\llvert \int_{\cyl(A_i,h)} \,d\vec\mu_n \cdot
\vv_i - \frac{2h
\Psi
(\vec\mu_n, \cyl(A_i,h), \vv_i)}{n^{d-1}} \biggr\rrvert \leq\frac
{\eta
}{ \sum_{i=1}^{\N} \H^{d-1} (A_i)}
\L^d \bigl(\cyl(A_i,h) \bigr).
\]
Thus there exists $h_3 (P, \eta) = \min_{1\leq i \leq\N}
(\widetilde
h_i) >0$ such that for all $h \leq h_3$, for all $n$ large enough (how
large depending on $P,\eta$), we have
%
%
\begin{equation}
\label{eqmax4} \Biggl\llvert \sum_{i=1}^{\N}
\int_{\cyl(A_i,h)} \,d\vec\mu_n \cdot \vv_i
- \frac{2h}{n^{d-1}} \sum_{i=1}^{\N} \Psi
\bigl(\vec\mu_n, \cyl (A_i,h), \vv _i \bigr)
\Biggr\rrvert \leq\eta h.
\end{equation}

\emph{Step \textup{4:} From $\sum_{i=1}^{\N} \Psi(\vec\mu_n,
\cyl
(A_i,h), \vv_i)$ to $\operatorname{flow}^{\mathrm{disc}}_n(\vec\mu_n)$}.
By construction of $P$ we know that $\G^1_n \subset P$ and $\G_n^2
\subset(\RR^d \setminus P)$ at least for $n$ large enough. Since
the stream $f_n$ satisfies the node law, we know that $\operatorname
{flow}^{\mathrm{disc}}_n(\vec\mu_n)$
is equal to the flow that goes out of $P$, that is,
\[
\operatorname{flow}^{\mathrm{disc}}_n(\vec\mu_n) =
\frac
{1}{n^{d-1}}\sum_{e=[a,b], a\in P, b
\notin P} f_n(e)
\vec{e}\cdot(n \overrightarrow{ab} ).
\]
Notice that $ \vec{e}\cdot( n \overrightarrow{ab})$ equals $+1$ or $-1$,
and $f_n(e) = 0$ if $e\notin\Pi_n$. We define
\[
E(P) = \bigl\{ e=[a,b] | e\in\Pi_n, a\in P, b \notin P \bigr\}.
\]
Thus
\[
\operatorname{flow}^{\mathrm{disc}}_n(\vec\mu_n) =
\frac
{1}{n^{d-1}}\sum_{e \in E(P)} f_n(e)
\vec{e} \cdot(n \overrightarrow{ab} ).
\]
For all $i\in\{ 1,\ldots,\N\}$, for all $h>0$, the set of edges
\[
E_i = \bigl\{ e\subset\cyl(A_i,h) | e\in E(P) \bigr\}
\]
is a cutset in $\cyl(A_i, h)$ from the lower half part of its boundary
to the upper half part of its boundary in the direction $\vv_i$; this
can be proved exactly as in Proposition~\ref{propflux1}. Thus
\[
\Psi \bigl(\vec\mu_n, \cyl(A_i,h), \vv_i
\bigr) = \sum_{e = [a,b]\in E_i,
a\in P, b\notin P} f_n(e) \vec{e}
\cdot(n \overrightarrow{ab}).
\]
Since the sets $E_i$ are disjoint, this implies that
\[
\Biggl\llvert \operatorname{flow}^{\mathrm{disc}}_n(\vec
\mu_n) - \frac{1}{n^{d-1}} \sum_{i=1}^{\N}
\Psi \bigl(\vec\mu _n, \cyl(A_i,h),\vv_i
\bigr) \Biggr\rrvert \leq\frac{M}{n^{d-1}} \card \Biggl(E(P) \setminus
\bigcup_{i=1}^{\N} E_i \Biggr).
\]
It remains to control $\card(E(P) \setminus\bigcup_{i=1}^{\N} E_i
)$. The edges that belong to this set are included in $\V_\infty
((\partial P
\cap\O') \setminus\bigcup_{i=1}^{\N} A_i, 2/n)$, thus
\begin{eqnarray*}
\card \Biggl(E(P) \setminus\bigcup_{i=1}^{\N}
E_i \Biggr)& \leq& 2d \card \Biggl(\V_\infty \Biggl( \bigl(
\partial P \cap\O' \bigr) \setminus\bigcup
_{i=1}^{\N} A_i, 2/n \Biggr) \cap
\ZZ^d_n \Biggr)
\\
& \leq&2d n^d \L^d \Biggl(\V_\infty \Biggl(
\overline{ \bigl(\partial P \cap\O' \bigr) \setminus \bigcup
_{i=1}^{\N} A_i}, 3/n \Biggr)
\Biggr).
\end{eqnarray*}
The set $ \overline{(\partial P \cap\O') \setminus\bigcup_{i=1}^{\N}
A_i}$ is a closed $(d-1)$-rectifiable set. Thus its $(d-1)$ dimensional
Minkowski content defined by
\[
\lim_{r\rightarrow0} \frac{1}{2r} \L^d \Biggl(
\V_2 \Biggl(\overline{ \bigl(\partial P \cap\O' \bigr)
\setminus\bigcup_{i=1}^{\N}
A_i}, r \Biggr) \Biggr)
\]
exists and is equal to $\H^{d-1}( \overline{(\partial P \cap\O')
\setminus\bigcup_{i=1}^{\N} A_i} )$ that is smaller than $\eta$ by
construction. Thus there exists a constant $\kappa(d)$ such that, for
$n$ large enough,
\[
\card \Biggl(E(P) \setminus\bigcup_{i=1}^{\N}
E_i \Biggr) \leq \kappa\eta n^{d-1}.
\]
For all $n$ large enough we get
%
%
\begin{equation}
\label{eqmax5} \Biggl\llvert \operatorname{flow}^{\mathrm{disc}}_n(
\vec\mu_n) - \frac{1}{n^{d-1}} \sum_{i=1}^{\N}
\Psi \bigl(\vec \mu_n, \cyl(A_i,h),\vv_i
\bigr) \Biggr\rrvert \leq\kappa\eta M.
\end{equation}

\emph{Step \textup{5:} Conclusion}.
Combining inequalities (\ref{eqmax12}), (\ref{eqmax3}), (\ref{eqmax4})
and (\ref{eqmax5}) for a $h\leq\min(h_1, h_2,h_3)$, we obtain that
for all $n$ large enough (how large depending on everything else)
\[
\bigl\llvert \operatorname{flow}^{\mathrm{cont}}(\vec\sigma) -
\operatorname{flow}^{\mathrm{disc}}_n(\vec\mu_n) \bigr
\rrvert \leq \eta \bigl(C \|\vec\sigma\|_\infty+ 1 + \kappa M\bigr),
\]
and this completes the proof of Proposition \ref{propmax2}.
\end{pf*}

\begin{pf*}{Proof of Proposition \ref{propAi}}
Let $\eta>0$. The sets $P$ and $ \O'$ are polyhedral; that is, their
boundaries are included in a finite number of hyperplanes. For any
$x\in\partial P \cap\partial\O'$, let us denote by $\theta(x)\in
[0,\pi]$ the
angle between the exterior unit normal to $\partial P$ at $x$ and the
exterior unit normal to $\partial\O'$ at $x$. Thus
\[
\theta_1 = \inf \bigl\{ \theta(x) | x \in\partial P \cap\partial
\O' \bigr\} > 0,
\]
since there are only finitely many different values of $\theta(x)$
which are all positive because $\partial P$ is transversal to
$\partial\O'$. We
denote by $(F_l, l=1,\ldots, \M)$ the faces of $\partial P$ that
intersects $\O
'$, thus $\partial P \cap\O' = \bigcup_{l=1}^{\M} F_l$, and by
$\vv_l$ the
exterior unit vector normal to $P$ along $F_l$. We define
\[
\theta_2 = \min \bigl\{ \arccos(\vv_l \cdot
\vv_m) | l,m=1,\ldots,\M, l\neq m, F_l \cap
F_m \neq\varnothing \bigr\},
\]
the minimum of the angles between two adjacent faces of $\partial P
\cap\O
'$, that is, between faces that intersect. Thus $\theta_2 > 0$ since
again there are finitely many such angles. Let $\theta_0 = \min
(\theta
_1, \theta_2) >0$; see Figure~\ref{figcylAi_2}.
%
%
%
Let $E$ be the set of the edges of $\partial P \cap\O'$, that is,
\[
E = \bigcup_{l\neq m\leq\M} F_l \cap
F_m.
\]
\begin{figure}

\includegraphics{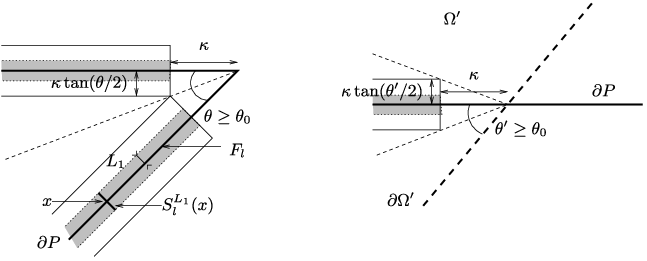}

\caption{Construction of the cylinders $(\cyl(A_i,h), i =1,\ldots,
\N)$.}
\label{figcylAi_2}
\end{figure}
There exists $\kappa>0$ small enough (how small depending on $\eta,
\O
', P$) so that we have
\[
\H^{d-1} \bigl( \partial P \cap\O' \cap\V_2
(E, \kappa) \bigr) \leq \frac{\eta}{2}.
\]
Let $F_l' = F_l \setminus\V_2 (E,\kappa)$ for $l\in\{1,\ldots
,\M\}
$. We define $L_1 =2^{-1} \kappa\tan(\theta_0 /2)$. By definition of
$\theta_0$ and $F_l'$, for all $l\in\{ 1,\ldots,\M\}$, for all $x$ in
$F_l'$, the set
\[
S_l^{L_1} (x) = \{ x + b \vv_l | -
L_1 \leq b\leq L_1 \}
\]
is included in $\O'$ and does not intersect $S^{L_1}_m (y)$ for any
$y\in F'_m$ and any $m$ such that $F_l$ and $F_m$ are adjacent; see
Figure~\ref{figcylAi_2}. Let $L_2$ be the infimum of the distances
between two nonadjacent faces of $\partial P $ (thus $L_2 >0$). Let
$L_0 =
\min(L_1, L_2/2) >0$. Then for any $l\in\{ 1,\ldots,\M\}$, for any $x
\in F'_l $, for any $z$ in $S_l^{L_0} (x) =\{ x + b \vv_l | - L_0
\leq b\leq L_0 \} $, we have $z \in\O'$, $d_2(z, \partial P
\setminus
F_l) > L_0 \geq d_2(z, F_l)$ and $z\notin S_m^{L_0} (y)$ for any $y \in
F_m'$ distinct from $x$ and any $m\in\{1,\ldots,\M\}$. We can now cover
$\bigcup_{l=1}^{\M} F_l'$ by a finite set of hyperrectangles $(A_i,
i=1,\ldots,\N)$ depending on $\O',P,\eta$ of disjoint interiors up
to a
surface of $\H^{d-1}$-measure less than $\eta/2$, that is,
\[
\H^{d-1} \Biggl( \bigcup_{l=1}^{\M}
F_l' \setminus\bigcup
_{i=1}^{\N} A_i \Biggr) \leq
\frac{\eta}{2}.
\]
This implies that
\[
\H^{d-1} \Biggl( \bigl(\partial P \cap\O' \bigr)
\setminus\bigcup_{i=1}^{\N}
A_i \Biggr) \leq\eta.
\]
Let us consider the cylinders $\cyl(A_i, h)$ for $h \leq L_0 (\O', P,
\eta)$, $i=1,\ldots,\N$. By construction, for all $i\in\{ 1,\ldots
,\N\}$,
$\cyl(A_i, h) \subset\O'$ and for all $x\in\cyl(A_i,h)$, $d_2 (x,
\partial
P \setminus A_i) > d_2 (x, A_i)$. To complete the proof, it
remains to control $\L^d  ( (\V_2 (\partial P, h) \cap\O')
\setminus\bigcup_{i=1}^{\N} \cyl(A_i,h)  )$. We remark that
%
%
\begin{eqnarray}
\label{eqlast1} \bigl(\V_2 (\partial P, h) \cap\O'
\bigr) \setminus\bigcup_{i=1}^{\N
}
\cyl(A_i,h) & \subset& \Biggl( \bigl(\V_2 (\partial P, h)
\cap\O' \bigr) \setminus\bigcup_{l=1}^{\M}
\cyl \bigl(F'_l,h \bigr) \Biggr)
\nonumber
\\[-8pt]
\\[-8pt]
\nonumber
&&{} \cup \Biggl( \bigcup_{l=1}^{\M}
\cyl \bigl(F'_l,h \bigr) \setminus\bigcup
_{i=1}^{\N} \cyl(A_i,h) \Biggr)
.
\end{eqnarray}
On one hand,
%
%
\begin{equation}
\label{eqlast2} \L^d \Biggl( \bigcup_{l=1}^{\M}
\cyl \bigl(F'_l,h \bigr) \setminus\bigcup
_{i=1}^{\N} \cyl(A_i,h) \Biggr)
\le2 \frac{\eta}{2} h.
\end{equation}
On the other hand,
\begin{eqnarray*}
\Biggl( \bigl(\V_2 (\partial P, h) \cap\O' \bigr)
\setminus\bigcup_{l=1}^{\M} \cyl
\bigl(F'_l,h \bigr) \Biggr)& \subset& \Biggl( \bigcup
_{l=1}^{\M} \cyl \bigl(F_l
\setminus F'_l, h \bigr) \Biggr)
\nonumber
\\[-8pt]
\\[-8pt]
\nonumber
&&{} \cup\V_2 (E, h) \cup\V_2 \bigl( \partial P \cap
\partial\O', h / \sin(\theta_0 /2) \bigr);
\end{eqnarray*}
see Figure~\ref{figcylAi_3}.
%
%
Thus
%
%
\begin{eqnarray}
\label{eqlastbis1} &&\L^d \Biggl( \bigl(\V_2 (\partial P,
h) \cap\O' \bigr) \setminus\bigcup
_{l=1}^{\M
} \cyl \bigl(F'_l,h
\bigr) \Biggr)
\nonumber
\\[-8pt]
\\[-8pt]
\nonumber
&&\qquad \leq2 \frac{\eta}{2} h + \L^d \bigl(
\V_2 (E, h) \bigr) + \L^d \bigl( \V_2 \bigl(
\partial P \cap\partial\O', h / \sin(\theta_0 /2)
\bigr) \bigr).
\end{eqnarray}
%
\begin{figure}

\includegraphics{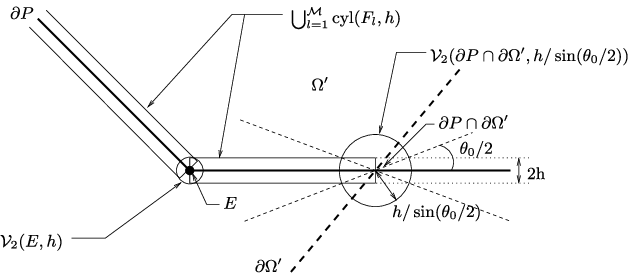}

\caption{Near $\partial P \cap\partial\O'$.}
\label{figcylAi_3}
\end{figure}

The sets $E$ and $\partial P \cap\partial\O'$ are finite unions of
$(d-2)$-closed
rectifiable subsets, whose $(d-2)$ dimensional Minkowski contents are
equal to their $\H^{d-2}$-measure, thus
\[
\limsup_{h\rightarrow0} \frac{\L^d ( \V_2 (E, h) ) }{\alpha_2 h^2} \leq\H^{d-2} (E),
\]
and we conclude that there exists $L_0' (\O',P)$ such that if $h\leq
L_0'$, we have
%
%
\begin{equation}
\label{eqlastbis2} \L^d \bigl( \V_2 (E, h) \bigr) \leq2
\alpha_2 h^2 \H^{d-2} (E).
\end{equation}
In the same way, we obtain that there exists $L_0'' (\O', P)$ such that
for all $h\leq L_0''$,
%
%
\begin{eqnarray}
\label{eqlastbis3}&& \L^d \bigl( \V_2 \bigl(\partial P
\cap\partial\O', h / \sin(\theta_0 /2) \bigr) \bigr)
\nonumber
\\[-8pt]
\\[-8pt]
\nonumber
&&\qquad\leq2 \alpha_2 h^2 \sin^{-2} (
\theta_0 /2) \H^{d-2} \bigl(\partial P \cap \partial
\O' \bigr).
\end{eqnarray}
Combining inequalities (\ref{eqlastbis1}), (\ref{eqlastbis2}) and
(\ref
{eqlastbis3}), we obtain that for $h \leq \min(L_0',L_0'')$,
\begin{eqnarray*}
&&\L^d \Biggl( \bigl(\V_2 (\partial P, h) \cap
\O' \bigr) \setminus\bigcup_{l=1}^{\M
}
\cyl \bigl(F'_l,h \bigr) \Biggr)\\
&&\qquad \leq2 \biggl[
\frac{1}{2} + \frac{h}{\eta} \alpha_2 \bigl(
\H^{d-2} (E)+ \sin^{-2} (\theta_0 /2)
\H^{d-2} \bigl(\partial P \cap\partial \O' \bigr) \bigr)
\biggr] \eta h.
\end{eqnarray*}
If $h \leq\eta$, we obtain
%
%
\begin{equation}
\label{eqlast3} \L^d \Biggl( \bigl(\V_2 (\partial P, h)
\cap\O' \bigr) \setminus\bigcup_{l=1}^{\M
}
\cyl \bigl(F'_l,h \bigr) \Biggr) \leq2 C\eta h,
\end{equation}
where $C = C(\O',P)$ is a constant depending on $\O',P$. Combining
inequalities~(\ref{eqlast1}), (\ref{eqlast2}) and (\ref{eqlast3}), we
obtain that for $h \leq(L_0',L_0'',\eta)$,
\[
\L^d \Biggl( \bigl( \V_2 (\partial P,h) \cap
\O' \bigr) \setminus\bigcup_{i=1}^{\N}
\cyl(A_i,h) \Biggr) \leq2 \bigl(1 + C \bigl(\O',P
\bigr) \bigr) \eta h.
\]
Finally, we fix $h_1 (\O',P, \eta) = h_1 (\O,P,\eta',\eta) = \min
(L_0,L_0',L_0'',\eta)$, and Proposition \ref{propAi} is proved.
\end{pf*}

\begin{rem}
In the proof of Proposition \ref{propmax2}, we could use a weaker
version of Proposition \ref{propAi} without defining the set $\O'$, and
with cylinders $\cyl(A_i, h)$ that almost cover $\V_2 (\partial P,
h)$ even
outside $ \O$. This weaker version of Proposition~\ref{propAi} would be
easier to prove, as we would not need to construct a set $P$ whose
boundary is transversal to $\G$, and then a set $\O'$. However, we will
use again Proposition~\ref{propAi} and its consequences in Section~\ref{secinegmfmn}, and at that point we will need Proposition \ref{propAi}
as it is stated.
\end{rem}


\section{Study of minimal cutsets}
\label{secmincut}

The study of the asymptotic behavior of minimal cutsets was almost done
in \cite{CerfTheret09infc}. However, it was not the goal of that
article to get information on minimal cutsets; thus the pieces of the
puzzle were not put together. This is what we do in this section. We
will not rewrite all the proofs, but we explain how to adapt them.

From now on, $(\E_n)_{n\geq1}$ denotes a sequence of $(\G^1_n, \G
^2_n)$-cutsets in $\O_n$, and $(\E_n^{\min})_{n\geq1}$ a sequence of
minimal $(\G^1_n, \G^2_n)$-cutsets in $\O_n$. We define as in Section~\ref{secdef1} the sets
\[
r(\E_n) = \bigl\{ x\in\O_n | \mbox{there exists a path
from } x \mbox{ to } \G^1_n \mbox{ in } \bigl(
\ZZ^d_n, \Pi_N \setminus
\E_n \bigr) \bigr\}
\]
and
\[
R(\E_n) = r(\E_n) + \frac{1}{2n}[-1,
1]^d,
\]
and we introduce the notation
\[
E_n = R(\E_n) \cap\O
\]
[the same definitions hold for $R(\E_n^{\min})$, $E_n^{\min}$]. We
recall that throughout the proofs, we suppose that the hypotheses \ref{hypo1} and \ref{hypo2} are fulfilled.


\subsection{\texorpdfstring{Restriction to $\O$}
{Restriction to Omega}}
\label{secEmin}

We prove that it is completely equivalent to study the convergence of
$(R(\E_n))_{n\geq1}$ or the convergence of $(E_n)_{n\geq1}$:
%
\begin{prop}
\label{propE}
Let $(\E_n)_{n\geq1}$ be a sequence of admissible discrete\break  $(\G^1_n,
\G
^2_n)$-cutsets in $\O_n$. We have
\[
\lim_{n\rightarrow\infty} \mathfrak d \bigl(E_n, R(
\E_n) \bigr) = 0.
\]
\end{prop}
%
\begin{rem}
This proposition implies that a subsequence of $(R(\E_n))_{n\geq1}$ is
convergent if and only if the corresponding subsequence of
$(E_n)_{n\geq1}$ is convergent, in which case they have the same
limit. Thus we can study the sequence $(E_n)_{n\geq1}$ instead of
$(R(\E_n))_{n\geq1}$.
\end{rem}

\begin{pf*}{Proof of Proposition \ref{propE}} For every $n\geq1$,
\[
\mathfrak d \bigl(R(\E_n), E_n \bigr) \leq
\L^d \bigl(\V_\infty(\O, 1/n) \setminus\O \bigr) \leq
\L^d \bigl(\V_{\infty} (\G, 1/n) \bigr).
\]
Since $\G$ is piecewise of class $\C^1$, $\G$ is a closed
$(d-1)$-rectifiable subset of $\RR^d$. Thus its $(d-1)$ dimensional
Minkowski content defined by
\[
\lim_{r\rightarrow0} \frac{1}{2r} \L^d \bigl(
\V_2 (\G, r) \bigr)
\]
exists and is equal to $\H^{d-1}(\G)$; see, for example, Appendix A in
\cite{Cerf-Pisztora}. This implies that
\[
\lim_{n \rightarrow\infty} \mathfrak d \bigl(R(\E_n),
E_n \bigr) \leq\lim_{n\rightarrow\infty} \L^d \bigl(
\V_\infty(\G, 1/n) \bigr) = 0.
\]
\upqed\end{pf*}


\subsection{Compactness}
\label{seccompmincut}

We prove the following result:
%
\begin{prop}
\label{propcompmincut}
We suppose that hypothesis \ref{hypo3} is also fulfilled. Let $(\E
_n^{\min})_{n\geq1}$ be a sequence of minimal discrete $(\G^1_n, \G
^2_n)$-cutsets in $\O_n$. Almost surely, for $n$ large enough, the
sequence $(E_n^{\min})_{n\geq1}$ takes its values in a deterministic
$\mathfrak d$ compact set that is included in $\{ F \subset\O|
\mathbh{1}_F \in \operatorname{BV}(\O) \}$.
\end{prop}
%
\begin{rem}
The previous proposition implies that a.s., any subsequence of
$(E_n^{\min})_{n\geq1}$ [thus of $(R(\E_n^{\min}))_{n\geq1}$] admits
a sub-subsequence which is convergent for the distance $\mathfrak d$,
and its limit $F$ is a subset of $\O$ that satisfies $\mathbh{1}_F
\in
\operatorname{BV} (\O)$.
\end{rem}
%
\begin{rem}
In the previous proposition, hypothesis \ref{hypo2} could be replaced by
the hypothesis that $\Lambda$ admits an exponential moment
\[
\exists\theta>0 \qquad \int_{\RR^+} e^{\theta x} \,d\Lambda(x) <
+ \infty.
\]
\end{rem}

\begin{pf*}{Proof of Proposition \ref{propcompmincut}} We
study the
sequence $(E_n^{\min})_{n\geq1}$ exactly as in \cite
{CerfTheret09infc}, Section~4. According to Theorem $1$ in \cite
{Zhang07}, adapted to our case as said in Remark $2$ in \cite{Zhang07},
we know that:
%
\begin{thmm}[(Zhang)]
\label{thmZhang}
If the law of the capacity of the edges admits an exponential moment,
and if hypothesis \ref{hypo3} is fulfilled, then there exist constants
$\beta_0 = \beta_0
(\Lambda,d)$, $C_i = C_i(\Lambda,d)$ for $i=1,2$ and $N=N(\Lambda,
d,\O,\G, \G^1, \G^2 )$ such that for all $\beta\geq
\beta_0$, for all $n\geq N$, we have
\[
\PP \bigl[\card \bigl(\E_n^{\min} \bigr) \geq\beta
n^{d-1} \bigr] \leq C_1 \exp \bigl(-C_2 \beta
n^{d-1} \bigr).
\]
\end{thmm}
%
\begin{rem}
The adaptation of Zhang's result in our setting involves one
difficulty: the cutsets we have to consider may not be connected.
However, we can get around this problem by considering the union of a
set $\E_n^{\min}$ with the edges that lie along $\G$: it is always
connected, and the number of edges we have added is bounded by $c
n^{d-1}$ for a constant $c$ depending only on the domain $\O$, since
$\G
$ is piecewise of class $\C^1$. Then the adaptation of Zhang's proof is
straightforward.
\end{rem}

If the capacities are bounded, their law admits an exponential moment.
Thus we can use Theorem \ref{thmZhang}. We obtain
\[
\sum_{n\geq1} \PP \bigl(\card \bigl(
\E_n^{\min} \bigr) \geq\beta_0 n^{d-1}
\bigr) \leq N + C_1 \exp \bigl(-C_2 \beta_0
n^{d-1} \bigr) < +\infty,
\]
and thus by the Borel--Cantelli lemma,
\[
\PP \Bigl( \limsup_{n\rightarrow\infty} \bigl\{ \card \bigl(
\E_n^{\min} \bigr) \geq \beta_0 n^{d-1}
\bigr\} \Bigr) = 0;
\]
that is, a.s. there exists $n_0$ such that for all $n\geq n_0$, $\card
(\E_n^{\min}) < \beta_0 n^{d-1}$. For $F\subset\RR^d$, we recall that
the perimeter of $F$ in $\O$ is
defined by
\begin{eqnarray*}
&&\P(F, \O)\\
&&\qquad = \sup \biggl\{\int_F \di\vec f(x) \,d
\L^d(x) | \vec f\in \C_c^\infty \bigl(\O,
\RR^d \bigr), \vec f (x) \in B(0,1) \mbox{ for all } x \in\O \biggr\}
.
\end{eqnarray*}
If $\card(\E_n^{\min}) \leq\beta_0 n^{d-1} $, then $\P(E_n^{\min
}, \O
)\leq\beta_0$. We define
\[
\C_{\beta_0} = \bigl\{ F\subset\O| \P(F, \O) \leq\beta_0
\bigr\}.
\]
Thus we have proved that
\[
a.s.\ \exists n_0\ \forall n\geq n_0\qquad
E_n^{\min} \in\C _{\beta_0}.
\]
We endow $\C_{\beta_0}$ with the pseudo-metric associated to the
distance $\mathfrak d$. Remember that $\mathfrak d (F, F') = \L^d(F
\,\triangle\, F')$, where $\triangle$ is the symmetric difference. For this
metric the set $\C_{\beta_0}$ is compact. Moreover $\C_{\beta_0}
\subset\{ F \subset\O| \mathbh{1}_F \in \operatorname{BV}(\O) \}$. This ends
the proof of Proposition \ref{propcompmincut}.
\end{pf*}


\subsection{Minimality}
\label{secmin}

We recall that for a set of edges $\E_n \subset\EE^d_n$,
\[
V(\E_n) = \sum_{e\in\E_n} t(e),
\]
and that for $F \subset\O$ of finite perimeter,
\begin{eqnarray*}
\capa(F) & =& \int_{\O\cap\partial^* F} \nu \bigl(\vv_F (x)
\bigr) \,d \H^{d-1} (x) + \int_{\G^2 \cap\partial^* F} \nu \bigl(
\vv_F (x) \bigr) \,d \H^{d-1} (x)
\\
&&{} + \int_{\G^1 \cap\partial^* (\O\setminus F)} \nu \bigl(\vv _{\O} (x) \bigr) \,d
\H^{d-1} (x).
\end{eqnarray*}
To complete the proof of Theorem \ref{thmpcpal2}, we must prove that
the random limit $F$ of a subsequence of minimal discrete cutsets $(\E
_n^{\min})_{n\geq1}$ satisfies
%
%
\begin{equation}
\label{eqminF} \capa(F) = \phi_{\O}^a\qquad \mbox{a.s.}
\end{equation}
In this section, we prove the following result:
%
\begin{prop}
\label{propmin}
Let $(\E_n)_{n\geq1}$ be a sequence of admissible discrete\break  $(\G^1_n,
\G
^2_n)$-cutsets in $\O_n$. If a subsequence $(R(\E_{\varphi
(n)}))_{n\geq1}$ converges for the distance $\mathfrak d$ toward a set
$F\subset\O$ of finite perimeter in $\O$, then almost surely
\[
\liminf_{n\rightarrow\infty} \frac{V(\E_{\varphi(n)})}{\varphi
(n)^{d-1}} \geq\capa(F).
\]
\end{prop}
%
\begin{rem}
As for the maximal streams, we will deduce equation (\ref{eqminF}) from
Proposition \ref{propmin} in Section~\ref{secccl}, using our study of
maximal streams and Lemma \ref{lemineg}.
\end{rem}

\begin{pf*}{Proof of Proposition \ref{propmin}} The idea of the
proof is the following. We almost cover $\partial F$ by a finite set of
disjoint balls $(B_i = B(x_i, r_i))$, small enough so that $\partial
F$ is
almost flat in each ball. ``Almost flat'' means that:
\begin{longlist}[(a)]
\item[(a)] the surface $\partial F \cap B_i$ is very close to the
flat disc
$\disc(x_i, r_i, \vv_i)$ where $\vv_i = \vv_F (x_i)$;
\item[(b)] $F \cap B_i$ is very close for the distance $\mathfrak d$ to
$B^-(x_i, r_i, \vv_i)$.
\end{longlist}
From (a) we deduce that $\capa(F)$ is very close to $\sum_i \alpha
_{d-1} r_i^{d-1} \nu(\vv_i)$, the sum of the capacities of the discs
$\disc(x_i, r_i, \vv_i)$. Since the balls are disjoint we get $V(\E_n)
\geq\sum_i V(\E_n \cap B_i)$, where $\E_n \cap B_i = \{ e\in\E_n
| e\subset B_i \}$. It remains to compare in any ball $B = B(x, r \vv
) \in(B_i)$ the quantities $V(\E_n \cap B)$ and $\alpha_{d-1} r^{d-1}
\nu(\vv)$. Since $\mathfrak d (E_n, F)$ goes to zero, by (b) we can
suppose that for large $n$, $E_n \cap B$ is very close to $B^- (x,r,\vv
)$. We can construct a cutset in $B$ from the upper half part of its
boundary to its lower half part by adding not too much edges to $\E_n
\cap B$---this is the difficult part. Thus $V(\E_n \cap B) \geq\tau_B$
up to an error term, where $\tau_B$ is, roughly speaking, the maximal
flow in $B$ from the upper half part of its boundary to its lower half
part. Using the known law of large numbers for the maximal flows~$\tau
$, we can prove that $\tau_B$ is equivalent to $\alpha_{d-1} r^{d-1}
\nu(\vv) n^{d-1}$ for large $n$, and this completes the proof.
\begin{figure}

\includegraphics{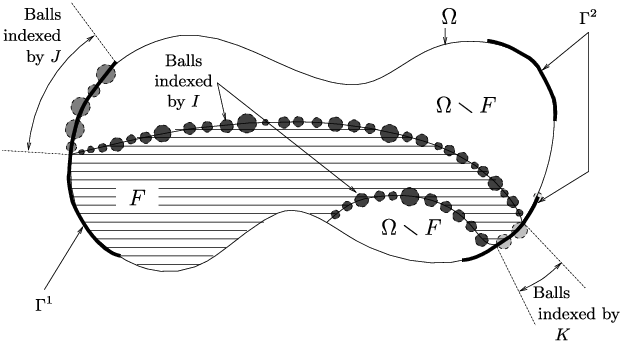}

\caption{Covering of $(\partial F \cap\O)\cup(\partial F \cap\G
^2) \cup(\partial(\O
\setminus F ) \cap\G^1)$ by balls.}
\label{figrecouvrement}
\end{figure}

We consider a subsequence of $(R(\E_n))_{n\geq}$ that converges toward
$F$, but we still denote it by $(R(\E_n))_{n\geq1})$ for simplicity.
If $\capa(F) =0$, there is nothing to prove. Thus we can suppose that
$\capa(F) >0$. In fact it has been proved in~\cite{CerfTheret09geoc}
that under hypotheses \ref{hypo1} and \ref{hypo2}, $\phi_{\O}^b > 0$ if
and only if hypothesis \ref{hypo3} is fulfilled. Thus it is indeed the
case that $\capa(F) >0$.

\emph{Step \textup{1:} From $V(\E_n)$ to $\sum_i V(\E_n \cap B_i)$
and from $\capa(F)$ to\break  $\sum_i \alpha_{d-1} r_i^{d-1} \times \nu(\vv_i)$}.
We consider a fixed realization of the capacities. We use Lemma 1 in
\cite{CerfTheret09infc} to cover $\partial F$ by a set of balls well chosen;
see Figure~\ref{figrecouvrement}:

%
\begin{lem}[(Lemma 1 in \cite{CerfTheret09infc})]
\label{chapitre7covering}
Let $F $ be a subset of $\O$ of finite perimeter. For every positive
constant $\delta$ and $\eta$, there exists a finite family of closed
disjoint balls
$(B_i)_{i\in I\cup J\cup K} $ where $B_i = (B(x_i, r_i), \vv_i)$
such that (the vector $\vv_i$ defines $B_i^-$)
\begin{eqnarray*} &&\forall i \in I, x_i \in\partial^*F
\cap\O, \vv_i = \vv_{F}(x_i), r_i
\in\,]0,1 \bigl[, B_i \subset\O, \mathfrak d \bigl(F\cap
B_i, B_i^- \bigr) \leq\delta\alpha_d
r_i^d,
\\
&&\forall i\in J, x_i \in\G^{1}\cap\partial^*(\O
\setminus F),\vv_i = \vv_{\O}(x_i),
r_i \in\, \bigr]0,1 \bigl[, \partial\O\cap B_i \subset
\G^1, \\
&&\qquad\mathfrak d \bigl( (\O\setminus F ) \cap B_i,
B_i^- \bigr) \leq \delta \alpha_d r_i^d
,
\\
&&\forall i \in K, x_i \in\G^{2} \cap\partial^* F,
\vv_i = \vv _{F}(x_i), r_i \in\,
\bigr]0,1[, \partial\O\cap B_i \subset\G^2,\\
&&\qquad \mathfrak d
\bigl(F\cap B_i, B_i^- \bigr) \leq\delta
\alpha_d r_i^d,
\end{eqnarray*}
and finally
\[
\biggl\llvert \capa(F) - \sum_{i \in I \cup J \cup K}
\alpha_{d-1} r_i^{d-1} \nu(\vv_i) \biggr
\rrvert \leq\eta.
\]
\end{lem}
%

%
%
We do not give the proof of Lemma \ref{chapitre7covering} here. It
relies on the Vitali covering theorem for $\H^{d-1}$ and the
Besicovitch differentiation theorem in $\RR^d$.
%
\begin{rem}
In fact, Lemma 1 in \cite{CerfTheret09infc} states the condition
$\mathfrak d (B_i\cap\O, B_i^-) \leq\delta$ instead of $
\mathfrak d (B_i\cap(\O\setminus F ), B_i^-) \leq\delta$
for $i \in J$. Both statements are true, since we can apply the same
techniques to $\O$ or $\O\setminus F$. However, Lemma 1 in \cite
{CerfTheret09infc} should have been written as Lemma \ref
{chapitre7covering} here since the property actually used in Section~5.2 of \cite{CerfTheret09infc} is in fact $ \mathfrak d (B_i\cap(\O
\setminus F ), B_i^-) \leq\delta$; there is a small mistake
in this section on page 653.
\end{rem}
We need to move these balls a little bit to obtain balls whose centers
have rational coordinates and with $\vv_i$ a rational direction. Let
$0<\eta, \delta<1$, and let $(B_i)_{i\in I\cup J \cup K}$ be the
family of balls associated to $(\eta/2, \delta/2)$. The function $\nu$
is continuous on $\SS^{d-1}$. Thus there exists $\theta_0>0$ such that
for all vectors $\vv, \vec w \in\SS^{d-1}$,
\[
\vv\cdot\vec w \geq\cos\theta_0\quad \Longrightarrow\quad\bigl|\nu(\vv) - \nu(\vec
w)\bigr| \leq\frac{\eta\nu_{\min}}{4 \capa(F)},
\]
where $\nu_{\min} = \inf_{\SS^{d-1}} \nu>0$. If for all $i \in
I\cup J
\cup K$, $\vv_i \cdot\vv_i ' \geq\cos\theta_0$, we get
\begin{eqnarray*}
&&\biggl\llvert \sum_{i\in I \cup J \cup K} \alpha_{d-1}
r_i^{d-1} \nu (\vv _i) - \sum
_{i\in I \cup J \cup K} \alpha_{d-1} r_i^{d-1}
\nu \bigl(\vv_i' \bigr) \biggr\rrvert\\
&&\qquad \leq
\frac{\eta\nu_{\min}}{4 \capa(F)} \sum_{i\in I
\cup J \cup K} \alpha_{d-1}
r_i^{d-1}
\\
&&\qquad \leq\frac{\eta\nu
_{\min
}}{4 \capa(F)} \frac{2 \capa(F)}{\nu_{\min}}
\\
&&\qquad \leq\frac{
\eta
}{2}.
\end{eqnarray*}
Moreover, for all $x, x', r, \vv, \vv'$ with $\vv\cdot\vv' = \cos
\theta$ (see Figure~\ref{figdeplacementboule}), we have
%
\begin{figure}

\includegraphics{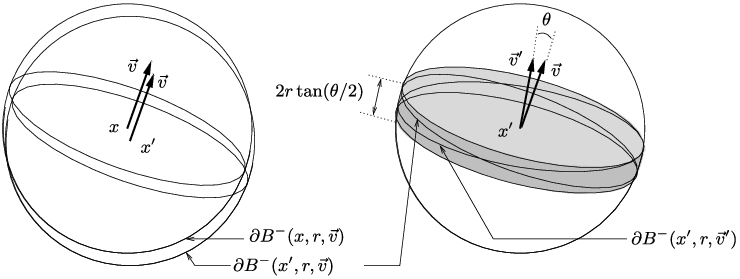}

\caption{Comparison between $B^- (x,r,\vv)$ and $B^- (x', r, \vv')$.}
\label{figdeplacementboule}
\end{figure}
%
%
%
\begin{eqnarray*}
\label{eqessai1}&& \L^d \bigl( B^-(x,r,\vv) \triangle B^-
\bigl(x',r,\vv' \bigr) \bigr) \\
&&\qquad \leq\L^d
\bigl( B^-(x,r,\vv) \triangle B^- \bigl(x',r,\vv \bigr) \bigr)+
\L^d \bigl( B^- \bigl(x',r,\vv \bigr) \triangle B^-
\bigl(x',r,\vv' \bigr) \bigr)
\\
&&\qquad \leq\L^d \bigl(\V_2 \bigl(\partial B^- (x,r,\vv),\bigl \| x
- x' \bigr\| \bigr) \bigr)
\\
&&\qquad\quad{} + \L^d \bigl(\cyl \bigl(\disc \bigl(x',r, \bigl(\vv+
\vv' \bigr)/2 \bigr), r \tan(\theta/2) \bigr) \bigr)
\\
&&\qquad \leq\L^d \bigl(\V_2 \bigl(\partial B^- (0,r,\vv), \bigl\| x
- x'\bigr \| \bigr) \bigr) + \alpha_{d-1} r^{d-1}
\times2 r \tan(\theta/2)
\\
&&\qquad \leq\frac{\delta}{2} \alpha_{d} r^{d},
\end{eqnarray*}
where the last inequality is valid as soon as
\[
\tan(\theta/2) \leq\frac{\delta\alpha_d}{8 \alpha_{d-1}}
\]
and
%
%
\begin{equation}
\label{eqcons2} \L^d \bigl(\V_2 \bigl(\partial B^- (0,r,
\vv),\bigl \| x - x' \bigr\| \bigr) \bigr) \leq\frac
{\delta}{4}
\alpha_{d} r^{d}.
\end{equation}
We know that $\partial B^{-} (0,r,\vv)$ is $(d-1)$-rectifiable. Thus its
$(d-1)$ dimensional Minkowski content exists and
\[
\lim_{R \rightarrow0} \frac{\L^d (\V_2 ( \partial B^{-} (0,r,\vv
),R))}{2R} = \H^{d-1} \bigl(
\partial B^{-} (0,r,\vv) \bigr) = K r^{d-1}
\]
for a constant $K$ depending only on the dimension. Thus for $x'$ close
enough to $x$,
\[
\L^d \bigl(\V_2 \bigl(\partial B^- (0,r,\vv), \bigl\| x -
x' \bigr\| \bigr) \bigr) \leq4 \bigl\| x - x'\bigr\| K
r^{d-1},
\]
and we obtain (\ref{eqcons2}) for $\| x' - x\| $ small enough (how
small depending on $d$, $r$ and $\delta$). Thus there exists $\theta_1
>0$ such that if $\vv_i \cdot\vv_i' \geq\cos\theta_1$ for all $i
\in
I\cup J \cup K$, and if $x_i'$ is close enough to $x_i$, we have
\[
\L^d \bigl( B^-(x_i,r_i,\vv_i)
\triangle B^- \bigl(x'_i,r_i,
\vv_i' \bigr) \bigr) \leq \frac{\delta}{2}
\alpha_{d} r_i^{d}.
\]
Since $\O$ is open and $\G^1$ and $\G^2$ are open in $\G$, we can
choose for all $i\in I\cup J \cup K$ a unit vector $\vv_i'$ that
defines a rational direction, that is, such that $\lambda_i \vv_i '$
has rational coordinates for a positive real number $\lambda_i$, and a
point $x_i'$ that has rational coordinates, such that
\begin{eqnarray*} &&\forall i \in I, r_i \in\,]0,1 [
, B \bigl(x'_i, r_i \bigr) \subset\O,
\mathfrak d \bigl(F\cap B \bigl(x_i', r_i
\bigr), B^- \bigl(x_i', r_i,
\vv_i' \bigr) \bigr) \leq\delta\alpha_d
r_i^d,
\\
&&\forall i\in J, r_i \in\, ]0,1 [, \partial\O\cap B
\bigl(x_i', r_i \bigr) \subset
\G^1, \mathfrak d \bigl( (\O\setminus F ) \cap B
\bigl(x_i',r_i \bigr),\\
&&\hspace*{193pt}\qquad B^-
\bigl(x_i', r_i, \vv_i'
\bigr) \bigr) \leq\delta \alpha_d r_i^d,
\\
&&\forall i \in K, r_i \in\, ]0,1[, \partial\O\cap B
\bigl(x_i', r_i \bigr) \subset
\G^2, \mathfrak d \bigl(F\cap B \bigl(x_i',
r_i \bigr), B^- \bigl(x_i',
r_i, \vv_i' \bigr) \bigr) \leq \delta
\alpha_d r_i^d.
\end{eqnarray*}
For simplicity of notation, we skip the prime and still denote this new
family of balls associated to $(\eta, \delta)$ by \mbox{$(B_i, i \in I\cup J
\cup K) = (B(x_i, r_i), \vv_i, i \in I \cup J \cup K)$}.

\begin{rem}
If $\mathfrak d (F\cap B(x,r), B^-(x,r,\vv))$ is small, $F$ ``looks
like''\break  $B^-(x,r,\vv)$ inside the ball $B(x,r)$. This means that the
volume of $(F\cap B(x,r))\times  \triangle B^-(x,r,\vv)$ is small; however $F$
(resp., $F^c$) might have some thin strands (of small volume, but that
can be long) that go deeply into $B^+(x,r,\vv)$ [resp., $B^-(x,r,\vv)$];
see Figure~\ref{figtrompes}. If $d\geq3$, this is not in contradiction
with the hypothesis that the capacity of $\partial F$ inside
$B(x,r,\vv)$ is
close to $\alpha_{d-1} r^{d-1} \nu(\vv)$.
%
\begin{figure}

\includegraphics{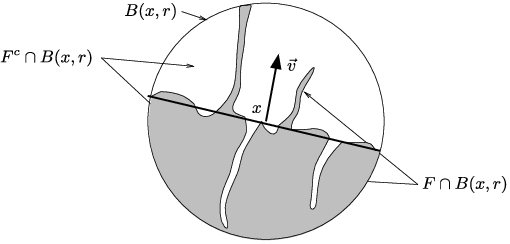}

\caption{The sets $F\cap B(x,r)$ and $B^-(x,r,\vv)$.}
\label{figtrompes}
\end{figure}
%
%
\end{rem}

Let $0<s<1$. We will prove that
\[
\liminf_{n\rightarrow\infty} \frac{V(\E_n)}{n^{d-1}} \geq\capa (F) (1-2s).
\]
We choose
\[
\eta= s \capa(F).
\]
We do not fix $\delta$ for the moment, and we consider the family of
balls $(B_i)_{i\in I\cup J \cup K}$ associated to $F$ by Lemma \ref
{chapitre7covering} (it depends on $\delta$) via the transformation we
did (thus $x_i$ and $\vv_i$ are rational for all $i$). By construction,
we get
%
%
\begin{equation}
\label{eqmin1} \biggl\llvert \capa(F) - \sum_{i \in I\cup J \cup K}
\alpha_{d-1} r_i ^{d-1} \nu(\vv_i)
\biggr\rrvert \leq s \capa(F).
\end{equation}
Since the capacities are nonnegative we have
%
%
\begin{equation}
\label{eqmin2} V(\E_n ) = \sum_{e\in\E_n}
t(e) \geq\sum_{i\in I \cup J \cup
K} V(\E_n \cap
B_i),
\end{equation}
where $\E_n \cap B_i = \{e\in\E_n | e\subset B_i\} $.

\emph{Step \textup{2:} From $V(\E_n \cap B_i)$ to $\alpha_{d-1}
r_i^{d-1} \nu(\vv_i)$}.

We define
\[
\eps= \eps(\delta) = \delta\min_{i \in I \cup J \cup K} \alpha_d
r_i ^d.
\]
Since $(R(\E_n))_{n\geq1}$ converges toward $F$, this implies that
%
%
\begin{equation}
\label{eqd} \exists n_0 \ \forall n\geq n_0\qquad \mathfrak{d}
\bigl(R(\E_n), F \bigr) \leq\eps.
\end{equation}
Let $(B(x,r),\vv) = (B_i (x_i, r_i ), \vv_i)$ for any $i \in I\cup J
\cup K$. Roughly speaking, we control the distance between $R(\E_n)$
and $F$ by (\ref{eqd}), and the distance between $F\cap B$ and $B^-$ by
construction of the balls. Thus we obtain a control on $\mathfrak d
(R(\E_n) \cap B, B^-)$, and since $R(\E_n) = r(\E_n) + [-1, 1]^d/(2n)$
[thus $ r(\E_n) =R(\E_n) \cap\ZZ^d_n$], this gives us a control on the
$\card(r(\E_n \cap B) \triangle(B^-\cap\ZZ^d_n))$.
%
\begin{figure}

\includegraphics{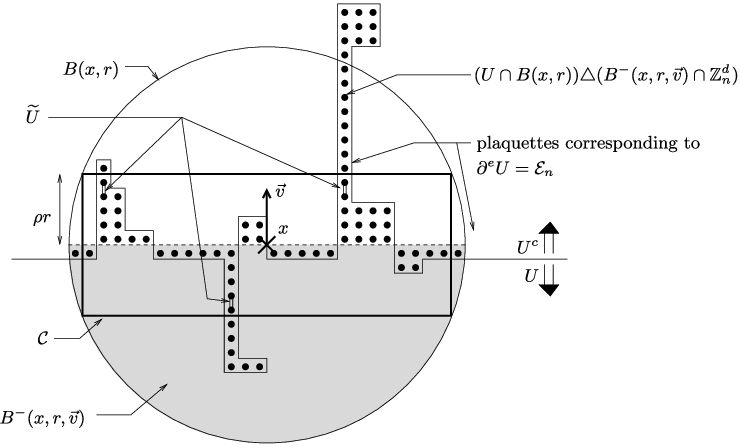}

\caption{The sets $U$ and $\widetilde U$.}
\label{figsetU}
\end{figure}
More precisely,
it is proved in Section~5 of \cite{CerfTheret09infc} that there exists
a set of vertices $U \subset\ZZ^d_n $ that satisfies
%
%
\begin{equation}
\label{eqlocal1} \card \bigl( \bigl(U \cap B(x, r) \bigr)\triangle \bigl(B^-
(x,r, \vv)\cap\ZZ^d_n \bigr) \bigr) \leq 4 \delta
\alpha_d r_i^d n^d
\end{equation}
and
%
%
\begin{equation}
\label{eqlocal2} \bigl(\partial^e U \bigr) \cap B =
\E_n \cap B,
\end{equation}
where we generalize the notation we have adopted for $\E_n \cap B_i$;
see Figure~\ref{figsetU}.
%
%
This statement is a bit more elaborated than expected because of the
slight difference between balls indexed in $I, J$ and $K$: we can
choose $U = r(\E_n)$ if $B=B_i$ for $i\in I \cup K$ and $U = \O_n
\setminus r(\E_n) $ if $B=B_i$ for $i\in J$. We define the
cylinder $\C\subset B (x,r)$ by
\[
\C= \cyl \bigl(\disc \bigl(x,r',\vv \bigr), \rho r \bigr),
\]
where $\rho$ is a positive constant we have to choose and $r' = r \cos
(\arcsin\rho)$. It is proved in Section~6 of \cite{CerfTheret09infc}
that there exists a set of edges $\widetilde U$ (denoted by $X^+ \cup
X^-$ in that paper) included in $B$ such that $( B \cap\partial^e U
) \cup
\widetilde U$ contains a cutset from the top to the bottom of $\C$ in
the direction $\vv$ (see Figure~\ref{figsetU}) and
\[
\card(\widetilde U) \leq C r^{d-1} \delta\rho^{-1}
n^{d-1},
\]
where $C= 10 d \alpha_d$ is a constant that depends only on the
dimension. We denote by $\phi_{\C}$ the maximal flow from the top to
the bottom of $\C$, that is, $\phi_{\C} = \phi(T(\disc(x,r',\vv),
\rho r), B(\disc(x,r',\vv), \rho r), \C)$. Thus, by the maximality
of $\phi_{\C}$ and thanks to equation (\ref{eqlocal2}),
%
%
\begin{eqnarray}
\label{eqlocal3} \phi_{\C} &\leq& V \bigl( \partial^e U
\cap B \bigr) + M C r^{d-1} \delta\rho^{-1} n^{d-1}
\nonumber
\\[-8pt]
\\[-8pt]
\nonumber
& =&
V( \E_n \cap B) + M C r^{d-1} \delta\rho^{-1}
n^{d-1}.
\end{eqnarray}
To complete the proof, it remains to compare $\phi_{\C}$ with $\alpha
_{d-1} r^{d-1} \nu(\vv) n^{d-1}$. This is done in Section~6 of \cite
{CerfTheret09infc} by almost covering $\disc(x, r', \vv)$ with a
finite family of disjoint closed hyperrectangles $(a_i)_{i \in\I}$
satisfying, for a constant $c=c(d)$ and chosen $\kappa>0$ as small as
we want,
\[
\sum_{i \in\I}\H^{d-1} (a_i) >
\alpha_{d-1} r'^{d-1} - \kappa \quad\mbox{and}\quad \sum
_{i \in\I} \H^{d-2} (\partial a_i) <
c r'^{d-2}.
\]
Thus the cylinders $(\cyl(a_i, \rho r))_{i\in\I}$ almost fill $\C$.
Since $x$ has rational coordinates and $\vv$ is a rational unit vector
(i.e., $\lambda\vv$ has rational coordinates for a positive real
number $\lambda$), we can choose the hyperrectangles $(a_i, i\in\I)$
with rational vertices. Indeed, we explained in Section~\ref{secconstr}
that there exist vectors $\vec u_i$, $i=2,\ldots,d$ that have integer
coordinates and such that $(\vv, \vec u_2,\ldots, \vec u_d)$ is an
orthogonal
basis of $\RR^d$. Then any hyperrectangle of the form $x + \sum_{i=2}^d
\lambda_i \vec u_i + \prod_{i=2}^d [0, \mu_i \vec u_i]$ with rational
$\lambda_i, \mu_i$ has rational vertices. We can choose the
hyperrectangles $(a_i, i\in\I)$ of this type; see Figure~\ref{figboule_rect}.
%
\begin{figure}

\includegraphics{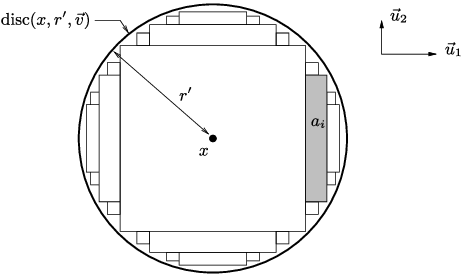}

\caption{The sets $\disc(x, r', \vv)$ and $(a_i, i\in\I)$.}
\label{figboule_rect}
\end{figure}

%

Let $h \in[\rho r, 2 \rho r] \cap\mathbb Q$. The cylinders $\cyl(a_i,
h)$, $i\in\I$, have rational vertices. If $\F_n$ is a cutset from the
top to the bottom of $\C$, then $\F_n \cap\cyl(a_i, \rho r)$ is a
cutset from the top to the bottom of $\cyl(a_i, h)$, and if we add to
$\F_n \cap\cyl(a_i, \rho r)$ some edges along the vertical sides of
$\cyl(a_i, h)$, we obtain a cutset in $\cyl(a_i, h)$ between the lower
half part and the upper half part of its boundary. More formally, if we define
\[
\P_i (n) = \cyl \bigl( \V_2 (\partial a_i,
2d/n) \cap\hyp(a_i), h \bigr),
\]
and if we denote by $P_i(n)$ the set of the edges included in $\P
_i(n)$, we get
\[
\sum_{i\in\I} \tau \bigl(\cyl(a_i, \rho
r), \vv \bigr) \leq\phi_{\C} + V \biggl(\bigcup
_{i\in\I} P_i(n) \biggr)
\]
for a complete proof, we refer to Section~6 of \cite{CerfTheret09infc}.
Moreover,
\[
\card \biggl(\bigcup_{i\in\I} P_i(n)
\biggr) \leq c' \rho r^{d-1} n^{d-1},
\]
where $c'$ is a constant depending on the dimension, thus
%
%
\begin{equation}
\label{eqlocal4} \phi_{\C} \geq\sum_{i\in\I}
\tau \bigl(\cyl(a_i, \rho r), \vv \bigr) - M c' \rho
r^{d-1} n^{d-1}.
\end{equation}
Combining inequalities (\ref{eqlocal3}) and (\ref{eqlocal4}), we get
%
%
\begin{equation}
\label{eqlocal5} V(\E_n \cap B ) \geq\sum
_{i\in\I} \tau \bigl(\cyl(a_i, h), \vv \bigr) - M
r^{d-1} c'' \bigl(\rho+ \delta
\rho^{-1} \bigr) n^{d-1}
\end{equation}
for a constant $c''$ depending on the dimension. Theorem \ref{thmnu}
states that for every cylinder $\cyl(A,h)$ with $A$ a nondegenerate
hyperrectangle normal to $\vv$ and $h>0$, we have
\[
\lim_{n\rightarrow\infty} \frac{\tau_n (\cyl(A,h), \vv)}{
n^{d-1} \H
^{d-1} (A)} = \nu(\vv) \qquad\mbox{a.s.}
\]
Thus these convergences hold a.s. simultaneously for all the rational
cylinders, that is, cylinders with rational vertices [like the
cylinders $(\cyl(a_i,h),\break  i \in\I)$]. We consider only realizations of
the capacities on which these convergences occur. Combined with
inequality (\ref{eqlocal5}), this implies that
\begin{eqnarray*}
&&\liminf_{n\rightarrow\infty} \frac{V(\E_n \cap B)}{n^{d-1}}\\
&&\qquad \geq \biggl( \sum
_{i\in\I} \H^{d-1} (a_i) \biggr) \nu(\vv) -
M r^{d-1} c'' \bigl(\rho+ \delta
\rho^{-1} \bigr)
\\
&&\qquad \geq \bigl(\alpha_{d-1} r'^{d-1} - \kappa
\bigr) \nu(\vv) - M r^{d-1} c'' \bigl(\rho+
\delta \rho^{-1} \bigr)
\\
&&\qquad \geq\alpha_{d-1} r^{d-1} \nu(\vv) - \alpha_{d-1}
\bigl(r^{d-1} - r'^{d-1} \bigr) \nu(\vv) - \kappa
\nu(\vv) \\
&&\qquad\quad{}- M r^{d-1} c'' \bigl(\rho+ \delta
\rho ^{-1} \bigr).
\end{eqnarray*}
We choose $\rho= \sqrt\delta$ and $\kappa= \alpha_{d-1} (r^{d-1} -
r'^{d-1})$. We define
\[
\nu_{\min} = \min_{\vv\in\SS^{d-1}} \nu( \vv) \quad\mbox{and}\quad
\nu_{\max} = \max_{\vv\in\SS^{d-1}} \nu(\vv).
\]
Under hypothesis \ref{hypo3}, we have $0 < \nu_{\min} \leq\nu_{\max
} <
+\infty$. Since $r' =\break   r \cos(\arcsin\rho)$, we get
%
%
\begin{eqnarray}
\label{eqmin3} &&\liminf_{n\rightarrow\infty} \frac{V(\E_n \cap B)}{n^{d-1}}\nonumber\\
&&\qquad\geq
\alpha_{d-1} r^{d-1} \nu(\vv)
\\
&&\qquad\quad{}- 2 \alpha_{d-1}
r^{d-1} \nu_{\max} \bigl[1 - \bigl(\cos(\arcsin\sqrt\delta)
\bigr)^{d-1} \bigr] - 2M c'' \sqrt\delta
r^{d-1}.\nonumber
\end{eqnarray}
Equations (\ref{eqmin1}), (\ref{eqmin2}) and (\ref{eqmin3}), and the
fact that $\sum_i \alpha_{d-1} r_i^{d-1} \leq\break \capa(F) (1 + s) \nu
_{\min} ^{-1}$, give
%
%
\begin{equation}
\label{eqminccl} \liminf_{n \rightarrow\infty} \frac{V (\E_n)}{n^{d-1}} \geq\capa
(F) (1- s - w ),
\end{equation}
where
%
%
\begin{equation}
w = 4 \frac{\nu_{\max}}{\nu_{\min}} \bigl[1- \bigl(\cos(\arcsin\sqrt \delta )
\bigr)^{d-1} \bigr] + \frac{4Mc'' }{\alpha_{d-1}} \sqrt\delta.
\end{equation}
For $\delta$ small enough, $w \leq s$, and we get
\[
\liminf_{n \rightarrow\infty} \frac{V (\E_n)}{n^{d-1}} \geq\capa (F) (1- 2s ).
\]
This completes the proof of Proposition \ref{propmin}.
\end{pf*}


\section{Continuous max-flow min-cut theorem}
\label{secccl}

\subsection
{Comparison between
continuous streams and cutsets}
\label{secinegmfmn}
\mbox{}
\begin{pf*}{Proof of Lemma \ref{lemineg}} Let $\vec\sigma$
be an
admissible continuous stream, that is:
\begin{itemize}
\item $\vec\sigma\in L^\infty \bigl(
\RR^d \rightarrow\RR^d, \L^d \bigr) \mbox{ and
} \vec\sigma=0\ \L^d \mbox {-a.e. on } \O^c,$
\item $\di\vec\sigma=0\ \L^d\mbox {-a.e. on } \O,$
\item $\vec\sigma\cdot\vv\leq\nu(\vv) \mbox{ for all } \vv\in\SS^{d-1}\
\L^d\mbox{-a.e. on }\O,$
\item $\vec\sigma\cdot\vv_{\O} \leq0\ \H^{d-1} \mbox {-a.e. on
} \G^1 \mbox{ and } \vec\sigma\cdot\vv_{\O} = 0\ \H
^{d-1} \mbox{-a.e. on } \G\setminus \bigl(\G^1\cup
\G^2 \bigr). $
\end{itemize}
As in Remark \ref{remdiv2} in Section~\ref{secdiv}, we obtain
\begin{eqnarray*}
&&\forall u \in\C_c ^{\infty} \bigl(\RR^d, \RR
\bigr)\\
 &&\qquad\int_{\RR^d} u \di\vec \sigma \,d\L^d  = -
\int_{\RR^d} \vec\sigma\cdot\vec\nabla u \,d \L^d \qquad
\mbox {by definition of } \di\vec\sigma
\\
&&\hspace*{70pt} \qquad =- \int_{\G} (\vec\sigma\cdot\vv_{\O}) u\, d
\H^{d-1} \qquad\mbox{ by definition of } \vec\sigma\cdot\vv_{\O},
\end{eqnarray*}
where for the last equality we have used the characterization of $\vec
\sigma\cdot\vv_{\O}$ given in equation (\ref{eqdefbord2}). Thus $
\di
\vec\sigma \,d\L^d $ is not only a distribution but also a measure, and
this measure is equal to $(\vec\sigma\cdot\vv_{\O}) \,d\H
^{d-1}\vert
_{\G}$.
Therefore we can apply to $\vec\sigma$ the techniques used in Section~\ref{secmax}. We consider a polyhedral set $P \subset\RR^d$ such that
\[
\overline{\G}^1 \subset\stackrel{\circ} {P},\qquad \overline{
\G}^2 \subset\,\stackrel{\circ} {\wideparen{\RR^d
\setminus P}} \quad\mbox{and}\quad \partial P \mbox{ is transversal to } \G.
\]
For any positive $\eta$, there exists $h_0 >0$ such that for all
$0<h<h_0$ we obtain inequality (\ref{eqmax12}),
\[
\Biggl\llvert \operatorname{flow}^{\mathrm{cont}}(\vec\sigma) -
\frac
{1}{2 h}\sum_{i=1}^{\N} \int
_{\cyl
(A_i,h)} \vec\sigma\cdot\vv_i \,d\L^d
\Biggr\rrvert \leq C \eta\| \vec\sigma\|_\infty.
\]
We have $\vec\sigma\cdot\vv\leq\nu(\vv)$ and $- \vec\sigma
\cdot\vv
\leq\nu(- \vv) = \nu(\vv)$ $\L^d$-a.e., thus $|\vec\sigma\cdot
\vv|
\leq\nu(\vv)$ $\L^d$-a.e. We obtain
\begin{eqnarray*}
\operatorname{flow}^{\mathrm{cont}}(\vec\sigma) & \leq&\frac
{1}{2h} \sum
_{i=1}^{\N} \int_{\cyl(A_i,
h)}
\vec\sigma\cdot\vv_i \,d\L^d + C \eta\| \vec\sigma\|
_{\infty}
\\
& \leq&\frac{1}{2h} \sum_{i=1}^{\N}
\nu(\vv_i) \L^d \bigl(\cyl(A_i, h ) \bigr) + C
\eta\| \vec\sigma\|_\infty
\\
& \leq&\frac{1}{2h} \sum_{i=1}^{\N}
\nu(\vv_i) 2h \H^{d-1} (A_i) + C \eta\| \vec
\sigma\|_\infty
\\
& \leq&\sum_{i=1}^{\N} \nu(
\vv_i) \H^{d-1} (A_i ) + C \eta\| \vec \sigma
\|_\infty
\\
& \leq&\int_{\partial P \cap\O'} \nu(\vv_P) \,d
\H^{d-1}+ C \eta \| \vec \sigma\|_\infty
\\
& \leq&\capa(P\cap\O) + \nu_{\max} \H^{d-1} \bigl(\partial P
\cap \bigl(\O' \setminus\O \bigr) \bigr) + C \eta\| \vec\sigma
\|_\infty
\\
& \leq&\capa(P\cap\O) + \nu_{\max} \widehat\eta+ C \eta\| \vec \sigma
\|_\infty,
\end{eqnarray*}
where we have used inequality (\ref{eqhateta}) to control $\H^{d-1}
(\partial
P \cap(\O' \setminus\O))$, and where $\nu_{\max} = \max_{\SS
^{d-1}} \nu$. Since $\eta$ and $\widehat\eta$ are arbitrarily small,
for all $P\subset\RR^d$ such that $\overline{\G}^1 \subset
\stackrel
{\circ}{P}$, $\overline{\G}^2 \subset\stackrel{\circ}{\wideparen
{\RR^d
\setminus P}}$ and $\partial P$ is transversal to $\G$, we obtain
%
%
\begin{equation}
\label{eqineg1} \operatorname{flow}^{\mathrm{cont}}(\vec\sigma) \leq\capa(P\cap
\O).
\end{equation}
The following result has been proved in \cite{CerfTheret09geoc}:
%
\begin{thmm}[(Theorem 11 in \cite{CerfTheret09geoc})]
\label{thmapproxpol}
We suppose that hypotheses \ref{hypo1} are fulfilled, and that the law
of the capacities is integrable,
\[
\int_{\RR^+} x \,d\Lambda(x) < +\infty.
\]
Let $F$ be a subset of $\O$ having finite perimeter in $\O$. For any
$\eps>0$, there exists a polyhedral set $P$ whose boundary $\partial
P$ is
transversal to $\G$ and such that
\begin{eqnarray*}
&\displaystyle \overline{\G}^1 \subset\stackrel{\circ} {P},\qquad \overline{
\G}^2 \subset\stackrel{\circ} {\wideparen{\RR^d
\setminus P}}, \qquad\L ^d \bigl(F \triangle(P\cap\O) \bigr) <\eps,&
\\
&\displaystyle\int_{\partial^* P \cap\O} \nu \bigl( \vv_P (x) \bigr) \,d
\H^{d-1} (x) = \capa(P\cap \O) \leq\capa(F) + \eps.&
\end{eqnarray*}
\end{thmm}
Combining inequality (\ref{eqineg1}) and Theorem \ref{thmapproxpol}, we
obtain that for all $F\subset\O$ such that $\mathbh{1}_F \in \operatorname{BV}(\O)$,
\[
\operatorname{flow}^{\mathrm{cont}}(\vec\sigma) \leq\capa(F),
\]
and thus Lemma \ref{lemineg} is proved.
\end{pf*}


\subsection{\texorpdfstring{End of the proofs of Theorems \protect\ref{thmpcpal}, \protect\ref{thmpcpal2}, \protect\ref{cormaxflowmincut} and \protect\ref{corLGNflow}}
{End of the proofs of Theorems 1.1, 1.2, 1.3 and 1.4}}
\label{secccl2}

We suppose first that hypothesis \ref{hypo3} is fulfilled. Let $(\vec
\mu
_n^{\max})_{n \geq1}$ be a sequence of discrete maximal streams and
$(\E
_n^{\min})_{n\geq1}$ be a sequence of discrete minimal cutsets. From a
subsequence of $(\vec\mu_n^{\max})_{n \geq1}$ that converges weakly
toward a measure $\vec\mu= \vec\sigma\L^d$, we can a.s. extract a
sub-subsequence\vspace*{1pt} $(\vec\mu_{\varphi(n)}^{\max})_{n\geq1}$ such that
$(R(\E_{\varphi(n)}^{\min}) )_{n\geq1}$ converges also for the
distance $\mathfrak d$ toward a set $F\subset\O$ of finite perimeter.
Conversely, from a subsequence of $(R(\E_n^{\min}))_{n\geq1}$ that
converges for the distance $\mathfrak d$ to a set $F\subset\O$ of
finite perimeter, we can extract a sub-subsequence $(R(\E_{\varphi
(n)}^{\min}) )_{n\geq1}$ such that $(\vec\mu_{\varphi(n)}^{\max
})_{n\geq1}$ converges weakly toward a measure $\vec\mu= \vec
\sigma\L
^d$. Combining Propositions~\ref{propmax2} and~\ref{propmin}, we obtain
that a.s.
%
%
\begin{eqnarray}
\label{eqccl1} \capa(F) &\leq&\liminf_{n\rightarrow\infty} \frac{V(\E^{\min
}_{\varphi(n)})}{\varphi(n)^{d-1}} =
\lim_{n\rightarrow\infty} \frac{\phi_{\varphi(n)}}{\varphi(n)^{d-1}}
\nonumber
\\[-8pt]
\\[-8pt]
\nonumber
 &=& \lim_{n\rightarrow
\infty}
\operatorname{flow}^{\mathrm{disc}}_{\varphi(n)} \bigl(\vec\mu
_{\varphi(n)}^{\max} \bigr) = \operatorname{flow}^{\mathrm{cont}}(\vec
\sigma ).
\end{eqnarray}
Since Lemma \ref{lemineg} implies that
\[
\capa(F) \geq\phi_{\O}^a \geq\phi_{\O}^b
\geq\operatorname {flow}^{\mathrm{cont}}(\vec \sigma),
\]
combining inequality (\ref{eqccl1}) and Lemma \ref{lemineg}, we obtain
that a.s.
\[
\capa(F) = \phi_{\O}^a = \phi_{\O}
^b = \operatorname {flow}^{\mathrm{cont}}(\vec\sigma).
\]
If hypothesis \ref{hypo3} is not fulfilled, then $\nu=0$, $\phi_{\O}^a
=0$, and Lemma \ref{lemineg} implies that $\phi_{\O}^b = \phi_{\O}^a
=0$, thus all the admissible continuous streams in $\Sigma^b$ are null
and all the admissible continuous cutsets have null capacity and are in
$\Sigma^a$. This completes the proofs of Theorems \ref{thmpcpal},
\ref
{thmpcpal2} and \ref{cormaxflowmincut}.

Let us prove Theorem \ref{corLGNflow}. We notice that if a subsequence
of $(\phi_n / n^{d-1})_{n\geq1}$ converges toward a real variable
$\phi
$, then we can extract a sub-subsequence $(\phi_{\varphi(n)} /
\varphi
(n)^{d-1})_{n\geq1}$ along which the maximal flows converge toward
$\phi$ and the maximal streams $(\vec\mu_{\varphi(n)}^{\max
})_{n\geq
1}$ converge toward a continuous stream $\vec\sigma$.\vspace*{1.5pt} Then by
Proposition \ref{propmax2} and Theorem \ref{thmpcpal} we know that a.s.
%
%
\begin{equation}
\label{eqccl2} \phi= \lim_{n\rightarrow\infty} \frac{\phi_{\varphi
(n)}}{\varphi(n)^{d-1}} =
\operatorname{flow}^{\mathrm{cont}}(\vec \sigma) = \phi_{\O}^b
.
\end{equation}
We claim that $(\phi_n / n^{d-1})_{n\geq1}$ takes its values in a
deterministic compact of $\RR^d$---this, together with equation (\ref
{eqccl2}), completes the proof of Theorem \ref{corLGNflow}. Indeed, let
$P$ be a polyhedral set of $\RR^d$ such that
\[
\overline{\G}^1 \subset\stackrel{\circ} {P} \quad\mbox{and}\quad \overline{
\G}^2 \subset\stackrel{\circ} {\wideparen{\RR^d
\setminus P}}.
\]
We define
\[
\F_n = \{ e \in\Pi_n | e\cap\partial P \neq\varnothing
\}.
\]
At least for $n$ large enough, $\F_n$ separates $\G^1_n$ from $\G^2_n$
in $\O_n$, thus
\[
\frac{\phi_n}{n^{d-1}} \leq V(\F_n) \leq M \frac{\card(\F
_n)}{n^{d-1}},
\]
and since the $(d-1)$ dimensional Minkowski content of $\partial P
\cap\O$
exists and is equal to $\H^{d-1}(\partial P \cap\O)$, for $n$ large enough
we get
\begin{eqnarray*}
\card(\F_n) &\leq&2d \card \bigl(\V_{\infty} (\partial P \cap
\O, 2/n) \cap\ZZ ^d_n \bigr) \leq2d n^d
\L^d \bigl(\V_\infty(\partial P \cap\O, 3/n) \bigr) \\
&\leq&
n^{d-1} K \H^{d-1}(\partial P \cap\O)
\end{eqnarray*}
for a constant $K$ that depends on $d$.

\section*{Acknowledgments} The second author would like to warmly
thank Antoine Lemenant and Vincent Millot for helpful discussions.


%
%



\printaddresses

\end{document}